\documentclass[12pt,leqno]{article}
\usepackage{amssymb,amsmath,amsfonts,amsthm}
\usepackage{relsize}
\usepackage[sort]{cite}
\usepackage{esvect}
\usepackage{txfonts}
\usepackage{graphicx}
\usepackage[usenames]{color}

\makeatletter
\renewcommand{\fnum@figure}{Fig. \thefigure}
\makeatother
\newtheorem{theorem}{Theorem}[section]
\newtheorem{definition}[theorem]{Definition}

\newtheorem{proposition}[theorem]{Proposition}
\newtheorem{remark}[theorem]{Remark}
\newtheorem{lemma}[theorem]{Lemma}

\newtheorem{problem}[theorem]{Problem}

\newtheorem{conjecture}[theorem]{Conjecture}

\newtheorem{claim}[theorem]{Claim}

\newcommand {\Hc}      {{\mathcal H}}

\newcommand {\Kc}      {{\mathcal K}}

\newcommand {\Mc}      {{\mathcal M}}

\newcommand {\Vc}      {{\mathcal V}}
\newcommand {\R}       {{\bf R}}

\newcommand {\RN}      {\R^n}

\newcommand {\RT}      {\R^2}

\newcommand {\tA}      {\widetilde{A}}

\newcommand {\tC}      {\widetilde{C}}

\newcommand {\tG}      {\widetilde{G}}

\newcommand {\tK}      {\widetilde{K}}

\newcommand {\tS}      {\widetilde{S}}
\newcommand {\tT}      {\widetilde{T}}

\newcommand {\tMc}     {\widetilde{{\mathcal M}}}

\newcommand {\tKc}     {\widetilde{\Kc}}

\newcommand {\tf}      {\tilde{f}}
\newcommand {\tg}      {\tilde{g}}

\newcommand {\tr}      {\tilde{r}}

\newcommand {\tu}      {\tilde{u}}
\newcommand {\tv}      {\tilde{v}}

\newcommand {\tgm}     {\tilde{\gamma}}
\newcommand {\mfM}      {{\mathfrak M}}
\newcommand {\mfC}      {{\mathfrak C}}

\newcommand {\mfR}      {{\mathfrak R}}
\newcommand {\ve}      {\varepsilon}

\newcommand {\MR}      {(\Mc,\rho)}
\newcommand {\MS}      {\mfM}
\newcommand {\vl}      {\vec{\lambda}}
\newcommand {\vf}      {\varphi}
\newcommand {\emp}     {\emptyset}

\newcommand {\CX}      {\mfC(\X)}

\newcommand {\RCT}     {\mfR(\RT)}

\newcommand {\HR}      {\Hc}

\newcommand {\slbig}   {\mathlarger{\mathlarger{/}}}

\newcommand {\cbig}    {\mathlarger{\mathlarger{\cap}}}
\newcommand {\cbg}    {\,\mathlarger{\mathlarger{\cap}}\,}

\newcommand {\cupbig}    {\mathlarger{\mathlarger{\cup}}}

\newcommand {\BX}      {B_{\X}}
\newcommand {\BXR}     {I_0}

\newcommand {\LTI}     {\ell^2_\infty}

\newcommand {\X}       {X}

\newcommand {\KM}      {\Kc_m(\X)}
\newcommand {\KX}      {\Kc(\X)}

\newcommand {\BAL}[3]  {{\cal BR\,}[{#1}\!:\!{#2};{#3}]}

\newcommand {\Lip}     {\operatorname{Lip}}

\newcommand {\ST}      {\operatorname{St}}

\newcommand {\dhf}     {\operatorname{d_H}}
\newcommand {\ds}      {\operatorname{d}}

\newcommand {\diam}    {\operatorname{diam}}
\newcommand {\dist}    {\operatorname{dist}}

\newcommand {\smsk}    {\smallskip}
\newcommand {\msk}     {\medskip}
\newcommand {\bsk}     {\bigskip}
\newcommand {\bx}      {\hspace{10mm}$\blacksquare$}
\newcommand {\rbx}     {\hspace{10mm}$\vartriangleleft$}

\newcommand {\nn}      {\nonumber}
\newcommand {\rf}[1]    {(\ref{#1})}      
\newcommand {\reff}[1] {\ref{#1}}         
\newcommand{\lbl}[1]      {\label{#1}}       
\newcommand{\be}          {\begin{eqnarray}}
\newcommand{\bel}[1]      {\begin{eqnarray} \label{#1}}
\newcommand{\ee}           {\end{eqnarray}}
\newcommand {\SECT}[2] {\section*{\centerline{\normalsize
{\bf #1}}} \setcounter{section}{#2}
\setcounter{theorem}{0}\setcounter{equation}{0}}

\begin{document}
\parindent 1em
\parskip 0mm
\medskip
\centerline{\large{\bf On the Core of a Low Dimensional  Set-Valued Mapping}}
\vspace*{12mm}
\centerline{By~ {\sc Pavel Shvartsman}}\vspace*{3 mm}
\centerline {\it Department of Mathematics, Technion - Israel Institute of Technology,}\vspace*{1 mm}
\centerline{\it 32000 Haifa, Israel}\vspace*{1 mm}
\centerline{\it e-mail: pshv@technion.ac.il}
\bsk\bsk
\renewcommand{\thefootnote}{ }
\footnotetext[1]{{\it\hspace{-6mm}Math Subject
Classification:} 46E35\smallskip
\par This research was supported by Grant No 2014055 from the United States-Israel Binational Science Foundation (BSF).}

%
%

\begin{abstract} Let $\MS=(\Mc,\rho)$ be a metric space and let $\X$ be a Banach space. Let $F$ be a set-valued mapping from $\Mc$ into the family $\Kc_m(\X)$ of all compact convex subsets of $\X$ of dimension at most $m$. The main result in our recent joint paper with Charles Fefferman (which is referred to as a ``Finiteness Principle for Lipschitz selections'') provides efficient conditions for the existence of a Lipschitz selection of $F$, i.e., a Lipschitz mapping $f:\Mc\to\X$ such that $f(x)\in F(x)$ for every $x\in\Mc$.
\par We give new alternative proofs of this result in two
special cases. When $m=2$ we prove it for $X=\mathbf{R}^{2}$, and when $m=1$ we prove it for all choices of $X$. Both of these proofs make use of a simple reiteration formula for the ``core'' of a set-valued mapping $F$, i.e., for a mapping $G:\Mc\to\Kc_m(\X)$ which is Lipschitz with respect to the Hausdorff distance, and such that $G(x)\subset F(x)$ for all $x\in\Mc$.
\end{abstract}

\SECT{1. Introduction.}{1}

\indent
\par Let $\mfM=(\Mc,\rho)$ be a {\it pseudometric space},
i.e., suppose that the ``distance function'' \mbox{$\rho:\Mc\times\Mc\to [0,+\infty]$} satisfies $\rho(x,x)=0$, $\rho(x,y)=\rho(y,x)$, and $\rho(x,y)\le \rho(x,z)+\rho(z,y)$ for all $x,y,z\in\Mc$. Note that $\rho(x,y)=0$ may hold with $x\ne y$, and $\rho(x,y)$ may be $+\infty$.
\par Let  $(\X,\|\cdot\|)$ be a real Banach space. Given a non-negative integer $m$ we let $\KM$ denote the family of all {\it non-empty compact convex subsets} $K\subset \X$ of dimension at most $m$. (We say that a convex subset of $\X$ has dimension at most $m$ if it is contained in an affine subspace of $\X$ of dimension at most $m$.) We let $\Kc(\X)=\bigcup\{\KM: m=0,1,...\}$ denote the family of all non-empty compact convex finite-dimensional subsets of $\X$.
\par By $\Lip(\Mc,\X)$ we denote the space of all Lipschitz mappings from $\Mc$ to $\X$ equipped with the Lipschitz seminorm
$$
\|f\|_{\Lip(\Mc,\X)}=\inf\,\{\,\lambda>0:\|f(x)-f(y)\|
\le\lambda\,\rho(x,y)~~\text{for all}~~x,y\in\Mc\,\}.
$$
\par In this paper we study the following problem.
\begin{problem} {\em Suppose that we are given a set-valued mapping $F$ which to each point $x\in\Mc$ assigns a set $F(x)\in\KM$. A {\it selection} of $F$ is a map $f:\Mc\to \X$ such that $f(x)\in F(x)$ for all $x\in\Mc$.
\par We want to know {\it whether there exists a selection $f$ of $F$ in the space $\Lip(\Mc,\X)$}. Such an $f$ is called a {\it Lipschitz selection} of the set-valued mapping $F:\Mc\to\KM$.
\par If a Lipschitz selection $f$ exists, then we ask {\it how small we can take its Lipschitz seminorm}. See Fig. 1.}
\end{problem}

\begin{figure}[h!]
\hspace{10mm}
\includegraphics[scale=0.6]{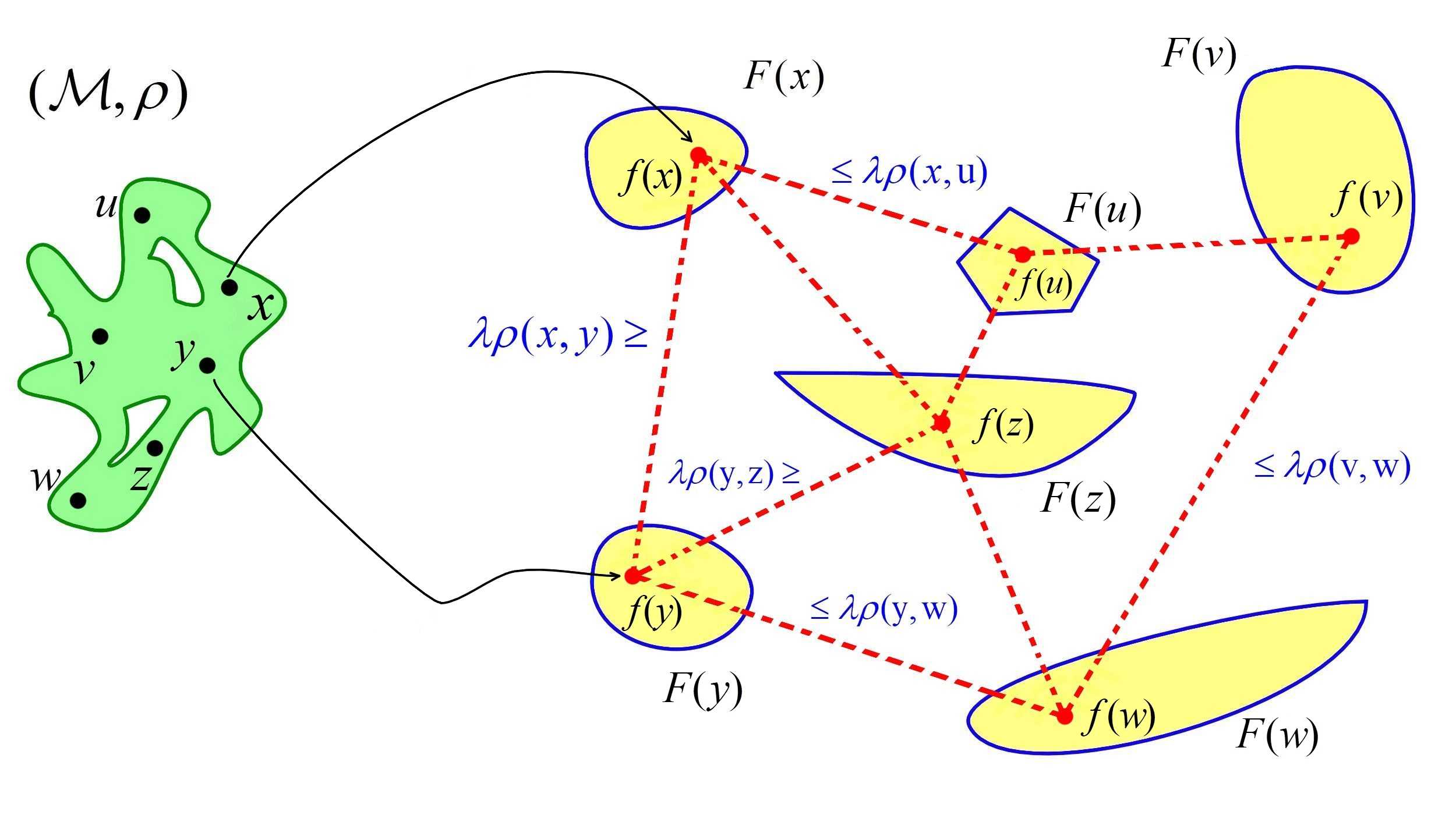}
\caption{$f:\Mc\to\RT$ is a Lipschitz selection of the set-valued mapping $F:\Mc\to\Kc(\RT)$.}
\end{figure}

\par The following result provides efficient conditions for the existence of a Lipschitz selection of an arbitrary set-valued mapping from a pseudometric space into the family $\KM$. We refer to it as a ``Finiteness Principle for Lipschitz selections'', or simply as a ``Finiteness Principle''.
\begin{theorem}\label{MAIN-FP} (Fefferman,Shvartsman \cite{FS-2018}) Fix $m\ge 1$. Let $(\Mc,\rho)$ be a pseudometric space, and let $F:\Mc\to\KM$ for a Banach space $\X$. Let
\bel{NMY-1}
N(m,\X)=2^{\ell(m,X)}~~~~~\text{where}~~~~~ \ell(m,X)=\min\{m+1,\dim\X\}.
\ee
\par Suppose that for every subset $\Mc'\subset\Mc$ consisting of at most $N=N(m,\X)$ points, the restriction $F|_{\Mc'}$ of $F$ to $\Mc'$ has a Lipschitz selection $f_{\Mc'}$ with Lipschitz  seminorm $\|f_{\Mc'}\|_{\Lip(\Mc',\X)}\le 1$.

\par Then $F$ has a Lipschitz selection $f$ with Lipschitz  seminorm
$\|f\|_{\Lip(\Mc,\X)}\le \gamma$ where $\gamma=\gamma(m)$ is a positive constant depending only $m$.
\end{theorem}

\par There is an extensive literature devoted to various
versions of Finiteness Principles for Lipschitz selections and related topics. We refer the reader to the papers \cite{Ar,AF,BL, FIL-2017,FP-2019,FS-2018,S-2001,PR,PY1,PY2,
S-2002,S-2004,S-2008} and references therein for numerous results in this direction.
\par We note that the ``finiteness number'' $N(m,\X)$ in Theorem \reff{MAIN-FP} is optimal; see \cite{S-2002}.
\par For the case of the trivial distance function $\rho\equiv 0$, Theorem \reff{MAIN-FP} agrees with the classical Helly's Theorem \cite{DGK}, except that the optimal finiteness constant for $\rho\equiv 0$ is
$$
n(m,X)=\ell(m,X)+1=\min\{m+2,\dim \X+1\}~~~~\text{in place of}~~~~N(m,X)=2^{\ell(m,X)}.
$$
Thus, Theorem \reff{MAIN-FP} may be regarded as a generalization of Helly's Theorem.
\par Our interest in Helly-type criteria for the existence of Lipschitz selections was initially motivated by some intriguing close connections of this problem with the
classical Whitney extension problem \cite{W1}, namely, the problem of characterizing those functions defined on a closed subset, say $E\subset\RN$, which are the restrictions to $E$ of $C^m$-smooth functions on $\RN$.
We refer the reader to the papers \cite{BS-1994,BS-1997,BS-2001,F-2005,F-2006,F-2009,
FI-2020,S-2008} and references therein for numerous results and techniques concerning this topic.
\par One of the main ingredients of the proof of Theorem \reff{MAIN-FP} is the construction of a special set-valued mapping $G:\Mc\to\KM$ introduced in \cite{FS-2018} which we call a {\it ``core''} of the set-valued mapping $F$. In fact each core is associated with a positive constant. Here are the relevant definitions.
\begin{definition}\label{D-CORE} {\em  Let $\gamma$ be a positive constant, and let $F:\Mc\to\KM$ be a set-valued mapping. A set-valued mapping $G:\Mc\to\KM$ is said to be a {\it $\gamma$-core} of $F$ if
\smsk
\par (i). $G(x)\subset F(x)$ for all $x\in\Mc$;
\smsk
\par (ii). $G$ is $\gamma$-Lipschitz with respect to  Hausdorff distance, i.e.,
$$
\dhf(G(x),G(y))\le \gamma\,\rho(x,y)~~~~~~\text{for all} ~~~x,y\in\Mc.
$$
}
\end{definition}
\msk
\par We refer to a map $G$ as a {\it core of $F$} if $G$ is a $\gamma$-core of $F$ for some $\gamma>0$.
\smsk
\par See Fig. 2, 3, 4.
\bsk

\begin{figure}[h!]
\hspace{40mm}
\includegraphics[scale=1.15]{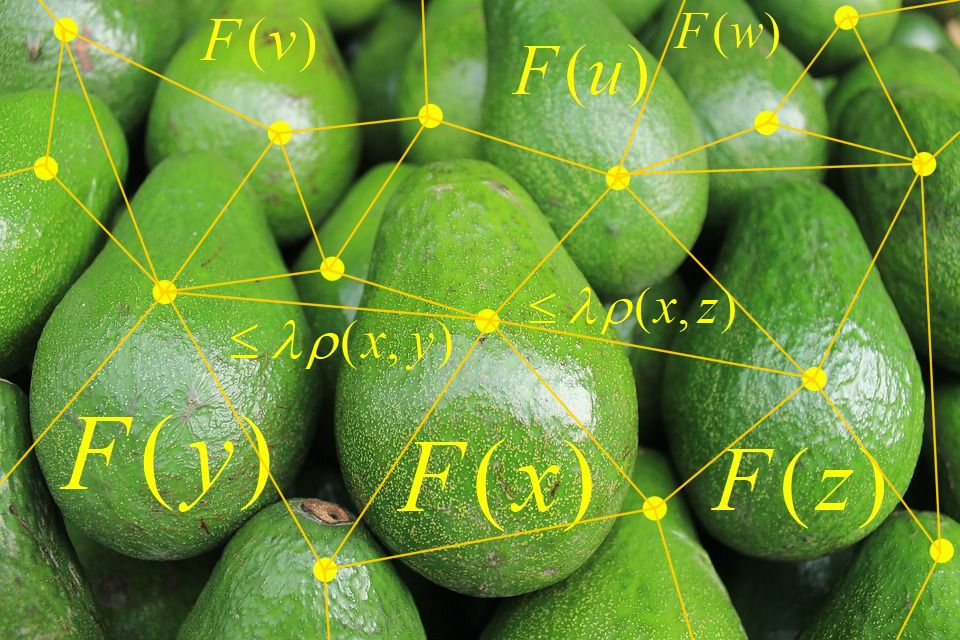}
\caption{A set-valued mapping $F$ into a family of avocados
and its Lipschitz selection
\\
\hspace*{11.5mm}
with Lipschitz seminorm at most $\lambda$.}
\end{figure}
\begin{figure}[h!]
\hspace{40mm}
\includegraphics[scale=0.97]{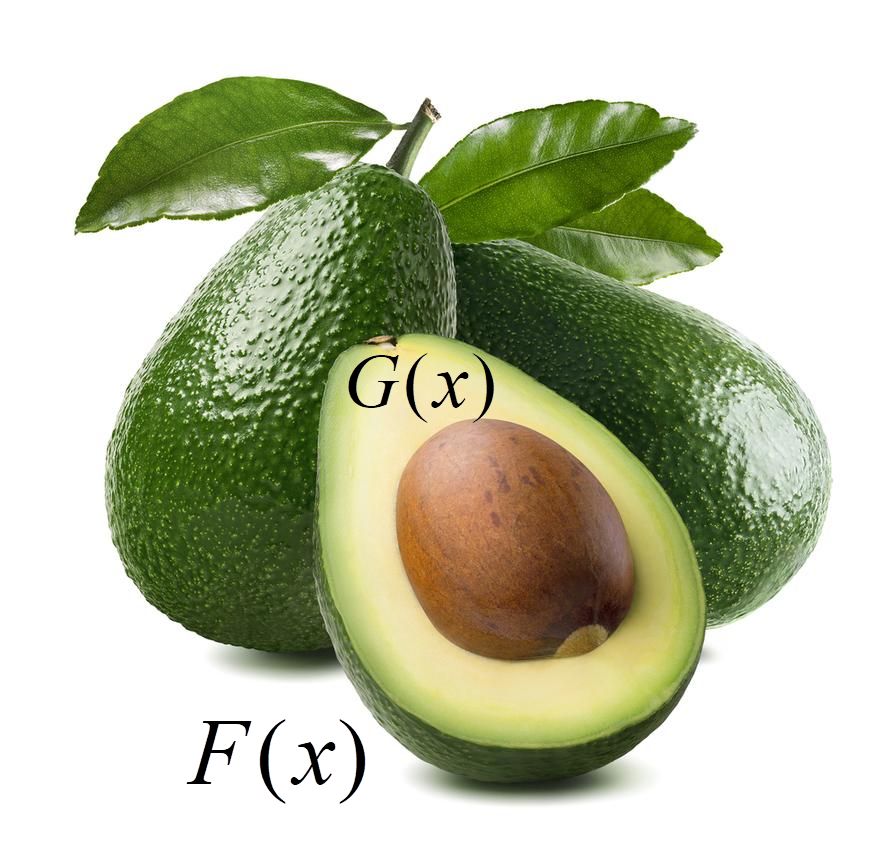}
\vspace*{-2mm}
\caption{The core $G(x)$ is a convex closed subset of $F(x)$.}
\end{figure}
\begin{figure}[h!]
\hspace{40mm}
\includegraphics[scale=0.31]{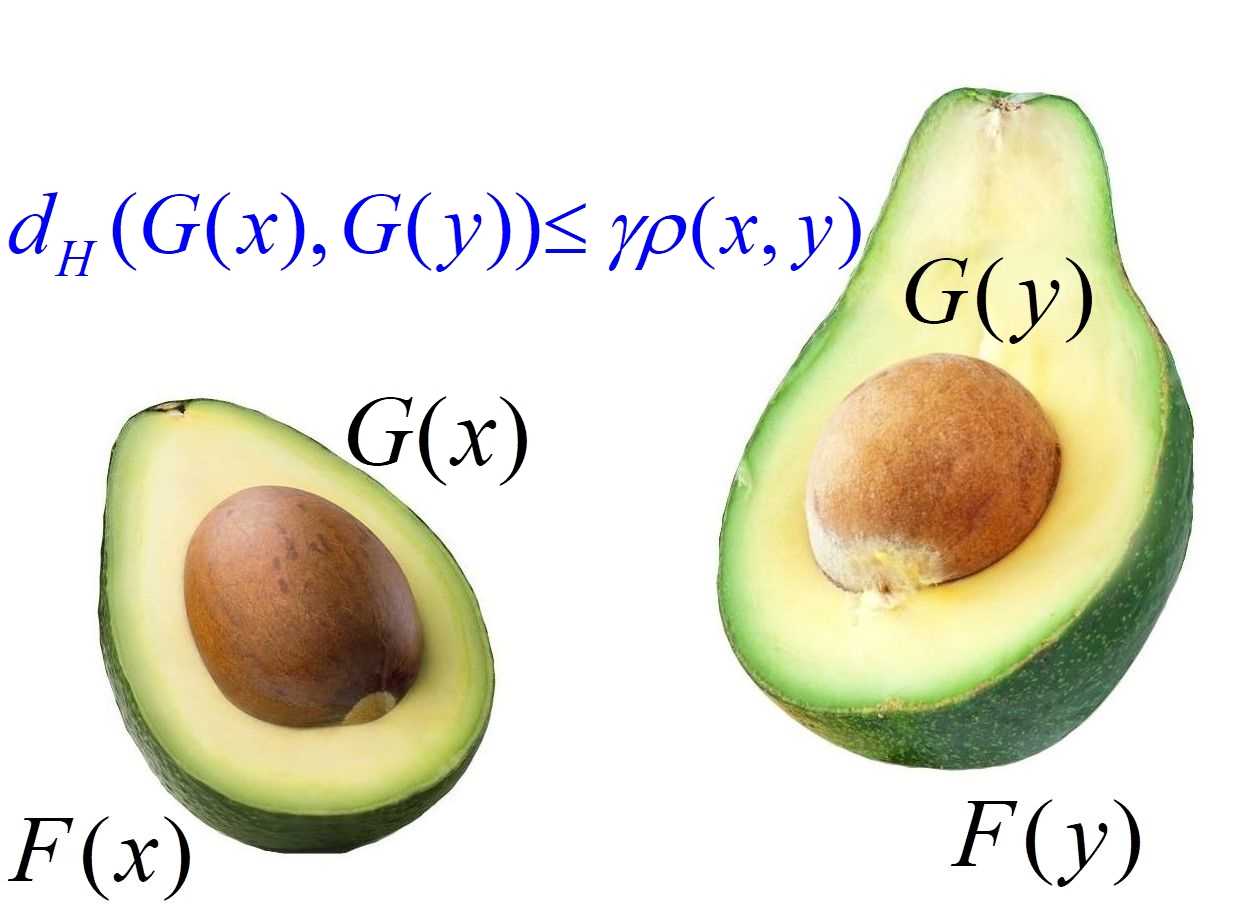}
\caption{The $\gamma$-core $G$ is $\gamma$-Lipschitz with respect to the Hausdorff distance.}
\end{figure}
\bsk\bsk

\par Recall that the Hausdorff distance $\dhf(A,B)$ between two non-empty bounded subsets $A$ and $B$ of $X$ is defined by
\begin{align}\lbl{HFD}
\dhf(A,B)=\inf\,\{r>0:~A+\BX(0,r)\supset B~~\text{and}~~
B+\BX(0,r)\supset A\}.
\end{align}
\par Here and throughout this paper, for
each $x\in X$ and $r>0$, we use the standard notation $\BX(x,r)$ for the closed ball in $X$ with center $x$ and radius $r$. We also let $\BX=\BX(0,1)$ denote the unit ball in $\X$, and we write $r\BX$ to denote the ball $\BX(0,r)$.

\msk
\par In Definition \reff{D-CORE} $m$ can be any non-negative integer not exceeding the dimension of the Banach space $X$. It can happen that a core $G:\mathcal{M}\to\mathcal{K}_{m}(X)$ of a given
set-valued mapping $F:\mathcal{M}\to\mathcal{K}_{m}(X)$ in fact maps $\mathcal{M}$ into the smaller collection $\mathcal{K}_{m'}(X)$ for some integer $m'\in[0,m$). The next claim shows that the existence of some core $G:\mathcal{M}\to\mathcal{K}_{m}(X)$
for $F$ implies the existence of a (possibly different) core which maps $\mathcal{M}$ into $\mathcal{K}_{0}(X)$. Since $\mathcal{K}_{0}(X)$ is identified with $X$, that core is simply a Lipschitz selection of $F$.
\bsk

\begin{claim} (\cite[Section 5]{FS-2018}) Let $\gamma$ be a positive constant, let $m$ be a non-negative
integer, and let $G:\Mc\to\KM$ be a $\gamma$-core of a set-valued mapping $F:\Mc\to\KM$ for some Banach space $X$. Then $F$ has a Lipschitz selection $f:\Mc\to \X$ with $\|f\|_{\Lip(\Mc,\X)}\le C\,\gamma$ where $C=C(m)$ is a constant depending only on $m$.
\end{claim}
\par In \cite{FS-2018} we showed that this claim follows  from Definition \reff{D-CORE} and the existence of the so-called {\it ``Steiner-type point''} map $\ST:\KM\to \X$ \cite{S-2004}. See formula \rf{FST-C}.
\msk
\par In \cite{FS-2018} given a set-valued mapping $F:\Mc\to\KM$ satisfying the hypothesis of Theorem \reff{MAIN-FP}, we constructed a $\gamma$-core $G$ of $F$ with a positive constant $\gamma$ depending only on $m$. We produced the core $G$ using a rather delicate and complicated procedure whose main ingredients are families  of {\it Basic Convex Sets} associated with $F$, metric spaces with bounded {\it Nagata dimension}, ideas and methods of work \cite{FIL-2017} related to the case $\Mc=\RN$, and Lipschitz selections on finite metric trees. See \cite{FS-2018} for more details.
\msk
\par In the present paper we suggest and discuss a different new geometrical method for pro\-du\-cing a core of a set-valued mapping. Its main ingredient is the so-called {\it balanced refinement} of a set-valued mapping which we define as follows.
\begin{definition} {\em  Let $\lambda\ge 0$, let $\MR$ be a pseudometric space, let $X$ be a Banach space, and let $F:\Mc\to\KM$ be a set-valued mapping for some non-negative integer $m$. For each $x\in\Mc$ we consider the subset of $F(x)$ defined by
$$
\BAL{F}{\lambda}{\rho}(x)=
\bigcap_{z\in\Mc}\,
\left[F(z)+\lambda\,\rho(x,z)\,\BX\right].
$$
\par We refer to the set-valued mapping
$\BAL{F}{\lambda}{\rho}:\Mc\to\KM\cup\{\emp\}$ as
the {\it $\lambda$-balanced refinement} of the mapping $F$.}
\end{definition}
\par We note that any Lipschitz selection $f$ of a set-valued mapping $F:\Mc\to\KM$ with $\|f\|_{\Lip(\Mc,\X)}\le \lambda$ is also a Lipschitz selection of the $\lambda$-balanced refinement of $F$, i.e.,
$$
f(x)\in\BAL{F}{\lambda}{\rho}(x)~~~~\text{for all}~~~~ x\in\Mc.
$$
\par Various geometrical parameters of the set $\BAL{F}{\lambda}{\rho}(x)$ (such as diameter and width, etc.) may turn out to be smaller than the same parameters for the set $F(x)$ which contains it. When attempting
to find Lipschitz selections of $F$ it may turn out to be convenient for our purposes to search for them in the more ``concentrated'' setting provided by the sets $\BAL{F}{\lambda}{\rho}(x)$. One can take this approach still further by searching in even smaller sets which can be obtained from {\it consecutive iterations of balanced refinements} of $F$, i.e., the sets which we describe in the following definition.
\begin{definition}\label{F-IT} {\em  Let $\ell$ be a positive integer, and let
$\vec{\lambda}=\left\{\lambda_k:1\le k\le\ell\right\}$ be a finite sequence of $\ell$ non-negative numbers $\lambda_k$. We set $F^{[0]}=F$, and, for every $x\in\Mc$ and integer $k\in[0,\ell-1]$, we define
\begin{align}\lbl{IT-F}
F^{[k+1]}(x)=\BAL{F^{[k]}}{\lambda_{k+1}}{\rho}(x)=
\bigcap_{z\in \Mc}\,
\left[F^{[k]}(z)+\lambda_{k+1}\,\rho(x,z)\,\BX\right].
\end{align}
\par We refer to the set-valued mapping $F^{[k]}:\Mc\to\KM\,\cupbig\,\{\emp\}$, $k\in[1,\ell]$, as
{\it the $k$-th order $(\vl,\rho)$-balanced refinement of $F$}.}
\end{definition}
\smsk
\par Clearly,
\begin{align}\lbl{CH-N}
F^{[k+1]}(x)\subset F^{[k]}(x)~~~~\text{on}~~~\Mc~~~ \text{for every}~~~k\in[0,\ell-1].
\end{align}
(Put $z=x$ in the right hand side of \rf{IT-F}.)

\begin{remark} {\em Of course, for each integer $k\in[1,\ell]$ the set $F^{[k]}(x)$ also depends on the sequence $\vec{\lambda}=\{\lambda_k:1\le k\le\ell\}$, on the pseudometric space $\mfM=\MR$ and the Banach space $\X$. However, in all places where we use $F^{[k]}$'s,  these objects, i.e., $\vec{\lambda}$, $\mfM$ and $\X$, are clear from the context. Therefore, in these cases, we omit any mention of $\vec{\lambda}$, $\mfM$ and $\X$ in the notation of $F^{[k]}$'s.\rbx}
\end{remark}
\smsk
\par We formulate the following
\begin{conjecture}\label{BR-IT} Let $\MR$ be a pseudometric space, and let $\X$ be a Banach space. Let $m$ be a fixed positive integer and (as in the formula \rf{NMY-1}
of Theorem \reff{MAIN-FP}) let $N(m,\X)$ denote the ``finiteness number'' $N(m,\X)=2^{\ell}$ where $\ell=\ell(m,X)=\min\{m+1,\dim\X\}$.
\par There exist a constant $\gamma\ge 1$ and a sequence
$\vec{\lambda}=\{\lambda_k:1\le k\le\ell\}$ of $\ell$ numbers $\lambda_{k}$ all satisfying $\lambda_{k}\ge1$ such that the following holds:
\smsk
\par Let $F:\Mc\to\KM$ be a set-valued mapping such that, for every $\Mc'\subset\Mc$ with $\#\Mc'\le N(m,\X)$, the restriction $F|_{\Mc'}$ of $F$ to $\Mc'$ has a Lipschitz selection $f_{\Mc'}:\Mc'\to \X$ with Lipschitz seminorm $\|f_{\Mc'}\|_{\Lip(\Mc',\X)}\le 1$.
\par Then the $\ell$-th order balanced refinement
of the mapping $F$, namely the set-valued mapping $F^{[\ell]}:\Mc\to\KM$ is a $\gamma$-core of $F$.
\smsk
\par Here $F^{[\ell]}$ is defined as in Definition \reff{F-IT} using the particular sequence $\vec{\lambda}$.
\end{conjecture}
\par Our main results, Theorem \reff{MAIN-RT} and Theorem \reff{X-LSGM} below, state that Conjecture \reff{BR-IT} holds in two special cases, when either (i) $m=2$ and $\dim X=2$, or (ii) $m=1$ and $X$ is an arbitrary Banach space. Note that in both of these cases the above mentioned  finiteness number $N(m,\X)$ equals $4$.
\begin{theorem}\label{MAIN-RT} Let $\mfM=\MR$ be a pseudometric space, and let $X$ be a two dimensional Banach space. Let $m=2$ so that the number $\ell(m,X)=2$. In this case Conjecture \reff{BR-IT} holds for every $\lambda_1,\lambda_2$ and $\gamma$ such that
\begin{align}\lbl{GM-FN}
\lambda_1\ge e(\mfM,X), ~~~~~~\lambda_2\ge 3\lambda_1,
~~~~~~
\gamma\ge \lambda_2\,(3\lambda_2+\lambda_1)^2/
(\lambda_2-\lambda_1)^2.
\end{align}
Here $e(\mfM,X)$ denotes the Lipschitz extension constant of $X$ with respect to $\mfM$. (See Definition \reff{LIP-C}.)
\par Thus, the following statement is true: Let $F:\Mc\to\Kc(\X)$ be a set-valued mapping from a pseudometric space $(\Mc,\rho)$ into the family $\Kc(X)$ of all non-empty convex compact subsets of $\X$. Given $x\in\Mc$ let
\begin{align}\lbl{F12-LTI}
F^{[1]}(x)=
\bigcap_{z\in\Mc}\,
\left[F(z)+\lambda_1\,\rho(x,z)\BX\right],~~~~
F^{[2]}(x)=\bigcap_{z\in\Mc}\,
\left[F^{[1]}(z)+\lambda_2\,\rho(x,z)\BX\right].
\end{align}
\par Suppose that for every subset $\Mc'\subset\Mc$ with $\#\Mc'\le 4$, the restriction $F|_{\Mc'}$ of $F$ to $\Mc'$ has a Lipschitz selection with Lipschitz seminorm at most $1$.
\par Then for every $\lambda_1,\lambda_2$ and $\gamma$ satisfying \rf{GM-FN} the set
\begin{align}\lbl{F2-NEM}
F^{[2]}(x)\ne\emp~~~~\text{for every}~~~~x\in\Mc.
\end{align}
Furthermore,
\begin{align}\lbl{HD-RT}
\dhf(F^{[2]}(x),F^{[2]}(y))\le \gamma\,\rho(x,y)
~~~~\text{for every}~~~~x,y\in\Mc.
\end{align}
\par If $X$ is a Euclidean two dimensional space, \rf{F2-NEM} and  \rf{HD-RT} hold when \rf{GM-FN} is replaced by the weaker requirements that
\begin{align}\lbl{X-HSP}
\lambda_1\ge e(\mfM,X), ~~~~~~\lambda_2\ge 3\lambda_1,
~~~~~~
\gamma\ge \,\lambda_2\left\{1+2\lambda_2\,\slbig
\left(\lambda_2^2-\lambda_1^2\right)^{\frac12}\right\}^2.
\end{align}
\end{theorem}
\smsk
\par In particular, in Section 3 we show that the mapping $F^{[2]}$ satisfies \rf{F2-NEM} and \rf{HD-RT} whenever $X$ is an {\it arbitrary} two dimensional Banach space and $\lambda_1=4/3$, $\lambda_2=4$, $\gamma=100$. If $X$ is also Euclidean, then one can set $\lambda_1=4/\pi$, $\lambda_2=12/\pi$ and $\gamma=38$. Furthermore, we prove that if $\Mc$ is a subset of a Euclidean space $E$, $\rho$ is the Euclidean metric in $E$, and $\X$ is a two dimensional Euclidean space, then properties \rf{F2-NEM} and \rf{HD-RT} hold for $\lambda_1=1$, $\lambda_2=3$, and $\gamma=25$.
\smsk
\par In Section 5 we prove Theorem \reff{LTI-M} which refines the result of Theorem \reff{MAIN-RT} for the space $X=\LTI$, i.e., for $\RT$ equipped with the norm $\|x\|=\max\{|x_1|,|x_2|\}$, $x=(x_1,x_2)$. More specifically, we show that in this case properties \rf{F2-NEM} and \rf{HD-RT} hold whenever
$\lambda_1\ge 1$, $\lambda_2\ge 3\lambda_1$, and $\gamma\ge \lambda_2\, (3\lambda_2+\lambda_1)/(\lambda_2-\lambda_1)$.
In particular, these properties hold for
$\lambda_1=1$, $\lambda_2=3$ and $\gamma=15$.
\smsk
\par Let us now explicitly formulate the above mentioned
second main result of the paper. We prove it in Section 4. It deals with set-valued mappings from a pseudometric space into the family $\Kc_1(\X)$ of all {\it bounded closed line segments} of an arbitrary Banach space $\X$.
\begin{theorem}\label{X-LSGM} Let $\MR$ be a pseudometric space. Let $m=1$ and let $X$ be a Banach space with $\dim X>1$; thus, $\ell(m,X)=2$, see \rf{NMY-1}. In this case Conjecture \reff{BR-IT} holds for every  $\lambda_1,\lambda_2$ and $\gamma$ such that
\begin{align}\lbl{GM-FN-1}
\lambda_1\ge 1, ~~~~~~\lambda_2\ge 3\lambda_1,
~~~~~~
\gamma\ge \lambda_2\,(3\lambda_2+\lambda_1)/
(\lambda_2-\lambda_1).
\end{align}
\par Thus, the following statement is true: Let $F:\Mc\to\Kc_1(\X)$ be a set-valued mapping such that for every subset $\Mc'\subset\Mc$ with $\#\Mc'\le 4$, the restriction $F|_{\Mc'}$ of $F$ to $\Mc'$ has a Lipschitz selection with Lipschitz seminorm at most $1$.
\par Let $F^{[2]}$ be the mapping defined by \rf{F12-LTI}. Then properties \rf{F2-NEM} and \rf{HD-RT} hold whenever $\lambda_1$, $\lambda_2$ and $\gamma$ satisfy \rf{GM-FN-1}. In particular, one can set $\lambda_1=1$, $\lambda_2=3$ and $\gamma=15$.
\par If $X$ is a Euclidean space, the same statement is also true whenever, instead of \rf{GM-FN-1}, $\lambda_1,\lambda_2$ and $\gamma$ satisfy the weaker  condition
\begin{align}\lbl{X1-H}
\lambda_1\ge 1, ~~~~~~\lambda_2\ge 3\lambda_1,
~~~~~~
\gamma\ge \,\lambda_2+2\lambda_2^2\,\slbig
\left(\lambda_2^2-\lambda_1^2\right)^{\frac12}.
\end{align}
\par In particular, in this case, \rf{F2-NEM} and \rf{HD-RT} hold whenever $\lambda_1=1$, $\lambda_2=3$ and $\gamma=10$.
\end{theorem}
\smsk
\par In Section 5.1 we note that Conjecture \reff{BR-IT} also holds for a one dimensional space $X$ and $m=1$. In this case the statement of the conjecture is true for every $\lambda_1\ge 1$ and $\gamma\ge 1$. See Proposition \reff{X-1DIM}.
\smsk
\par Note that Theorem \reff{MAIN-RT} tells us that for every set-valued mapping $F$ satisfying the hypothesis of this theorem, the mapping $F^{[2]}$ determined by \rf{F12-LTI} with $\lambda_1=4/3$ and $\lambda_2=4$ provides a $\gamma$-core of $F$ with $\gamma=100$. (See Definition \reff{D-CORE}.) In turn, Theorem \reff{X-LSGM} states that the mapping $F^{[2]}$ corresponding to the parameters $\lambda_1=1$ and $\lambda_2=3$ is a $15$-core of any $F$ satisfying the conditions of this theorem.
\smsk
\par We note that the proofs of Theorem
\reff{MAIN-RT} and Theorem \reff{X-LSGM} rely on Helly's Intersection Theorem and a series of auxiliary results about neighborhoods of intersections of convex sets. See Section 2.

\begin{remark} {\em Let us compare Conjecture \reff{BR-IT} (and Theorems \reff{MAIN-RT} and \reff{X-LSGM}) with
the Finiteness Principle (FP) formulated in Theorem \reff{MAIN-FP}. First we note that FP is invariant with respect to the transition to an equivalent norm on $X$, while the statement of Conjecture \reff{BR-IT} is not.
\par To express this more precisely, let $\|\cdot\|_1$ and $\|\cdot\|_2$ be two equivalent norms on $X$, i.e., suppose that for some $\alpha\ge 1$  the following inequality
$(1/\alpha)\,\|\cdot\|_1\le\|\cdot\|_2\le \alpha\,\|\cdot\|_1$ holds. Clearly, if FP holds for  $(X,\|\cdot\|_1)$ then it immediately holds also for $(X,\|\cdot\|_2)$ (with the constant $\alpha^2\gamma$  instead of $\gamma$). However the validity of Conjecture \reff{BR-IT} for the norm $\|\cdot\|_1$ does not imply its validity for an equivalent norm $\|\cdot\|_2$ on $X$ (at least we do not see any obvious way for obtaining such an implication). For example, the validity of Conjecture \reff{BR-IT} in $\ell^n_\infty$ (i.e., $\RN$ equipped with the uniform norm) does not automatically imply its validity in the space $\ell^n_2$ (i.e., $\RN$ with the Euclidean norm).
\smsk
\par  We also note the following: in a certain sense,
the result of Theorem \reff{MAIN-RT} is ``stronger'' than Theorem \reff{MAIN-FP} (i.e., FP for the case of a two dimensional Banach space $X$). Indeed, in this case, the hypotheses of FP and Theorem \reff{MAIN-RT} coincide. Moreover, Theorem \reff{MAIN-RT} ensures that the set-valued mapping $F^{[2]}$ is a core of $F$. This
property of $F^{[2]}$ implies, via arguments in \cite{FS-2018} that the function
\begin{align}\lbl{FST-C}
f(x)=\ST\,(F^{[2]})\,(x),~~~~~~x\in\Mc,
\end{align}
is a Lipschitz selection of $F$. Here $\ST:\KM\to\X$  is the Steiner-type point map \cite{S-2004}.
\par Thus, FP (in the two dimensional case) follows immediately from Theorem \reff{MAIN-RT}. However, it is absolutely unclear how the statement of Theorem \reff{MAIN-RT} can be deduced from FP. I would like to thank Charles Fefferman who kindly drew my attention to this interesting fact.\rbx}
\end{remark}
\par Let us reformulate Conjecture \reff{BR-IT} in a way  which {\it does not require the use of the notion of a core of a set-valued mapping}. We recall that the mapping $F^{[\ell]}:\Mc\to \KM$ which appears in Conjecture \reff{BR-IT} is a $\gamma$-core of $F$ if
$$
\dhf(F^{[\ell]}(x),F^{[\ell]}(y))\le \gamma\, \rho(x,y)~~~\text{for all}~~~x,y\in\Mc.
$$
See part (ii) of Definition \reff{D-CORE}. Hence, given $x\in\Mc$,
\begin{align}\lbl{IM-L}
F^{[\ell]}(x)\subset F^{[\ell]}(y)+
\gamma\,\rho(x,y) \BX~~~\text{for every}~~~y\in\Mc.
\end{align}
\par Let
$$
F^{[\ell+1]}(x)=\BAL{F^{[\ell]}}{\gamma}{\rho}(x)=
\bigcap_{y\in \Mc}\,
\left[F^{[\ell]}(y)+\gamma\,\rho(x,y)\,\BX\right].
$$
Cf. \rf{IT-F}. This and \rf{IM-L} imply the inclusion
$F^{[\ell+1]}(x)\supset F^{[\ell]}(x)$, $x\in\Mc$. On the other hand, \rf{CH-N} tells us that $F^{[\ell+1]}(x)\subset F^{[\ell]}(x)$ proving that $F^{[\ell+1]}=F^{[\ell]}$ on $\Mc$.
\smsk
\par These observations enable us to reformulate Conjecture \reff{BR-IT} as follows.
\begin{conjecture}\label{ST-IT} Let $\MR$ be a pseudometric space, and let $\X$ be a Banach space. Let $m$ be a fixed positive integer and let $\ell=\ell(m,X)$, see\rf{NMY-1}.
\par There exists a sequence $\vec{\lambda}=\left\{\lambda_{k}:1\le k\le\ell+1\right\}$ of $\ell+1$ numbers $\lambda_{k}$ all satisfying $\lambda_{k}\ge1$ such that, for every set-valued mapping $F:\Mc\to\KM$ satisfying the hypothesis of the Finiteness Principle (Theorem \reff{MAIN-FP}), the family $\{F^{[k]}: k=1,...,\ell+1\}$ of  set-valued mappings constructed by formula \rf{IT-F} has the following property:
\begin{align}\lbl{ST-PR}
F^{[\ell]}(x)\ne\emp~~~~~\text{and}
~~~~~F^{[\ell+1]}(x)=F^{[\ell]}(x)~~~~ \text{for all} ~~~
x\in\Mc.
\end{align}
\end{conjecture}
\par We refer to \rf{ST-PR} as a {\it Stabilization Property} of balanced refinements.
\smsk
\par Thus, Theorem \reff{MAIN-RT} and Theorem \reff{X-LSGM} tell us that a Stabilization Property of ba\-lanced refinements holds whenever $\dim\X=2$ or $m=1$ (and $\X$ is an arbitrary). More specifically, Theorem \reff{MAIN-RT} shows that if $m=2$ and $\dim \X=2$, Conjecture \reff{ST-IT} holds with $\ell=2$ and $\vec{\lambda}=\{4/3,4,10^2\}$.
\par In other words, in this case, $F^{[2]}(x)\ne\emp$ for each $x\in\Mc$ and $F^{[3]}=F^{[2]}$ on $\Mc$. In turn, Theorem \reff{X-LSGM} states that the same property holds whenever $X$ is an arbitrary Banach space, $m=1$, and $\vec{\lambda}=\{1,3,15\}$.
\smsk
\msk
\par {\bf Acknowledgements.} I am very thankful to Michael Cwikel for useful suggestions and remarks. I am also very grateful to Charles Fefferman for stimulating discussions and valuable advice.
\par The results of this paper were presented at the 12th Whitney Problems Workshop, August 2019, the University of Texas at Austin, TX. I am very thankful to all participants of that workshop for valuable conversations and useful remarks. I am also grateful to the Conference Board of the
Mathematical Sciences and the University of Texas for supporting the Austin workshop.

\SECT{2. Neighborhoods of intersections of convex sets
in a Banach space.}{2}


\indent
\par We first need to fix some notation.
\par Let $(\X,\|\cdot\|)$ be a Banach space, and let $\BX$ be the unit ball in $X$. Let $A$ and $B$ be non-empty subsets of $X$. We let $A + B=\{a + b: a\in A, b\in B\}$ denote the Minkowski sum of these sets. By $\CX$ we denote the family of all convex closed non-empty subsets of $\X$.
We equip this family with the Hausdorff distance.
\par Sometimes, for a given set $\Mc$, we will be looking simultaneously at two distinct pseudometrics on $\Mc$, say $\rho$ and $\delta$. In this case we will speak of a $\rho$-Lipschitz selection and $\rho$-Lipschitz seminorm, or a $\delta$-Lipschitz selection and $\delta$-Lipschitz seminorm to make clear which pseudometric we are using. Furthermore, given a mapping $f:\Mc\to X$ we will write
$\|f\|_{\Lip((\Mc;\rho),\X)}$ or  $\|f\|_{\Lip((\Mc;\delta),\X)}$ to denote the Lipschitz seminorm of $f$ with respect to the pseudometric $\rho$ or $\delta$ respectively.
\par For each finite set $S$, we let $\#S$ denote the number of elements of $S$.

\smsk
\par Given a Banach space $\X$, Przes{\l}awski and Yost \cite{PY2} introduced an important geometrical characteristic of $\X$, the so-called {\it modulus of squareness} of $\X$. Let us recall its definition.
\par We observe that for any $x,y\in\X$ with $\|y\|<1<\|x\|$ there exists a {\it unique} $z=z(x,y)$ with $\|z\|=1$ which belongs to the line segment $[x,y]$. We set
\bel{D-OM-XY}
\omega(x,y)=\frac{\|x-z(x,y)\|}{\|x\|-1}
\ee
and define a function $\xi:[0,1)\to [1,\infty)$ by
\bel{D-XI-B}
\xi_{\X}(\beta)=\sup\,\{\omega(x,y): x,y\in X,\,\|y\|\le \beta<1<\|x\|\,\}.
\ee
We also put
\bel{FI-PS}
\vf(\beta)=(1+\beta)/(1-\beta)~~~~~\text{and}~~~~~ \psi(\beta)=(1-\beta^2)^{-\frac12},~~~~\beta\in[0,1).
\ee
\par It is shown in \cite{PY2}, that for any Banach space $\X$
\bel{F-B}
\xi_X(\beta)\le \vf(\beta)~~~~\text{for every}~~~\beta\in[0,1),
\ee
and
\bel{PS-H}
\xi_X(\beta)=\psi(\beta)~~~~\text{for every}~~~\beta\in[0,1),
\ee
provided $X$ is a Euclidean space.
\msk
\begin{theorem}(\cite[Theorem 4]{PY2})\lbl{PY-M} Let $(S,\delta)$ be a metric space, let $\X$ be a Banach space, and let $f:S\to\X$ and $g:S\to[0,\infty)$ be Lipschitz mappings. Let $F:S\to\CX$ be a Lipschitz (with respect to the Hausdorff distance) set-valued mapping.
\msk
\par Suppose that there exists a constant $\gamma>1$ such that $g(x)\ge\gamma\dist(f(x),F(x))$ for every $x\in S$. Then the intersection mapping $G:S\to\CX$ defined by
$$
G(x)=F(x)\cap \BX(f(x),g(x))
$$
is Lipschitz continuous on $S$ (with respect to $\dhf$) with Lipschitz seminorm
$$
\|G\|_{\Lip(S,\CX)}\le \|F\|_{\Lip(S,\CX)}+
(\,\|F\|_{\Lip(S,\CX)}+\|f\|_{\Lip(S,\X)}+\|g\|_{\Lip(S,\R)})
\,\xi(1/\gamma).
$$
\end{theorem}
\par This theorem implies the following
\begin{proposition}\lbl{HD-IM} Let $X$ be a Banach space, $a\in \X$, $r\ge 0$, and let $C\subset \X$ be a convex set. 
\par Suppose that $C\cap \BX(a,r)\ne\emp$. Then for every $s>0$ and every $L>1$ the following inequality
$$
\dhf\left(C\cap\BX(a,Lr),(C+s\BX)\cap\BX(a,Lr+s)\right)\le
\left(1+2\,\xi_{\X}\left(\tfrac{1}{L}\right)\right)\,s
$$
holds.
\end{proposition}
\par {\it Proof.} Let $S=\{x,y\}\subset\R$ where $x=0$ and $y=s$, and let $\delta(x,y)=s$.
\par We define a mapping $f:S\to X$ and a function $g:S\to\R$ by letting $f(x)=f(y)=a$ and $g(x)=Lr, g(y)=Lr+s$. Clearly, $\|f\|_{\Lip(S,\X)}=0$, and $\|g\|_{\Lip(S,\R)}=1$.
\par We put $\gamma=L$. We note that $C\cap \BX(a,r)\ne\emp$ so that
$$
\dist(f(x),F(x))=\dist(a,C)\le r.
$$
Hence, $g(x)=Lr=\gamma\,r\ge \gamma\,\dist(f(x),F(x))$.
\smsk
\par Then we define a mapping $F:S\to\CX$ by setting  $F(x)=C$ and $F(y)=C+s\BX$. Clearly,
$$
\dhf(F(x),F(y))\le s=\delta(x,y)~~~~\text{so that}~~~~~
\|F\|_{\Lip(S,\CX)}\le 1.
$$
\par Thus, the conditions of Theorem \reff{PY-M} are satisfied for the metric space $(S,\delta)$ and the mappings $f$, $g$ and $F$. This theorem tells us that
the mapping $G:S\to\CX$ defined by
$$
G(u)=F(u)\cap \BX(f(u),g(u)),~~~~~~u\in S,
$$
is Lipschitz on $S$ with respect to the Hausdorff distance. Furthermore,
$$
\|G\|_{\Lip(S,\CX)}\le \|F\|_{\Lip(S,\CX)}+
(\,\|F\|_{\Lip(S,\CX)}+\|f\|_{\Lip(S,\X)}+\|g\|_{\Lip(S,\R)})
\,\xi(1/\gamma)\le 1+2\,\xi_{\X}\left(\tfrac{1}{L}\right).
$$
Hence,
\be
\dhf\left(C\cap\BX(a,Lr),(C+s\BX)\cap\BX(a,Lr+s)\right)
&=&
\dhf(G(x),G(y))\nn\\
&\le& \|G\|_{\Lip(S,\CX)}\,\delta(x,y)
\le
\left(1+2\,\xi_{\X}\left(\tfrac{1}{L}\right)\right)\,s
\nn
\ee
proving the proposition.\bx
\msk
\par Proposition \reff{HD-IM} implies the following important
\begin{theorem}\lbl{N-S} Let $X$ be a Banach space, and let $C\subset \X$ be a convex set. Let $a\in \X$ and let $r\ge 0$. Suppose that 
\bel{C-PR}
C\cap \BX(a,r)\ne\emp.
\ee
\par Then for every $s>0$ and $L>1$
\bel{IN-MN}
[C\cap \BX(a,Lr)]+\theta(L)\,s\,\BX
\supset (C+s\BX)\cap \BX(a,Lr+s)
\ee
where
\bel{TH-L}
\theta(L)=(3L+1)/(L-1).
\ee
\par If $X$ is a Euclidean space then \rf{IN-MN} holds with
\bel{TH-L-H}
\theta(L)=1+\frac{2L}{\sqrt{L^2-1}}\,.
\ee
\end{theorem}
\par {\it Proof.} Let
$$
G=C\cap \BX(a,Lr)~~~~\text{and}~~~~
\tG=(C+s\BX)\cap \BX(a,Lr+s).
$$
\par Definition \rf{HFD} tells us that $\tG\subset G+\dhf(G,\tG)\,\BX$. In turn, Proposition \reff{HD-IM} states that
$$
\dhf(G,\tG)\le
\left(1+2\,\xi_{\X}\left(\tfrac{1}{L}\right)\right)\,s.
$$
Hence,
$$
\tG\subset G+\Theta(L)\,s\,\BX
~~~~~\text{where}~~~~~
\Theta(L)=1+2\,\xi_{\X}\left(\tfrac{1}{L}\right).
$$
\par Now, let $X$ be an arbitrary Banach space. In this case, thanks to \rf{FI-PS} and \rf{F-B}, we have
$$
\Theta(L)= 1+2\,\xi_{\X}\left(\tfrac{1}{L}\right)\le 1+
2\frac{1+1/L}{1-1/L}=\frac{3L+1}{L-1}.
$$
This inequality and \rf{TH-L} imply inclusion \rf{IN-MN} in the case under consideration.
\msk
\par Finally, let $X$ be a Euclidean space. In this case, from \rf{FI-PS}, \rf{PS-H} and \rf{TH-L-H}, we have
$$
\Theta(L)=1+2\,\xi_{\X}\left(\tfrac{1}{L}\right)=
1+2\,(1-(1/L)^2)^{-\frac12}=1+\frac{2L}{\sqrt{L^2-1}}
=\theta(L).
$$
\par The proof of the theorem is complete.\bx
\msk

\par For the case  of an arbitrary Banach space $\X$, Theorem \reff{N-S} was proved by Przes{\l}awski and Rybinski \cite[p. 279]{PR}. For the case of a Euclidean space $\X$ see Przes{\l}awski, Yost \cite[Theorem 4]{PY2}. For similar results we refer the reader to \cite{Ar}, \cite[p. 369]{AF} and \cite[p. 26]{BL}.

\msk
\par For the sake of completeness, and for the reader's convenience, below we give

\msk
\par {\it A direct proof of Theorem \reff{N-S}.} We follow  the proof of Lemma 5.3 from \cite[p. 279]{PR}. If $r=0$ then \rf{IN-MN} holds trivially, so we assume that $r>0$. Without loss of generality, we may also assume that $a=0$. Thus we should prove that
$$
[C\cap (Lr\BX)]+\theta s\,\BX\supset
(C+s\BX)\cap(Lr\BX+s\BX)
$$
provided $r>0$, $s>0$, $L>1$. Let
\bel{IN-Z}
z\in(C+s\BX)\cap(Lr\BX+s\BX)=(C+s\BX)\cap[(Lr+s)\BX].
\ee
Prove that
\bel{F-Z}
z\in [C\cap(Lr\BX)]+\theta s\,\BX.
\ee
\par Thanks to \rf{IN-Z}, $z\in(C+s\BX)$ so that there exists an element $v\in C$ such that
\bel{VV-Z}
\|v-z\|\le s.
\ee
\par If $\|v\|\le Lr$, then $v\in C\cap(Lr\BX)$ proving \rf{F-Z}.
\msk
\par Suppose that
\bel{V-Z}
\|v\|>Lr.
\ee
Property \rf{IN-Z} tells us that $\|z\|\le Lr+s$ so that
\bel{V-NR}
\|v\|\le \|z\|+s\le Lr+2s.
\ee
In turn, assumption \rf{C-PR} tells us that there exists an element $v'\in C$ such that
\bel{VP-NR}
\|v'\|\le r<Lr.
\ee
\par Choose $\lambda\in(0,1)$ such that the element
$$
\tv=\lambda v'+(1-\lambda)v
$$
has the norm $\|\tv\|=Lr$. See Fig. 5.
\bsk\msk

%
\begin{figure}[h!]
\centering{\includegraphics[scale=0.20]{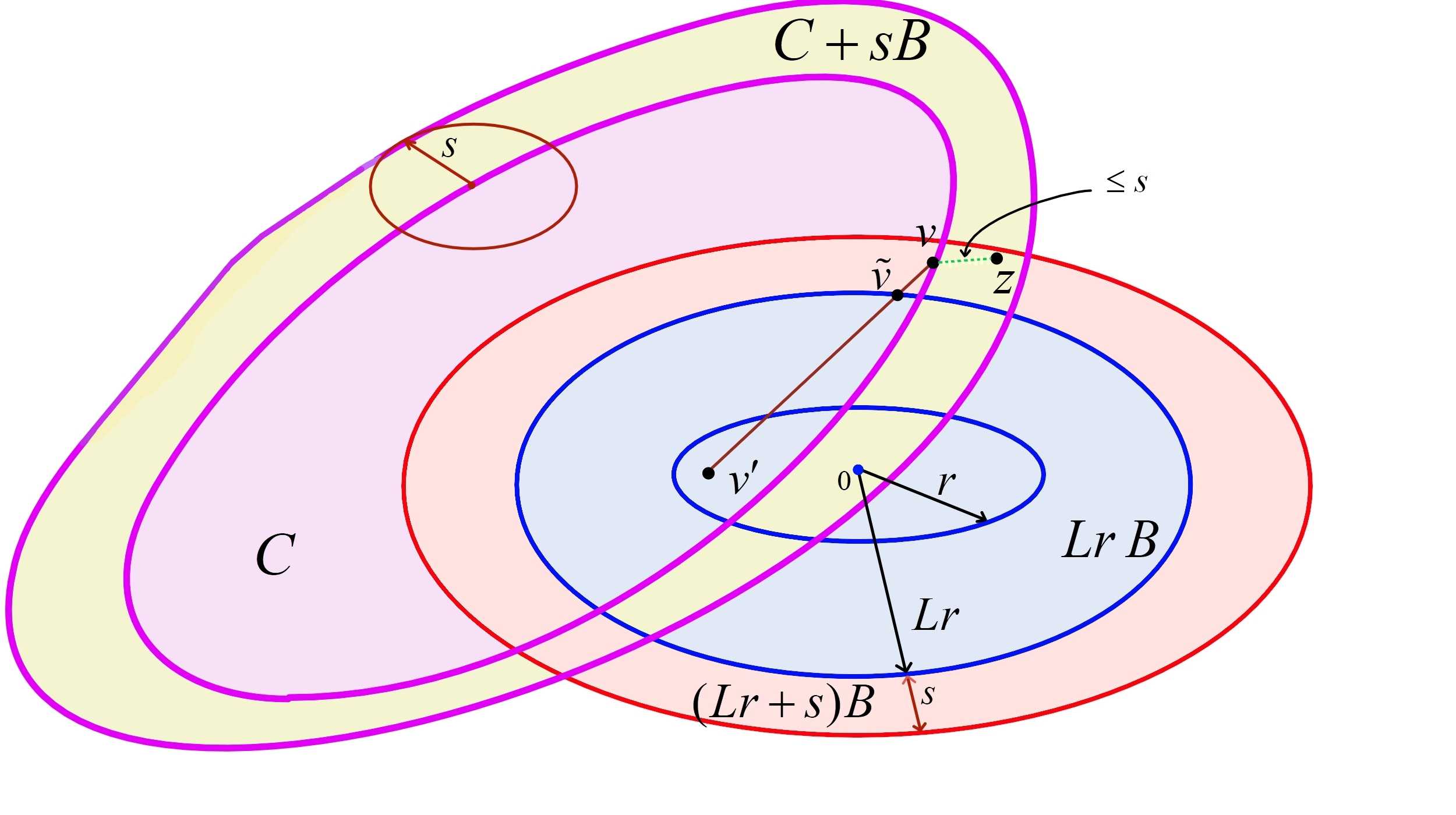}}
\caption{The points $z$, $v$, $\tv$ and $v'$.}
\end{figure}
\msk

\par We know that $C$ is convex so that $[v',v]\subset C$ proving that
\bel{VP-D}
\tv\in C\cap (Lr\BX).
\ee
\par Thanks to \rf{VP-NR}, \rf{V-NR} and the triangle inequality,
$$
Lr=\|\tv\|=\|\lambda v'+(1-\lambda)v\|\le \lambda r+(1-\lambda)(Lr+2s)
$$
proving that
$$
\lambda\le \frac{2s}{(L-1)r+2s}.
$$
Consequently, thanks to this inequality, \rf{VP-NR} and \rf{V-NR}
$$
\|v-\tv\|=\lambda\|v-v'\|\le  \lambda(\|v\|+\|v'\|)
\le \frac{2s}{((L-1)r+2s)}\cdot(Lr+2s+r)\le 2s(L+1)/(L-1).
$$
From this inequality and \rf{VV-Z} we have
$$
\|z-\tv\|\le \|z-v\|+\|v-\tv\|\le s+2s(L+1)/(L-1)
=\theta(L)\, s
$$
which together with \rf{VP-D} implies \rf{F-Z}.
\msk
\par Let now $X$ be a {\it Euclidean space}. We modify the above proof after \rf{VP-D} as follows.
\par We put $\beta=1/L$, and
\bel{XYW}
x=\mathlarger{\tfrac{1}{L r}}\,v,
~~~~~~y=\mathlarger{\tfrac{1}{L r}}\,v',~~~~~~
w=\mathlarger{\tfrac{1}{L r}}\,\tv.
\ee
\par Then, thanks to \rf{V-Z} and \rf{VP-NR},
\bel{XYB}
\|y\|\le \beta<1<\|x\|.
\ee
\par We note that for any $u,\tu\in X$ such that $\|\tu\|<1<\|u\|$, there exists a unique point $w=w(u,\tu)\in[u,\tu]$ with $\|w\|=1$. Hence, thanks to
\rf{D-OM-XY},
\bel{D-OM}
\omega(u,\tu)=\frac{\|u-w(u,\tu)\|}{\|u\|-1}.~~~~~~~~
\text{See Fig. 6.}
\ee
\begin{figure}[h!]
\centering{\includegraphics[scale=0.3]{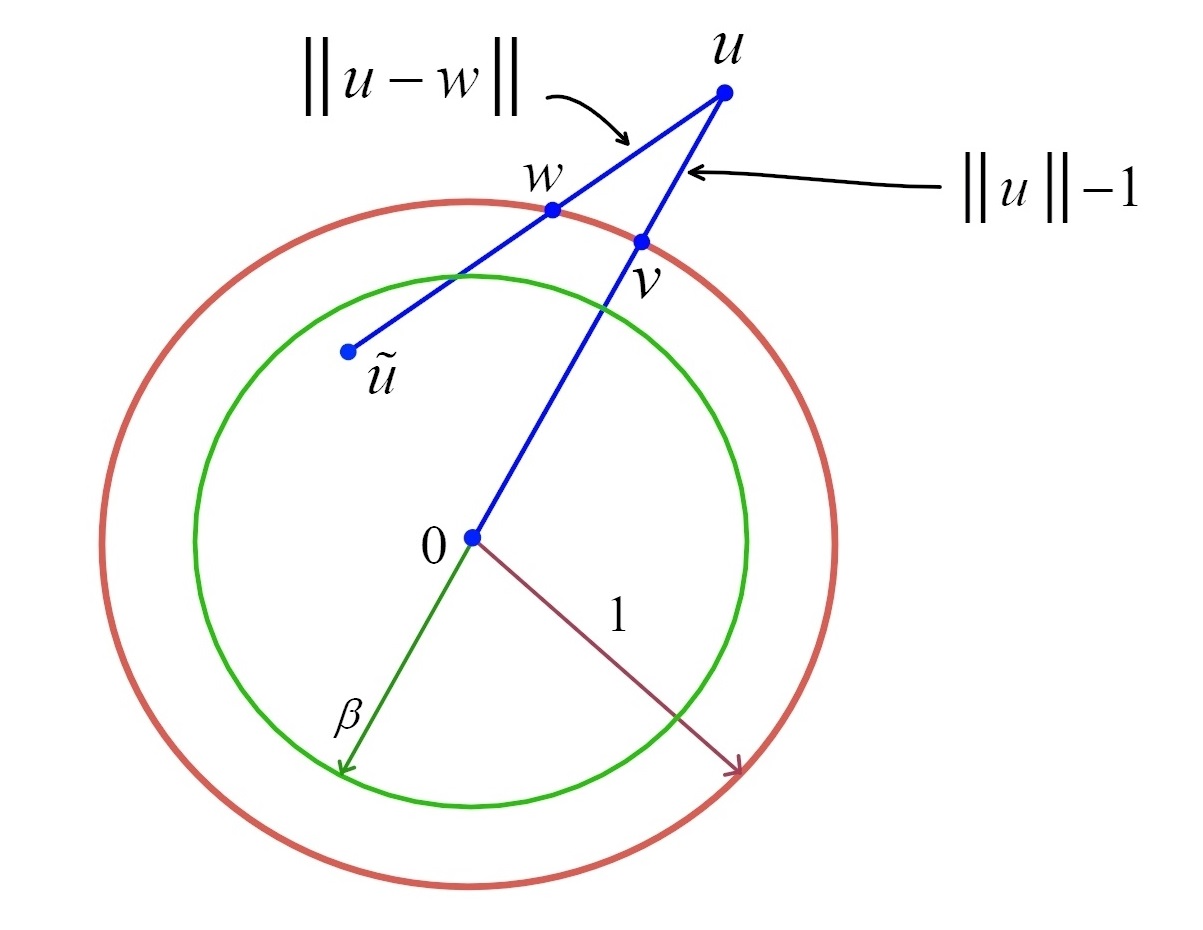}}
\caption{The points $u$, $w$ and $\tu$.}
\end{figure}
\par We also recall the definition of the function $\xi_X$, see \rf{D-XI-B}:
\bel{D-XI}
\xi_X(\beta)=\sup\,\{\omega(u,\tu): u,\tu\in X,\,\|\tu\|\le \beta<1<\|u\|\,\}.
\ee
\par Prove that
\bel{XI}
\xi_X(\beta)=(1-\beta^2)^{-\frac12}.
\ee
\par In fact, fix $u$ with $\|u\|>1$. One can easily see that $\sup\,\{\omega(u,\tu): \|\tu\|\le \beta\}$ is attained for some $\tu$ with $\|\tu\|=\beta$, and the line segment $[\tu,u]$ touches the sphere centered at the origin with radius $\beta$. Thus, $u-\tu$ is perpendicular to $\tu$. See Fig. 7.
\begin{figure}[h!]
\centering{\includegraphics[scale=0.30]{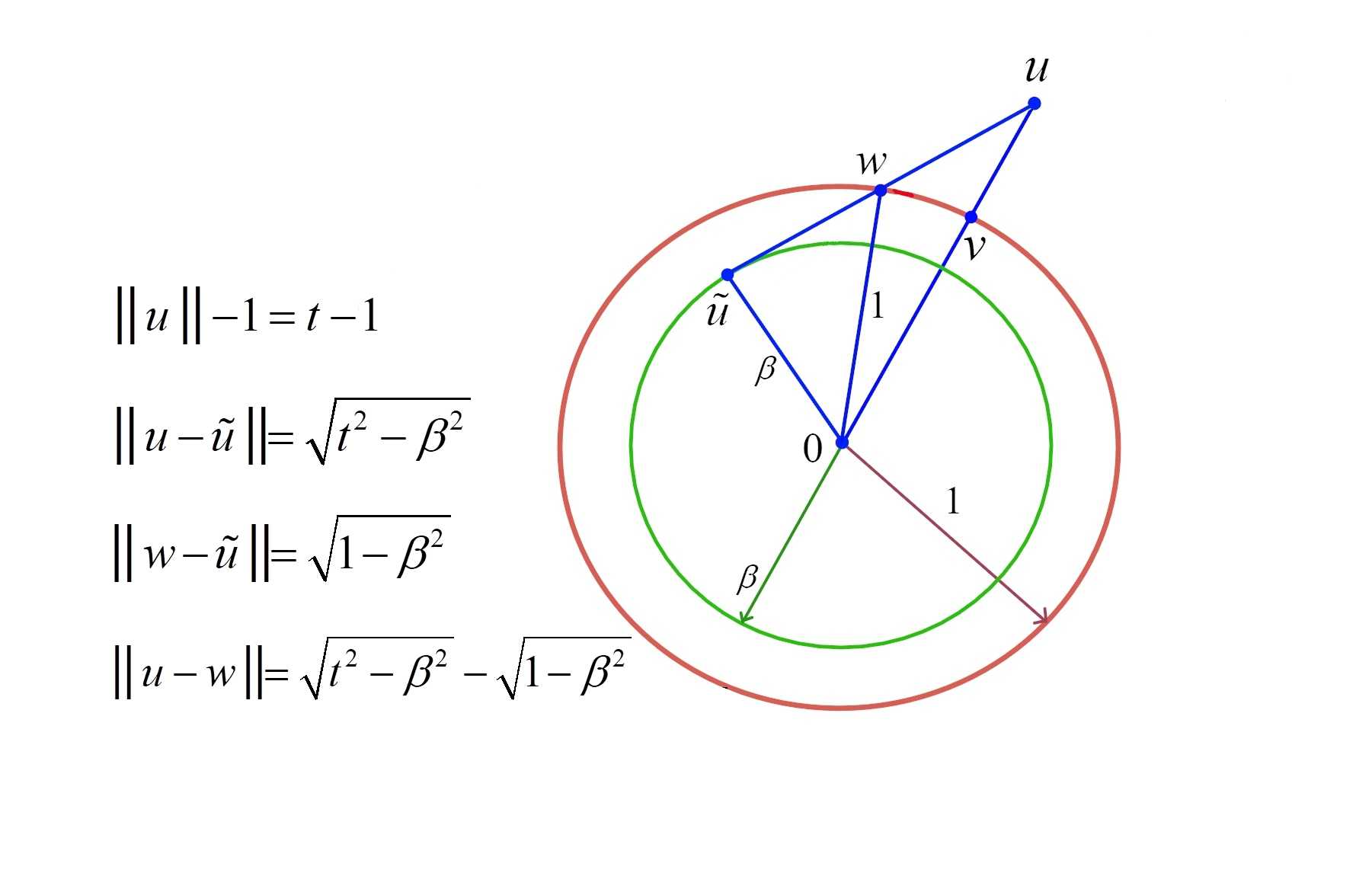}}
\vspace*{-15mm}
\caption{The optimal position of the point $\tu$ for given $u$ with $\|u\|=t>1$.}
\end{figure}
\par Hence, $\sup\{\omega(u,\tu): \|\tu\|\le \beta\}=f(\|u\|)$ where
$$
f(t)=
\frac{\sqrt{t^2-\beta^2}-\sqrt{1-\beta^2}}{t-1},~~~~
t>1.
$$
\par The function $f$ is decreasing on $(1,+\infty)$ so that
$$
\sup_{t>1} f(t)=\lim_{t\to 1} f(t)=(1-\beta^2)^{-\frac12}
$$
proving \rf{XI}.
\msk
\par We apply formula \rf{XI} to the points $x,y,w$ defined
by \rf{XYW}, and to $\beta=1/L$. We have
$$
\frac{\|v-\tv\|}{\|v\|-Lr}=\frac{\|x-w\|}{\|y\|-1}=
\omega(x,y).
$$
See \rf{D-OM}. Thanks to \rf{XYB} and \rf{D-XI},
$$
\frac{\|v-\tv\|}{\|v\|-Lr}\le
\sup\,\{\omega(u,\tu): u,\tu\in X,\,\|\tu\|\le \beta<1<\|u\|\,\}=\xi_X(\beta)=\xi_X(1/L)
$$
so that, thanks to \rf{XI},
$$
\frac{\|v-\tv\|}{\|v\|-Lr}\le
\xi_X(1/L)=\frac{L}{\sqrt{L^2-1}}.
$$
In turn, thanks to \rf{V-NR}, $\|v\|-Lr\le 2s$, so that
$$
\|v-\tv\|\le \frac{L}{\sqrt{L^2-1}}\,(\|v\|-Lr)
\le \frac{2sL}{\sqrt{L^2-1}}.
$$
This inequality and \rf{VV-Z} imply the following:
$$
\|z-\tv\|\le \|z-v\|+\|v-\tv\|\le s+\frac{2sL}{\sqrt{L^2-1}}=
\left(1+2L/\sqrt{L^2-1}\,\right)\,s.
$$
\par This and \rf{VP-D} imply \rf{F-Z} with $\theta=\theta(L)$ defined by \rf{TH-L-H}.
\par The proof of Theorem \reff{N-S} is complete.\bx
\par For the case  of an arbitrary Banach space $\X$, Theorem \reff{N-S} was proved by Przes{\l}awski and Rybinski \cite[p. 279]{PR}. For the case of a Euclidean space $\X$ see Przes{\l}awski, Yost \cite[Theorem 4]{PY2}. For similar results we refer the reader to \cite{Ar}, \cite[p. 369]{AF} and \cite[p. 26]{BL}.
\begin{remark} {\em In general, assumption \rf{C-PR} in Theorem \reff{N-S} cannot be omitted. For instance, let $X=\ell^2_\infty$ and let $L=2$; then
$\theta(L)=7$. Let $r\gg s>0$, $a\in\RT$, and  let $C$ be a half-plane shown in Fig. 8. (In Fig. 8  -- 11 we omit $X$ in the notation of $\ell^2_\infty$-balls.)
\begin{figure}[h!]
\centering{\includegraphics[scale=0.60]{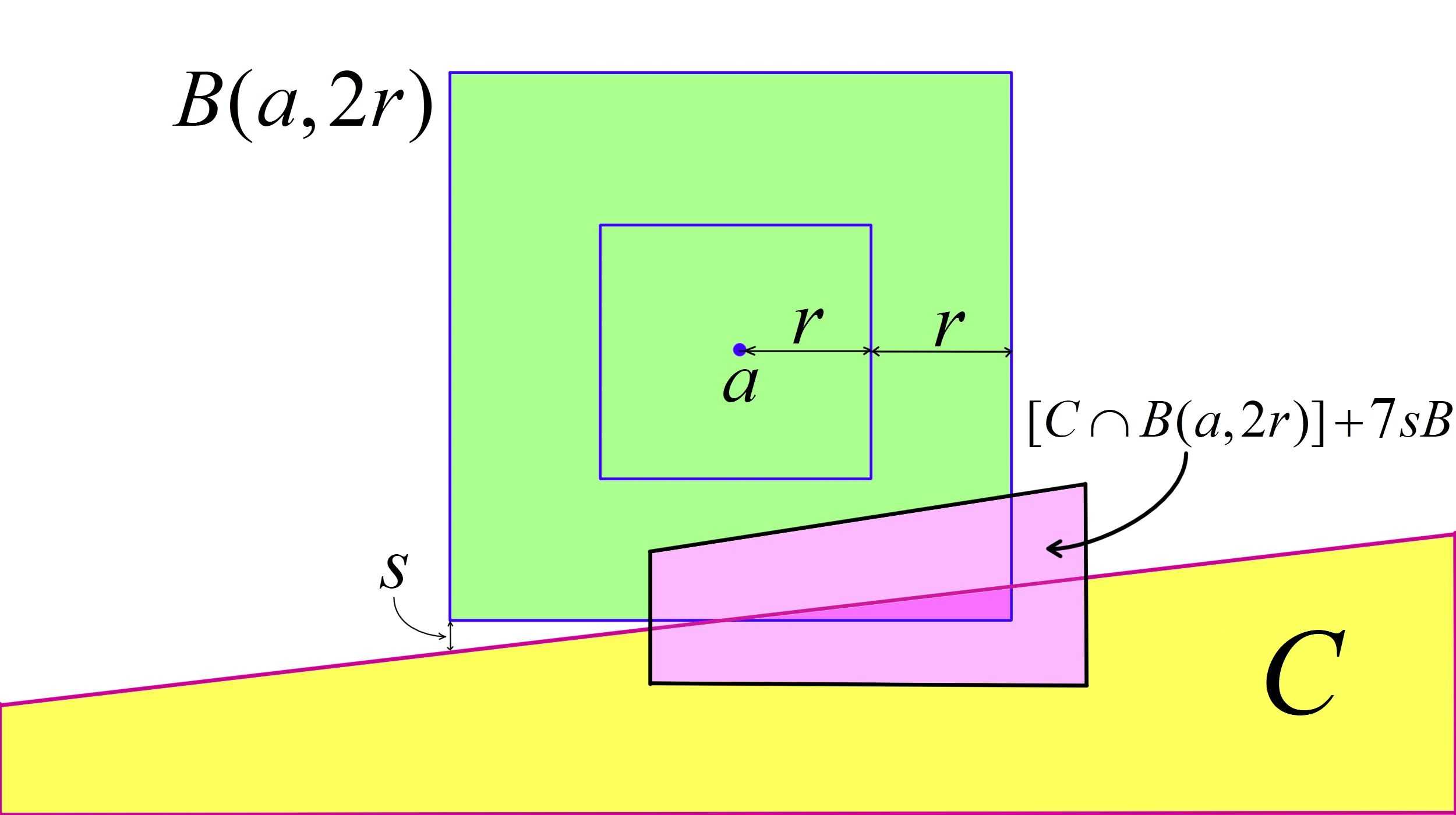}}
\caption{The half-plane $C$ does not intersect the square $B(a,r)$.}
\end{figure}
\par The reader can easily see that in this case imbedding \rf{IN-MN} does not hold provided
$r/s\to\infty$. Cf. Fig. 8 with Fig. 9 below.
\begin{figure}[h!]
\centering{\includegraphics[scale=0.60]{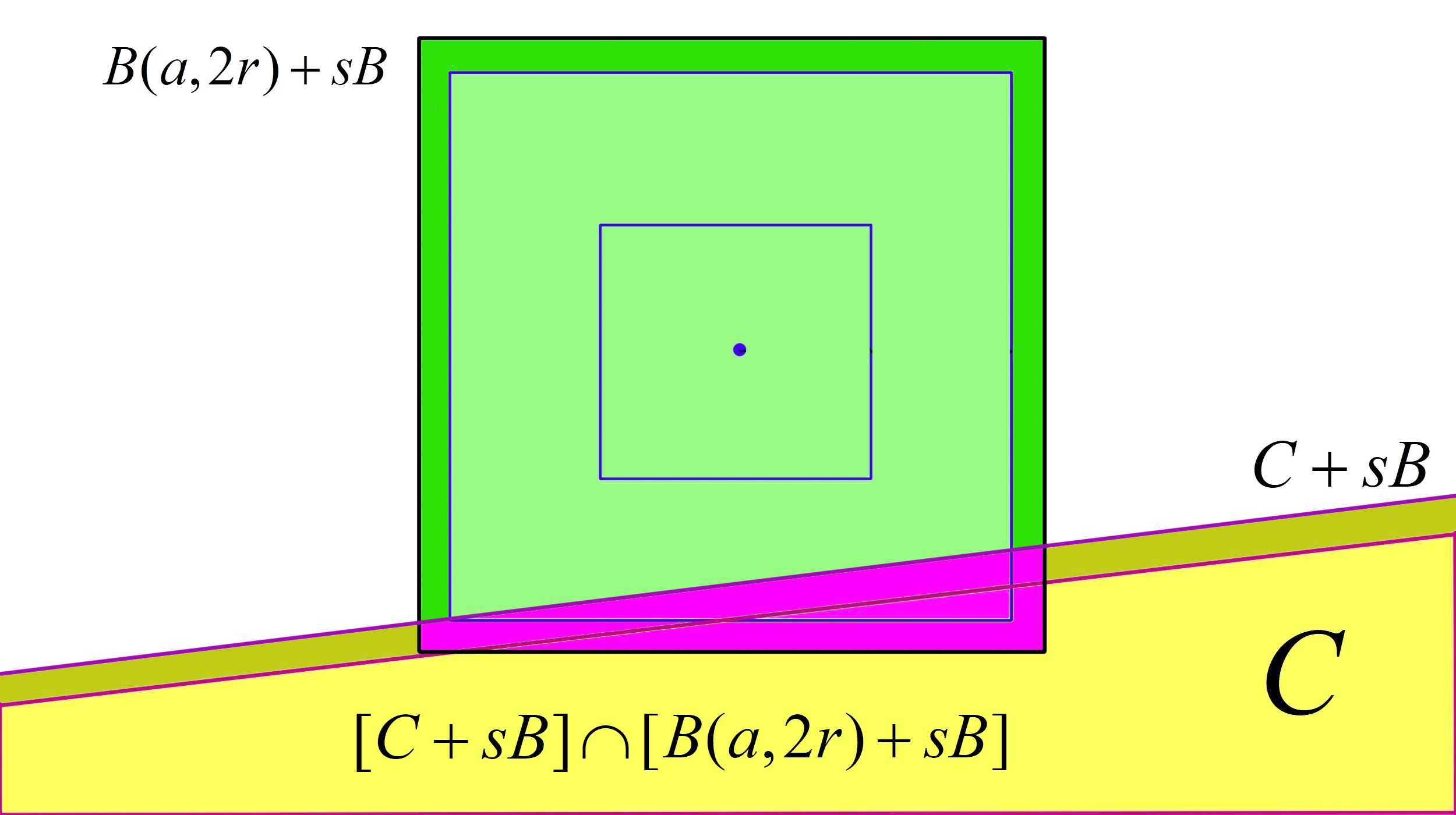}}
\caption{$\diam\{[C+s\BX]\cap [B(a,2r)+sB]\}\gg$
$\diam\{[C\cap B(a,2r)]+7sB\}$.}
\end{figure}
\par The obvious reason for this is the fact that if condition \rf{C-PR} is not met, then the angle between the boundaries of $C$ and $\BX(a,2r)$ at the point of their intersection is very small.
\par However, if $B_X(a,r)\cap C\ne\emp$, than, thanks to \rf{IN-MN}, the sets $A_1=C\cap \BX(a,2r)+7s\BX$ and $A_2=(C+s\BX)\cap (\BX(a,Lr)+s\BX)$ are ``compatible'', i.e., $A_1\supset A_2$.
\smsk
\par Note that in this case the angle between the boundaries of $C$ and $\BX(a,2r)$ at the point of their intersection is not small. Cf. Fig. 10 with Fig. 11 below.
\begin{figure}[h!]
\centering{\includegraphics[scale=0.65]{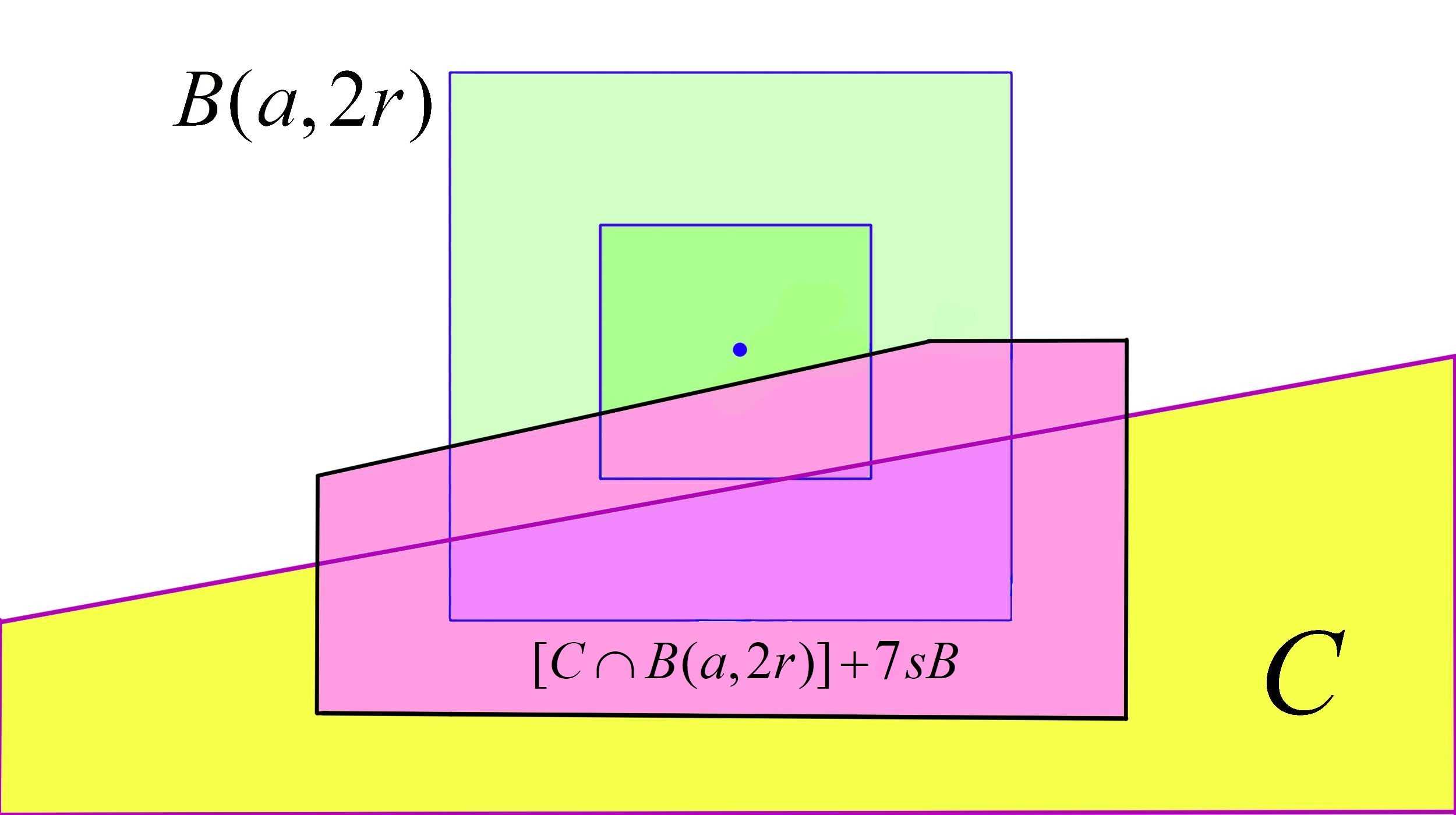}}
\caption{The set $[C\cap B(a,2r)]+7sB$.}
\end{figure}

\begin{figure}[h!]
\centering{\includegraphics[scale=0.60]{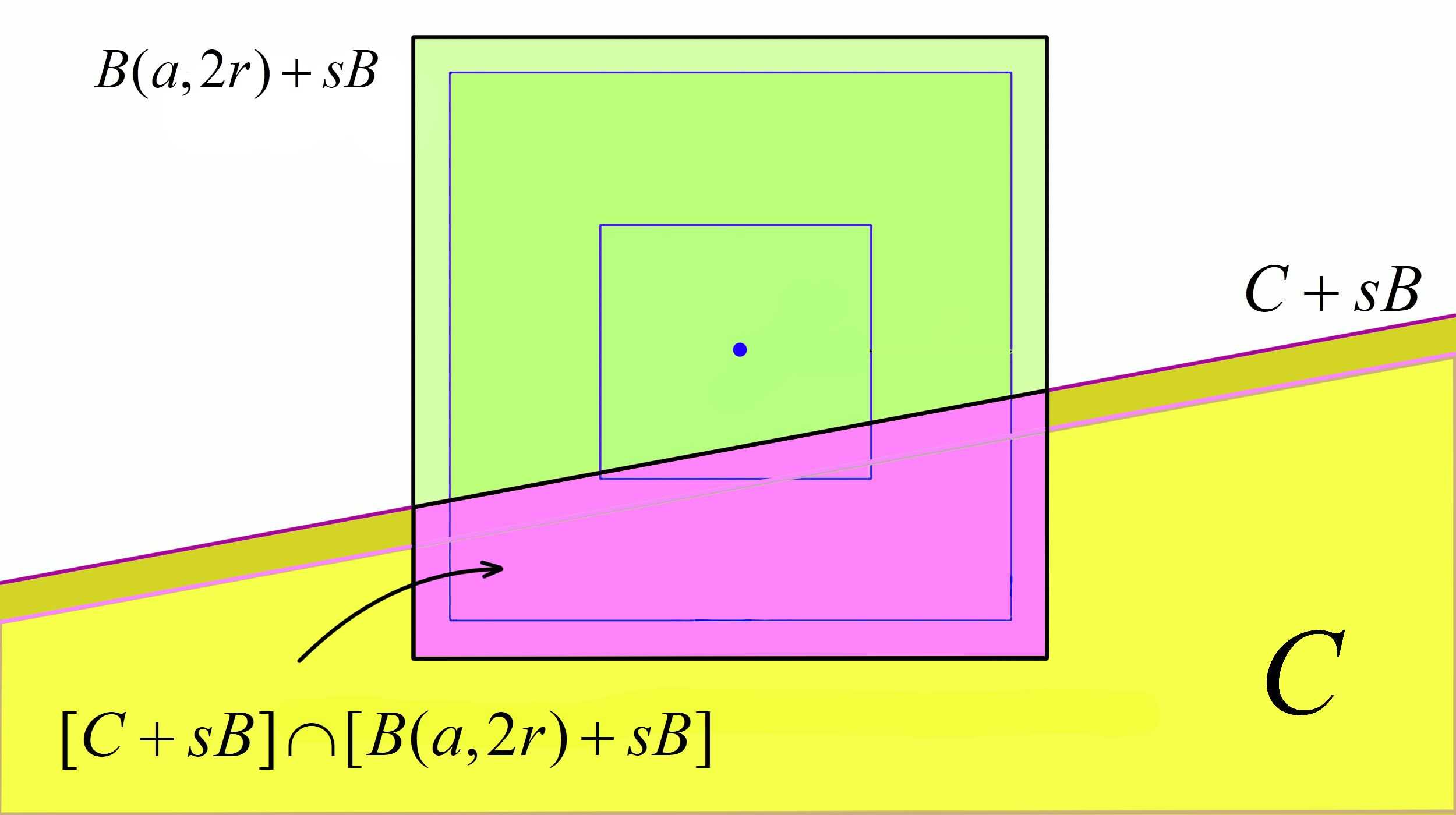}}
\caption{The set $[C+s\BX]\cap [B(a,2r)+sB]$.}
\end{figure}
}
\end{remark}

\msk
\par  Let us recall the classical Helly Intersection Theorem for {\it two dimensional} Banach spaces. It can be formulated as follows:
\begin{theorem}\label{H-TH} Let $\Kc$ be a collection of convex closed subsets of a two dimensional Banach space $\X$. Suppose that $\Kc$ is finite or at least one member of the family $\Kc$ is bounded.
\par If every subfamily of $\Kc$ consisting of at most three elements has a common point then there exists a point common to all of the family $\Kc$.
\end{theorem}
\par We conclude this section by stating and proving one
more result which will be one of the main tools for proving our main theorems.
\begin{proposition}\label{P-F3} Let $\X$ be a two dimensional Banach space. Let $C,C_1,C_2\subset \X$ be convex subsets, and let $r>0$. Suppose that
\begin{align}\lbl{A-PT}
C_1\cbg C_2\cbg(C+r\BX)\ne\emp.
\end{align}
\par Then for every $L>1$ and every $\ve>0$ the following inclusion
\begin{align}
&[\{(C_1\cbg C_2)+L\,r\BX\}\cbg C]+\theta(L)\,\ve\BX
\supset\nn\\
&[(C_1\cbg C_2)+(L\,r+\ve) \BX]
\cbg
[\{(C_1+r\BX)\cbg C\}+\ve \BX]
\cbg
[\{(C_2+r\BX)\cbg C\}+\ve \BX]\nn
\end{align}
holds. The function $\theta$ in the above inclusion
is as defined in Theorem \reff{N-S}. I.e., for an arbitrary Banach space $X$ that inclusion holds for $\theta(L)=(3L+1)/(L-1)$, and if $X$ is Euclidean it also holds for $\theta(L)=1+2L/\sqrt{L^{2}-1}$.
\end{proposition}
\par {\it Proof.} Let
$$
A=
[(C_1\cbg C_2)+(L\,r+\ve) \BX]
\cbg
[\{(C_1+r\BX)\cbg C\}+\ve \BX]
\cbg
[\{(C_2+r\BX)\cbg C\}+\ve \BX]
$$
and let $a\in A$. Prove that
\begin{align}\lbl{A-FN1}
a\in
[\{(C_1\cbg C_2)+L\,r\BX\}\cbg C]+\theta(L)\,\ve\BX.
\end{align}
\par Let us show that if
\begin{align}\lbl{N-EM1}
C_1\cbg C_2\cbg(C+r\BX)\cbg \BX(a,Lr+\ve)\ne\emp
\end{align}
then \rf{A-FN1} holds. Indeed, \rf{N-EM1} provides the existence of a point $x\in \X$ such that
\begin{align}\lbl{N-6}
x\in C_1\cbg C_2\cbg (C+r\BX)\cbg \BX(a,Lr+\ve)\,.
\end{align}
\par In particular, $x\in C+r\BX$ so that
$\BX(x,r)\cbg C\ne\emp$ proving that condition \rf{C-PR} of Theorem \reff{N-S} holds. This theorem tells us that
\begin{align}\lbl{L-STR}
[\BX(x,L\,r)\cbg C]+\theta(L)\,\ve\BX
\supset
(C+\ve \BX)\cbg \BX(x,L\,r+\ve)\,.
\end{align}
\par Since $a\in A$, we have $a\in C+\ve\BX$. From \rf{N-6} we learn that and $a\in \BX(x,Lr+\ve)$. Thus, the point $a$ belongs to the set  $(C+\ve\BX)\cbg \BX(x,Lr+\ve)$. Therefore, by \rf{L-STR},
$$
a\in [\BX(x,Lr)\cbg C]+\theta(L)\,\ve\BX=
[(x+Lr\BX)\cbg C]+\theta(L)\,\ve\BX.
$$
But $x\in C_1\cbg C_2$, see \rf{N-6}, and the required inclusion \rf{A-FN1} follows.
\msk
\par Thus, it remains to prove \rf{N-EM1}. Helly's Theorem \reff{H-TH} tells us that this property holds provided any three sets in the left hand side of \rf{N-EM1} have a common point. Note that, thanks to \rf{A-PT}, this is true for $C_1$, $C_2$ and $C+r\BX$.  Since the point $a\in A$, we have  $a\in [C_1\cbg C_2]+(Lr+\ve)\BX$,
so that $C_1\cbg C_2\cbg \BX(a,Lr+\ve)\ne\emp$.
\par Let us prove that
\begin{align}\lbl{N-5}
C_1\cbg (C+r\BX)\cbg \BX(a,Lr+\ve)\ne\emp\,.
\end{align}
\par The point $a\in A$ so that
$a\in (C_1+r\BX)\cbg C +\ve\BX$. Let $b$ be a point nearest to $a$ on $(C_1+r\BX)\cbg C$, and let $b_1\in C_1$ be a point nearest to $b$ on $C_1$. Prove that
\begin{align}\lbl{N-5A}
b_1\in C_1\cbg (C+r\BX)\cbg \BX(a,Lr+\ve)\,.
\end{align}
\par Indeed, we have $\|a-b\|\le \ve$ and $\|b_1-b\|\le r$. Thus,
\begin{align}\lbl{N-5B}
b_1\in C_1~~~~\text{(by definition)~~ and}
~~~~b_1\in C+r\BX~~~\text{(because $b\in C$)}.
\end{align}
Furthermore,
$$
\|a-b_1\|\le \|a-b\|+\|b-b_1\|\le \ve+r\le \ve+Lr
$$
proving that $b_1\in \BX(a,Lr+\ve)$. Combining this property with \rf{N-5B} we obtain \rf{N-5A} and \rf{N-5}. In a similar way we show that
$C_2\cbg (C+r\BX)\cbg \BX(a,Lr+\ve)\ne\emp$.
\par The proof of the proposition is complete.\bx

\SECT{3. The main theorem for two dimensional Banach spaces.}{3}

\indent
\par In this section we prove Theorem \reff{MAIN-RT}
\par First, let us recall the notion of the Lipschitz extension constant $e(\mfM,\X)$ which we use in the formulation of this theorem.
\begin{definition}\label{LIP-C} {\em Let $\mfM=(\Mc,\rho)$ be a pseudometric space, and let $X$ be a Banach space.
We define the Lipschitz extension constant $e(\mfM,\X)$ of $X$ with respect to $\mfM$ as the infimum of the constants $\lambda>0$ such that for every subset $\Mc'\subset \Mc$, and every Lipschitz mapping $f:\Mc'\to X$, there exists a Lipschitz extension $\tf:\Mc\to X$ of $f$ to all of $\Mc$ such that $\|\tf\|_{\Lip(\Mc,\X)}\le
\lambda \|f\|_{\Lip(\Mc',\X)}$.}
\end{definition}
\begin{remark}\label{EXT-CN} {\em Recall several results about Lipschitz extension constants which we use in this paper. In particular, thanks to the McShane-Whitney extension theorem, $e(\mfM,\R)=1$ for every pseudometric space $\mfM=\MR$. Hence, $e(\mfM,\LTI)=1$ as well.
\par It follows from \cite{Rif} and \cite{CL} that
$e(\mfM,X)\le 4/3$ provided $X$ is an {\it arbitrary two dimensional Banach space}. See also \cite{Bas}. Furthermore, $e(\mfM,X)\le 4/\pi$ whenever $X$ is an arbitrary  {\it two dimensional Euclidean space}. See \cite{Rif} and \cite{Gr}.
\par We also note that, thanks to Kirszbraun's extension theorem \cite{Kir}, $e(\mfM,\X)=1$ provided {\it $\X$ is a Euclidean space, $\Mc$ is a subset of a Euclidean space $E$}, and $\rho$ is the Euclidean metric in $E$.\rbx}
\end{remark}
\par {\it Proof of Theorem \reff{MAIN-RT}.} Let $\mfM=(\Mc,\rho)$ be a pseudometric space, and let $\X$ be a two dimensional Banach space. Let $F:\Mc\to \Kc(X)$ be a set-valued mapping satisfying the hypothesis of Theorem \reff{MAIN-RT}. Recall that, by this hypothesis, {\it for every subset $S\subset\Mc$ with $\#S\le 4$, the restriction $F|_{S}$ of $F$ to $S$ has a $\rho$-Lipschitz selection $f_{S}:S\to \X$ with $\|f_{S}\|_{\Lip((S,\,\rho),\X)}\le 1$.}
\smsk

\par Let $\lambda_1$ and $\lambda_2$ be positive constants satisfying inequalities \rf{GM-FN}, i.e.,  $\lambda_1\ge e(\mfM,X)$ and $\lambda_2\ge 3\lambda_1$. We set
$L=\lambda_2/\lambda_1$. Thus,
the following inequalities
\begin{align}\lbl{L-GE3}
\lambda_1\ge e(\mfM,\X),~~~~~~L\ge 3,
\end{align}
hold. Then we introduce a new pseudometric on $\Mc$ defined by
\begin{align}\lbl{D-AL}
\ds(x,y)=\lambda_1\,\rho(x,y),~~~~~x,y\in \Mc.
\end{align}
\par This definition, Definition \reff{LIP-C}, the above  hypothesis of Theorem \reff{MAIN-RT} and the inequality $\lambda_1\ge e(\mfM,\X)$ imply the following claim.
\begin{claim}\label{CL-D} Let $\tS\subset \Mc$ be a finite set, and let $S\subset\tS$ be a set with $\#S\le 4$.
Then there exists a $\ds$-Lipschitz mapping $\tf_S:\tS\to X$ with $\|\tf_S\|_{\Lip((\tS,\ds),\X)}\le 1$ such that
$\tf_S(x)\in F(x)$ for every $x\in S$.
\end{claim}
%
%
%

\par We introduce set-valued mappings
\begin{align}\lbl{F-1G}
F^{[1]}(x)=
\bigcap_{z\in\Mc}\,
\left[F(z)+\ds(x,z)\,\BX\right],~~~~~x\in\Mc,
\end{align}
and
\begin{align}\lbl{F-2G}
F^{[2]}(x)=\bigcap_{z\in\Mc}\,
\left[F^{[1]}(z)+L\ds(x,z)\,\BX\right],~~~x\in\Mc.
\end{align}
\par Thus, $F^{[1]}$ and $F^{[2]}$ are the first and the second order $(\{1,L\},\ds)$-balanced refinements of $F$ respectively. See Definition \reff{F-IT}.
\smsk
\par Our aim is to show that, if $L$ and $\lambda_1$ satisfy inequality \rf{L-GE3}, then (i) $F^{[2]}(x)\ne\emp$ on $\Mc$, and (ii) the mapping $F^{[2]}$ is $\ds$-Lipschitz with respect to the Hausdorff distance. We prove the statements (i) and (ii) in Proposition \reff{N-EM} and Proposition \reff{HD-G1} respectively.
\par We begin with the property (i). Its proof relies on a series of auxiliary lemmas.
\begin{lemma}\label{H-IN} Let $X$ be a two dimensional Banach space, and let $\Kc\subset\KX$ be a collection of convex compact subsets of $X$ with non-empty intersection.
\par  Given $\tau>0$, let $B=\tau\,\BX$. Then
\begin{align}\lbl{H-SM}
\left(\,\bigcap_{K\in\,\Kc} K\right) +B
=\bigcap_{K,K'\in\,\Kc}\,
\left\{\,\left(\,K \,\cbig\, K'\right)+B\,\right\}.
\end{align}
\end{lemma}
\par {\it Proof.} Obviously, the right hand side of \rf{H-SM} contains its left hand side. Let us prove the converse statement. Fix a point
\begin{align}\lbl{X-IN}
x\in
\bigcap_{K,K'\in\,\Kc}\,
\left\{\,\left(\,K \,\cbig\, K'\right)+B\,\right\}
\end{align}
and prove that $x\in \cbg\{K:K\in\,\Kc\} +B$.
\par Clearly, it is true if
$\BX(x,\tau)\cap\,(\cap\{K:K\in\Kc\})\ne\emp$.
\par Let $S=\Kc\,\cupbig\,\{\BX(x,\tau)\}$. Helly's intersection Theorem \reff{H-TH} tells us that this property holds provided $\cbg\{K:K\in\,S'\}\ne\emp$ for every subfamily $S'\subset S$ consisting of at most three elements. Clearly, this is true if  $\BX(x,\tau)\notin S'$ because there exists a point common to all of the sets from $\Kc$.
\par Suppose that $\BX(x,\tau)\in S'$. Then $S'=\{\BX(x,\tau),K,K'\}$ for some $K,K'\in\Kc$. Thanks to \rf{X-IN}, $x\in\left(\,K \cbig K'\right)+B$ so that
$\BX(x,\tau)\,\cbig\,K\,\cbig\, K'\ne\emp$.
\par Thus, $\BX(x,\tau)\cap\,(\cap\{K:K\in\Kc\})\ne\emp$, and the proof of the lemma is complete.\bx
\smsk
\begin{lemma}\label{G-NE1} For each $x\in\Mc$ the set $F^{[1]}(x)\in\Kc(X)$, i.e., $F^{[1]}(x)$ is a non-empty convex compact subset of $X$. Furthermore, for every $x,z\in\Mc$ we have
$$
F^{[1]}(z)+L\ds(x,z)\BX=\bigcap_{y',y''\in\Mc}
\left\{\left[(F(y')+\ds(z,y')\BX)\cbg (F(y'')+\ds(z,y'')\BX)\right]+L\ds(x,z)\BX\right\}.
$$
\end{lemma}
\par {\it Proof.} Prove that $F^{[1]}(x)\ne\emp$ for every $x\in\Mc$. Indeed, formula \rf{F-1G} and Helly's Theorem \reff{H-TH} tell us that $F^{[1]}(x)\ne\emp$ provided
\begin{align}\lbl{N-3S}
[F(z_1)+\ds(x,z_1)\BX]\,\cbg\, [F(z_2)+\ds(x,z_2)\BX]\,\cbg\,[F(z_3)+\ds(x,z_3)\BX]\ne\emp
\end{align}
for every $z_1,z_2,z_3\in\Mc$.
\par To prove this property, we apply the hypothesis of Theorem \reff{MAIN-RT} to the set $S=\{x, z_1,z_2,z_3\}$. Thanks to this hypothesis, the restriction $F|_S$ has a $\rho$-Lipschitz selection $f_{S}:S\to X$ with $\rho$-Lipschitz seminorm at most $1$. Hence,
$$
\|f_{S}(x)-f_{S}(z_i)\|\le \rho(x,z_i)\le \lambda_1\rho(x,z_i)=\ds(x,z_i)~~~\text{for every}~~~i=1,2,3,
$$
proving that $f_{S}(x)$ belongs to the left hand side of \rf{N-3S}. Thus, \rf{N-3S} holds for arbitrary $z_i\in\Mc$, $i=1,2,3$, so that $F^{[1]}(x)\ne\emp$.
\par Clearly, $F^{[1]}(x)$ is a convex bounded subset of $X$. See \rf{F-1G}.
\par Finally, the second statement of the lemma immediately follows from this property, Lemma \reff{H-IN} and formula \rf{F-1G}. The proof of the lemma is complete.\bx
\begin{lemma}\label{REP-G} For every $x\in\Mc$ the set $F^{[2]}(x)$ admits the following representation:
$$
F^{[2]}(x)=\bigcap_{u,u',u''\in\Mc}
\left\{\left[(F(u')+\ds(u',u)\BX)\cbg (F(u'')+\ds(u'',u)\BX)\right]+L\ds(u,x)\BX\right\}.
$$
\end{lemma}
\par {\it Proof.} The lemma is immediate from \rf{F-2G} and Lemma \reff{G-NE1}.\bx
\par Given $x,u,u',u''\in\Mc$ we set
\begin{align}\lbl{H-D}
T_x(u,u',u'')=
[(F(u')+\ds(u',u)\BX)\cbg (F(u'')+\ds(u'',u)\BX)]+L\ds(u,x)\,\BX
\,.
\end{align}
%
\begin{remark} {\em
Let us recall the standard definitions
of a {\it metric tree}. Let $T=(V,E)$ be a finite (graph theoretic) tree, where $V$ denotes the set of nodes of $T$, and $E$ denotes the set of edges. Suppose we assign a positive number $\Delta(e)$ to each edge $e\in E$. Then for $x,y\in V$ we can define their {\it distance} $\ds_T(x,y)$ to be the sum of $\Delta(e)$ over all the edges $e$ in the ``minimal path'' joining $x$ to $y$ as in Fig. 12.
\begin{figure}[h!]
\hspace*{30mm}
\includegraphics[scale=0.35]{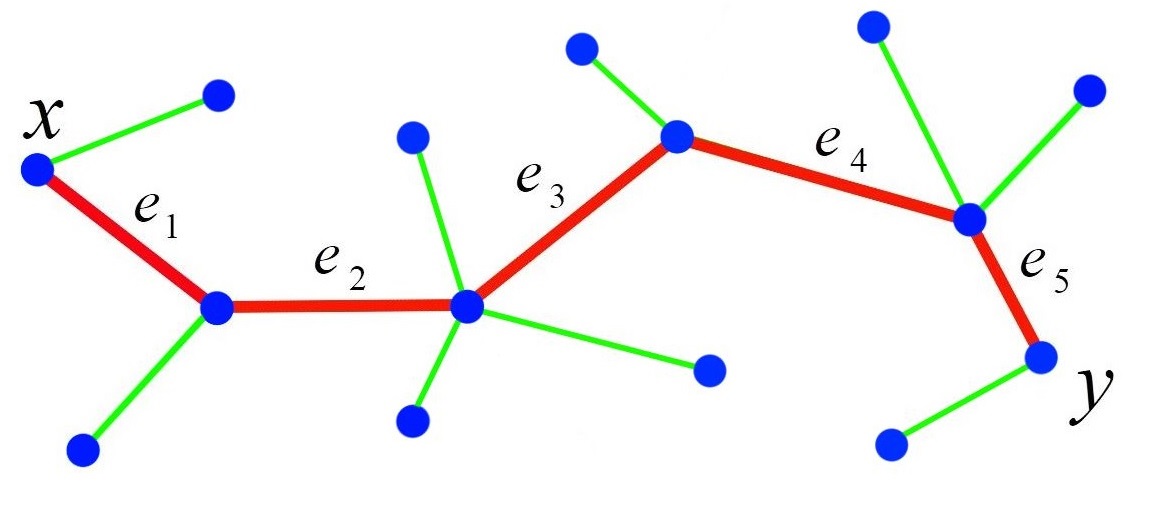}
\caption{A minimal path joining nodes $x$ and $y$ in a tree. In this case, $d(x,y)=\Delta(e_1)+\Delta(e_2)+...+\Delta(e_5).$}
\end{figure}
\par We call $\ds_T$ a {\it tree metric}; $(V,\ds_T)$ is a {\it metric tree}. We say that a mapping $g:V\to X$ {\it agrees} with the metric tree $\mfM_T=(V,\ds_T)$ if the Lipschitz seminorm $\|g\|_{\Lip((V,\ds_T),X)}\le 1$.
\par In these settings, given $x,u,u',u''\in\Mc$, the set
$T_x(u,u',u'')$ defined by \rf{H-D} admits the following equivalent definition: Let $T$ be the metric tree with the family of nodes $V=\{x,u,u',u''\}$ and the set of edges $E=\{e_1=[u,x],e_2=[u',u],e_3=[u'',u]\}$ and the weights $w(e_1)=L\,\ds(u,x)$, $w(e_2)=\ds(u',u)$ and $w(e_3)=\ds(u'',u)$, as it shown in Fig. 13 below.
\begin{figure}[h!]
\centering{\includegraphics[scale=0.47]{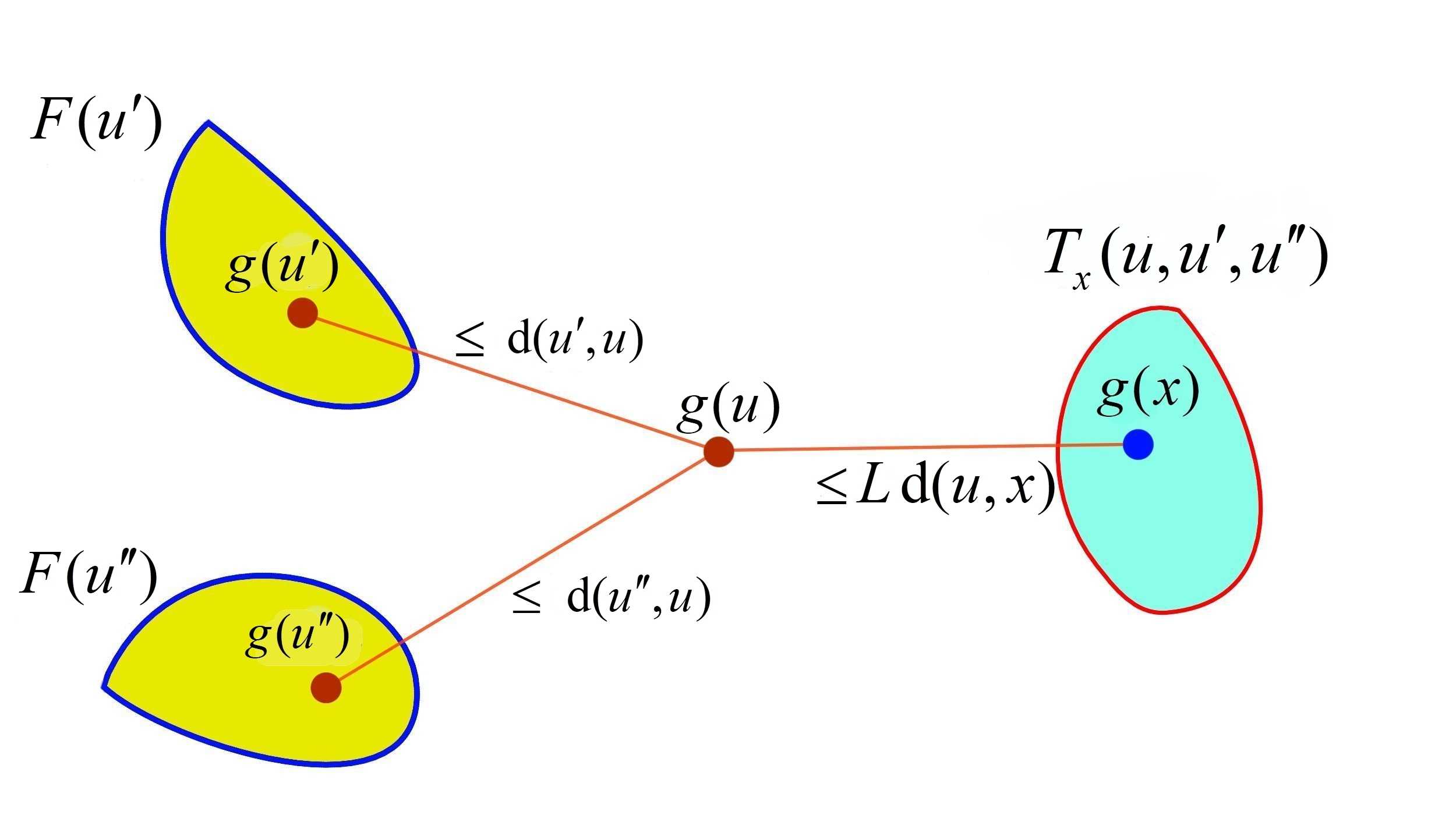}}
\caption{$T_x(u,u',u'')$ is the orbit of $x$ with respect to this diagram.}
\end{figure}
\par Then $T_x(u,u',u'')=\{g(x)\}$ where $g$ runs over all Lipschitz mappings which agree with the metric tree $\mfM_T=(V,\ds_T)$ such that $g(u')\in F(u')$ and $g(u'')\in F(u'')$. We say that the set $T_x(u,u',u'')$ is {\it the orbit of $x$ with respect to the diagram shown in Fig. 13.}
\rbx}
\end{remark}
\par In these settings, Lemma \reff{REP-G} reformulates as follows:
\begin{align}\lbl{G-XP}
F^{[2]}(x)=\bigcap_{u,u',u''\in\Mc} T_x(u,u',u'').
\end{align}
\begin{proposition}\label{N-EM} For every $x\in\Mc$ the set $F^{[2]}(x)$ is non-empty.
\end{proposition}
\par {\it Proof.} Formula \rf{G-XP} and Helly's Theorem \reff{H-TH} tell us that $F^{[2]}(x)\ne\emp$ provided for every choice of elements $u_i,u'_i,u''_i\in\Mc$, $i=1,2,3$, we have
\begin{align}\lbl{H-INT}
T_x(u_1,u'_1,u''_1)\,\cbg\, T_x(u_2,u'_2,u''_2)\,\cbg\,
T_x(u_3,u'_3,u''_3)\ne\emp.
\end{align}
\par We note that \rf{H-INT} holds provided there exists a mapping $g$ which agrees with the diagram shown in Fig. 14 below.
\begin{figure}[h!]
\centering{\includegraphics[scale=0.27]{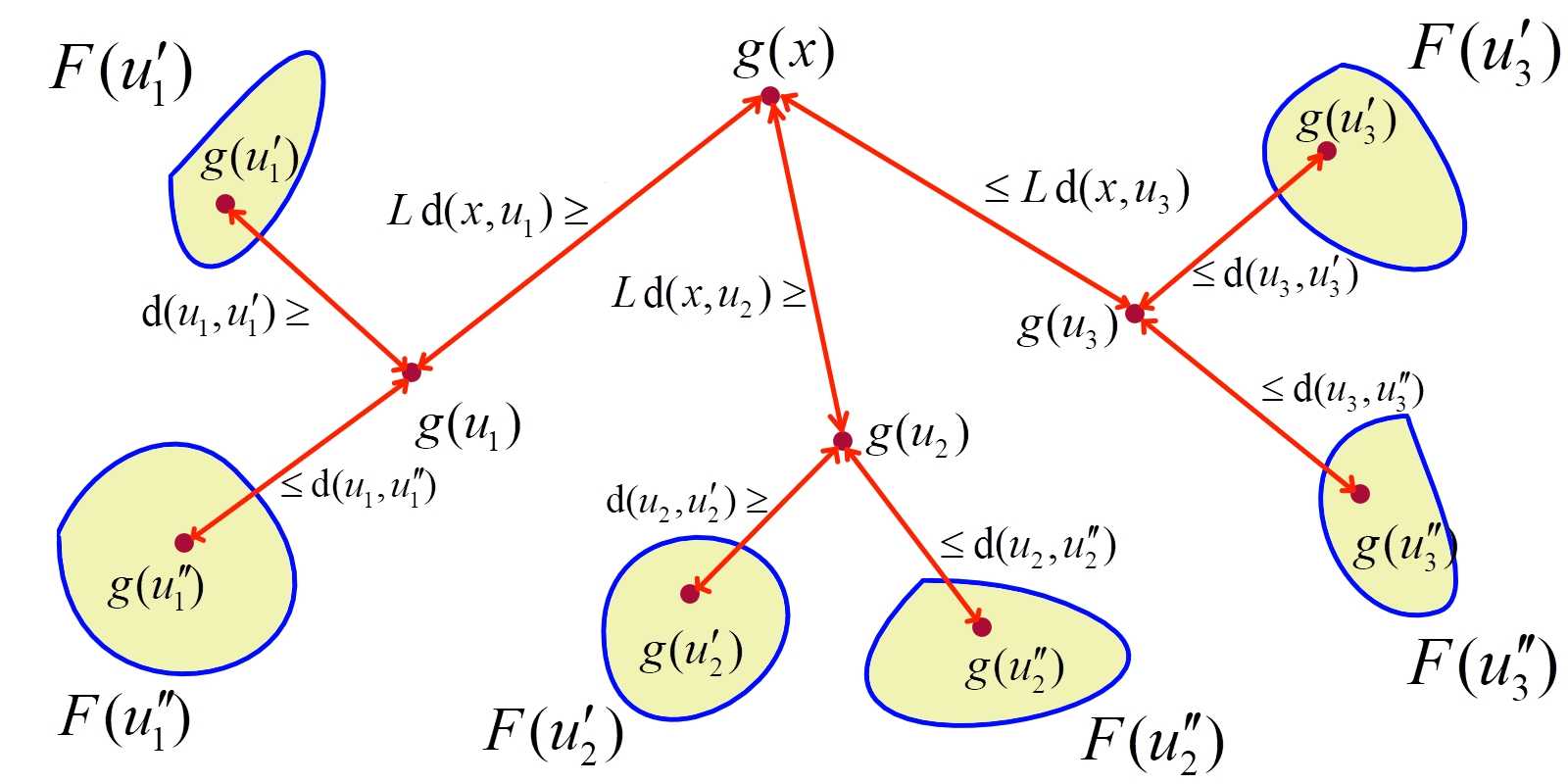}}
\caption{The mapping $g$ agrees with this diagram.}
\end{figure}
\par We set $r_i=\ds(x,u_i)$, $i=1,2,3$. Without loss of generality, we may assume that $r_1\le r_2\le r_3$. For each $i\in\{1,2,3\}$ we also set
\begin{align}\lbl{G-UI}
G(u'_i)=F(u'_i)+\ds(u_i',u_i)\BX~~~~\text{and}~~~~
G(u''_i)=F(u''_i)+\ds(u_i'',u_i)\BX.
\end{align}
\par Note that the existence of a mapping $g$ on the set $\{x,u_i,u'_i,u''_i:i=1,2,3\}$ which agrees with the diagram in Fig. 14 is equivalent to the existence of a mapping $g$ on $\{x,u_1,u_2,u_3\}$ which agrees with the diagram shown in the following Fig. 15.
\begin{figure}[h!]
\centering{\includegraphics[scale=0.27]{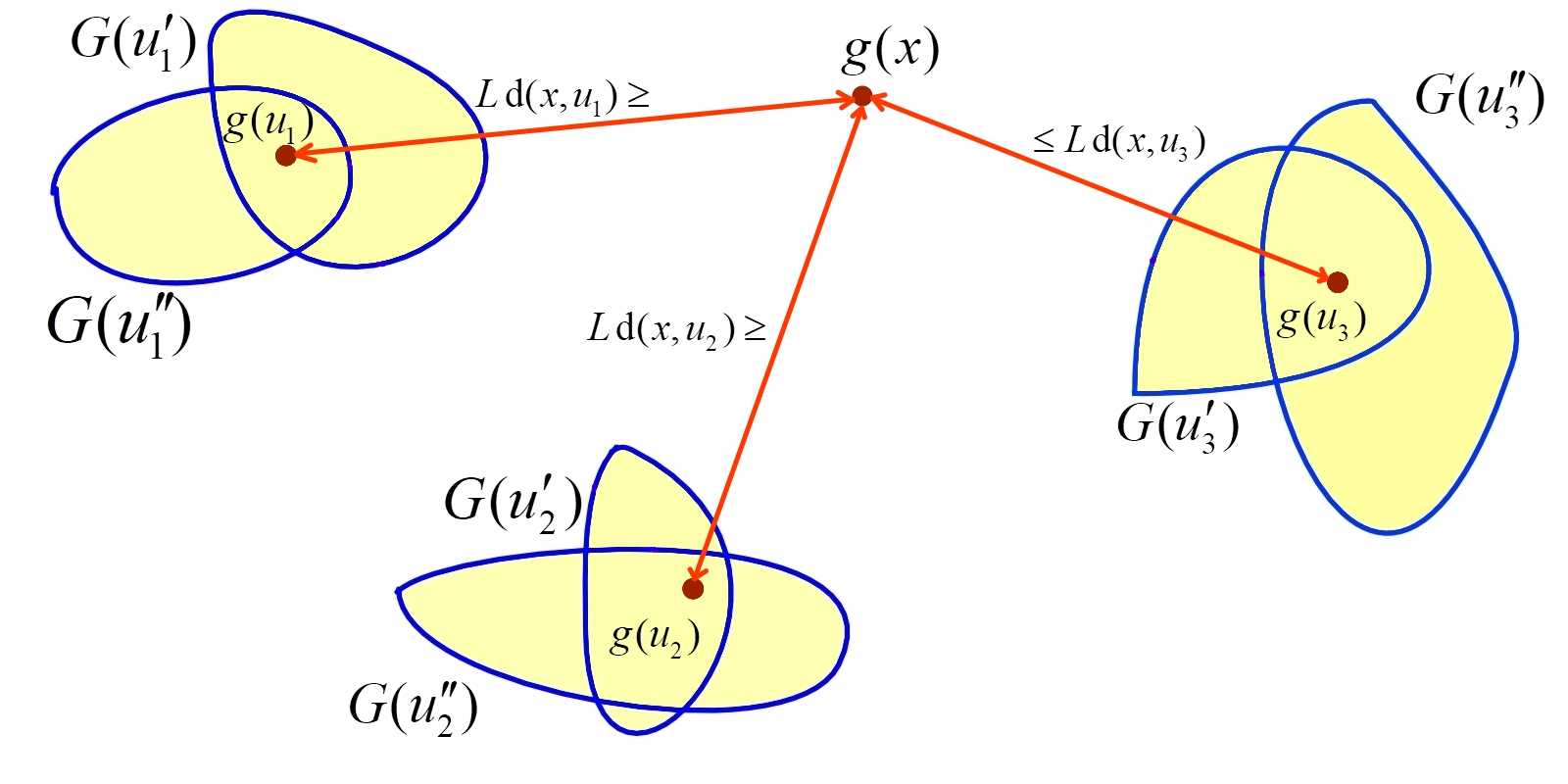}}
\caption{The mapping $g:\{x,u_1,u_2,u_3\}\to X$ agrees with this diagram.}
\end{figure}
\par Let us prove that there exist points $y_i\in\X$, $i=1,2,3$, such that
\begin{align}\lbl{YK-S}
y_i\in G(u'_i)\,\cbg\, G(u''_i)~~~~\text{for every}~~~~
i=1,2,3,
\end{align}
and
\begin{align}\lbl{Y-GA}
\|y_1-y_2\|\le r_1+r_2~~~~\text{and}~~~~
\|y_1-y_3\|\le r_1+2r_2+r_3.~~~~~~~~~\text{See Fig. 16.}
\end{align}
%
%
\begin{figure}[h!]
\centering{\includegraphics[scale=0.27]{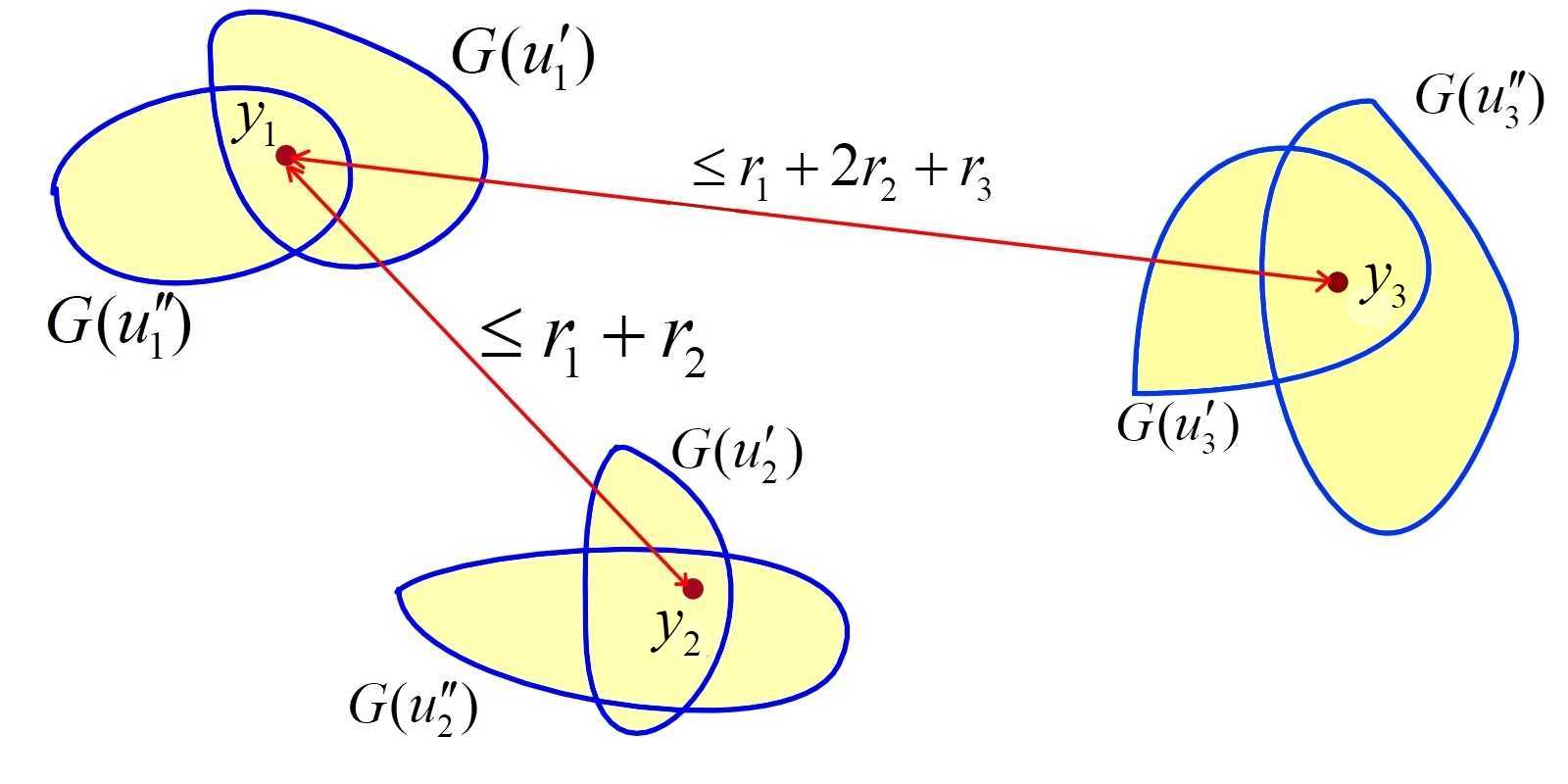}}
\caption{The existence of points $\{y_i:i=1,2,3\}$ implies property \rf{H-INT}. }
\end{figure}
\par Let us see that the existence of the points $y_i$ with these properties implies \rf{H-INT}. Indeed, since $\|y_1-y_3\|\le r_1+2r_2+r_3$, there exists  $z\in[y_1,y_3]$ such that $\|z-y_1\|\le r_1$ and $\|z-y_3\|\le 2r_2+r_3$. (For instance, one can set
$z=y_1+\tau(y_3-y_1)$ with $\tau=r_1/(r_1+2r_2+r_3)$.)
\par Hence,
$$\|y_2-z\|\le \|y_2-y_1\|+\|y_1-z\|\le r_1+r_2+r_1=2r_1+r_2.$$
\par Recall that $r_i=\ds(x,u_i)$ and $r_1\le r_2\le r_3$.
From this and above inequalities, we have
\begin{align}\lbl{Z-YI}
\|z-y_i\|\le 3r_i=3\ds(x,u_i),~~~~i=1,2,3.
\end{align}
\par Let us prove that $z\in T_x(u_i,u'_i,u''_i)$ for each
$i\in\{1,2,3\}$. In fact, we know that $L\ge 3$, see \rf{L-GE3}. Furthermore, by \rf{YK-S}, $y_i\in G(u'_i)\cbg G(u''_i)$, so that, thanks to \rf{G-UI}, \rf{Z-YI} and definition \rf{H-D},
$$
z\in [G(u'_i)\,\cbg\, G(u''_i)]+3\ds(x,u_i)\BX\subset
[G(u'_i)\,\cbg\, G(u''_i)]+L\ds(x,u_i)\BX=T_x(u_i,u'_i,u''_i)
$$
proving \rf{H-INT}.
\msk
\par Thus, our aim is to prove the existence of points $y_i$ satisfying \rf{YK-S} and \rf{Y-GA}. We will do this in three steps.
\msk

\par {\it STEP 1}. We introduce sets $W_i\subset\X$, $i=1,...,4$, defined by
\begin{align}\lbl{W-123-D}
W_1=G(u'_1),~~~~~W_2=G(u''_1),~~~~
W_3=[G(u'_2)\cbg G(u''_2)]+(r_1+r_2)\BX,
\end{align}
and
\begin{align}\lbl{W-4}
W_4=[G(u'_3)\cbg G(u''_3)]+(r_1+2r_2+r_3)\BX.
\end{align}
\par Obviously, there exist the points $y_i$ satisfying \rf{YK-S} and \rf{Y-GA} whenever
$$
W_1\cbg W_2\cbg W_3\cbg W_4\ne \emp.
$$
See Fig. 17. 

\begin{figure}[h!]
\centering{\includegraphics[scale=0.37]{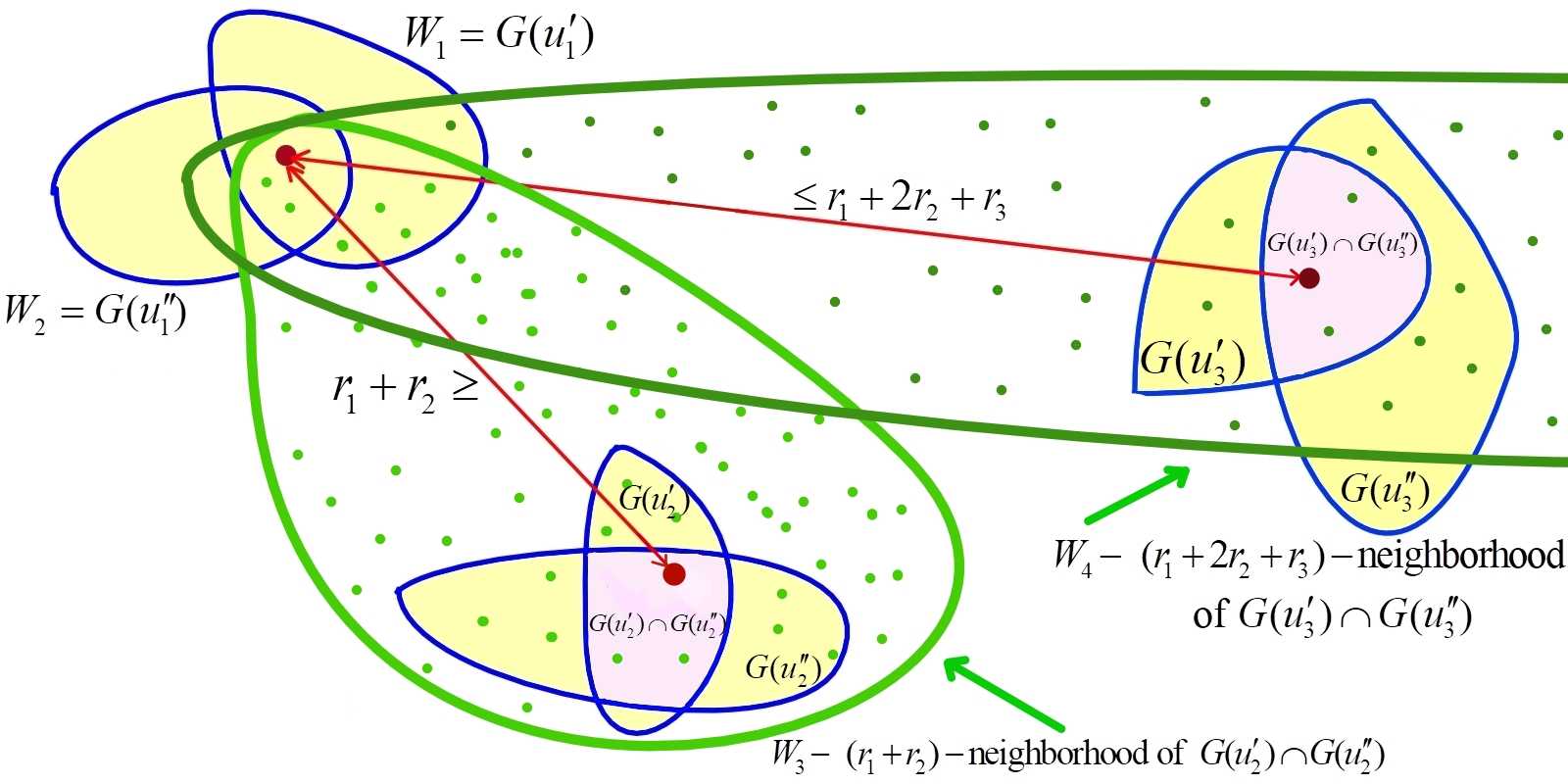}}
\caption{The sets $W_1$, $W_2$, $W_3$ and $W_4$.}
\end{figure}
\par By Helly's Theorem \reff{H-TH}, this property holds provided any three members of the family of sets $\{W_1,W_2,W_3,W_4\}$ have a common point.
\par {\it STEP 2.} Prove that $W_1\cbg W_3\cbg W_4\ne\emp$.
To see this, we set
\begin{align}\lbl{V-123}
V_1=G(u'_1)+(r_1+r_2)\BX,~~~~~V_2=G(u'_2),~~~~V_3=G(u''_2),
\end{align}
and
\begin{align}\lbl{V-4}
V_4=[G(u'_3)\cbg G(u''_3)]+(r_2+r_3)\BX.
\end{align}
\par Let us show that if
\begin{align}\lbl{VI-INT}
V_1\cbg V_2\cbg V_3\cbg V_4\ne \emp,
\end{align}
then $W_1\cbg W_3\cbg W_4$ is non-empty as well. Indeed, this property, definitions \rf{V-123} and \rf{V-4} imply the existence of points $z_1\in G(u'_1)$, $z_2\in G(u'_2)\cbg G(u''_2)$, $z_3\in G(u'_3)\cbg G(u''_3)$ such that $\|z_1-z_2\|\le r_1+r_2$ and $\|z_2-z_3\|\le r_2+r_3$. Hence,
$$
\|z_1-z_3\|\le\|z_1-z_2\|+\|z_2-z_3\|\le (r_1+r_2)+(r_2+r_3)=r_1+2r_2+r_3.
$$
\par Thus, thanks to \rf{W-123-D} and \rf{W-4}, the point
$z_1$ belongs to $W_1\cbg W_3\cbg W_4$ proving that this set is non-empty.
smsk
\par Let us prove \rf{VI-INT}. Helly's Theorem \reff{H-TH}
tells us that \rf{VI-INT} holds whenever every three members of the family $\Vc=\{V_1,V_2,V_3,V_4\}$ have a common point.
\par Let us prove this property. First, let us show that
\begin{align}\lbl{V-124}
V_1\cbg V_2\cbg V_4\ne \emp.
\end{align}
\par Let $S=\{u'_1,u'_2,u'_3,u''_3\}$ and let $\tS=\{u_1',u_2,u_2',u_3,u_3',u_3''\}$. Note that \rf{V-124} holds if and only if there exists a mapping $\tg$ on $\tS$ which agrees with the diagram in Fig. 18 below.
\begin{figure}[h!]
\centering{\includegraphics[scale=0.35]{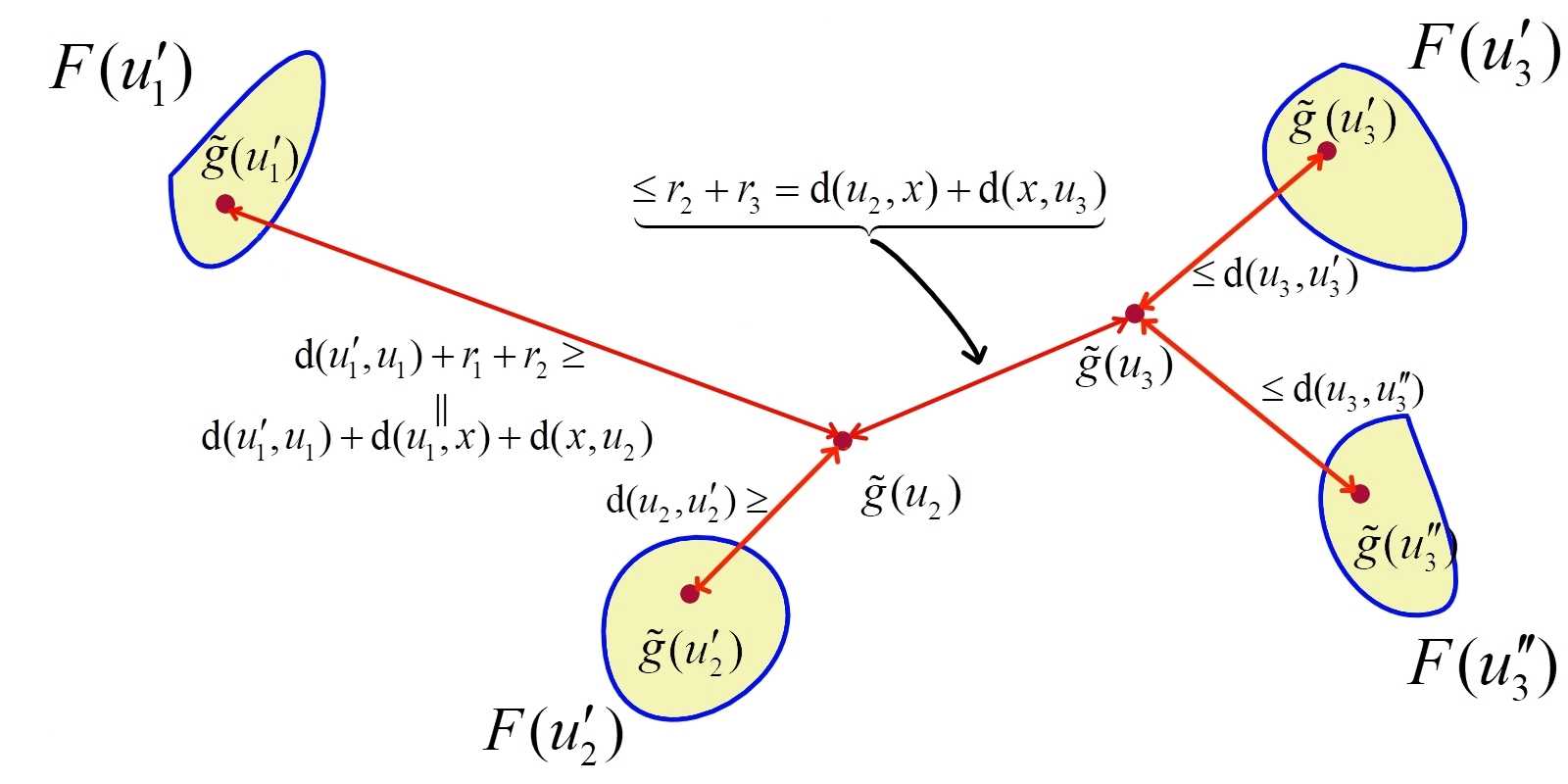}}
\caption{The mapping $\tg:\tS\to X$ agrees with this diagram.}
\end{figure}
\par Note that, thanks to Claim \reff{CL-D}, there exists a $\ds$-Lipschitz mapping $\tf_{S}:\tS\to\X$ with $\|\tf_{S}\|_{\Lip((\tS,\,\ds),\X)}\le 1$ such that
$$
\tf_{S}(u'_1)\in F(u'_1),~~~~
\tf_{S}(u'_2)\in F(u'_2),~~~~
\tf_{S}(u'_3)\in F(u'_3),~~~~\text{and}~~~~
\tf_{S}(u''_3)\in F(u''_3).
$$
\par Prove that
\begin{align}\lbl{F1-AV}
\tf_{S}(u_2)\in V_1\cbg V_2\cbg V_4.
\end{align}
\par Indeed, $\tf_{S}(u'_2)\in F(u'_2)$ and $\|\tf_{S}(u'_2)-\tf_{S}(u_2)\|\le \ds(u'_2,u_2)$, so that
$$
\tf_{S}(u_2)\in F(u'_2)+\ds(u'_2,u_2)\BX=G(u'_2)=V_2.
$$
\par In the same way we prove that $\tf_S(u_1)\in G(u'_1)$.
\par Note that $\|\tf_S(u_1)-\tf_S(u_2)\|\le\ds(u_1,u_2)$ so that $\tf_S(u_2)\in G(u'_1)+\ds(u_1,u_2)\BX$.
By the triangle inequality,
$$
\ds(u_1,u_2)\le \ds(u_1,x)+\ds(x,u_2)=r_1+r_2
$$
proving that
$$
\tf_S(u_2)\in G(u'_1)+(r_1+r_2)\BX=V_1.
$$
\par It remains to show that $\tf_S(u_2)$ belongs to $V_4$. We know that $\tf_S(u'_3)\in F(u'_3)$ and $\tf_S(u''_3)\in F(u''_3)$. Furthermore,
$$
\|\tf_S(u_3)-\tf_S(u'_3)\|\le\ds(u_3,u'_3),~~~~
\|\tf_S(u_3)-\tf_S(u''_3)\|\le\ds(u_3,u''_3).
$$
Hence,
$$
\tf_S(u_3)\in [F(u'_3)+\ds(u'_3,u_3)\BX]\cbg
[F(u''_3)+\ds(u''_3,u_3)\BX]=G(u'_3)\cbg G(u''_3).
$$
Furthermore, $\|\tf_S(u_2)-\tf_S(u_3)\|\le\ds(u_2,u_3)$. These properties of $\tf_S(u_3)$ and the triangle inequality
$$
\ds(u_2,u_3)\le \ds(u_2,x)+\ds(x,u_3)=r_2+r_3
$$
imply the following:
$$
\tf_S(u_2)\in [G(u'_3)\cbg G(u''_3)]+\ds(u_2,u_3)\BX
\subset [G(u'_3)\cbg G(u''_3)]+(r_2+r_3)\BX=V_4.
$$
\par Thus, $\tf_S(u_2)\in V_1\cbg V_2\cbg V_4$ proving \rf{V-124}.
\smsk
\par In the same fashion we show that $V_1\cbg V_3\cbg V_4\ne\emp$.
%

\par Next, let us prove that
\begin{align}\lbl{V-234}
V_2\cbg V_3\cbg V_4=G(u'_2)\cbg G(u''_2)
\cbg\{[G(u'_3)\cbg G(u''_3)]+(r_2+r_3)\BX\}\ne \emp.
\end{align}
\par Following the scheme of the proof of \rf{V-124}, we put
$$
S=\{u'_2,u''_2,u'_3,u''_3\}~~~~\text{and}~~~~ \tS=\{u_2,u'_2,u''_2,u_3,u'_3,u''_3\}.
$$
\par We note that \rf{V-234} holds provided there exists a mapping $\tg:\tS\to X$ which agrees with the diagram on Fig. 19 below.

\hspace*{-4mm}
\begin{figure}[h!]
\centering{\includegraphics[scale=0.36]{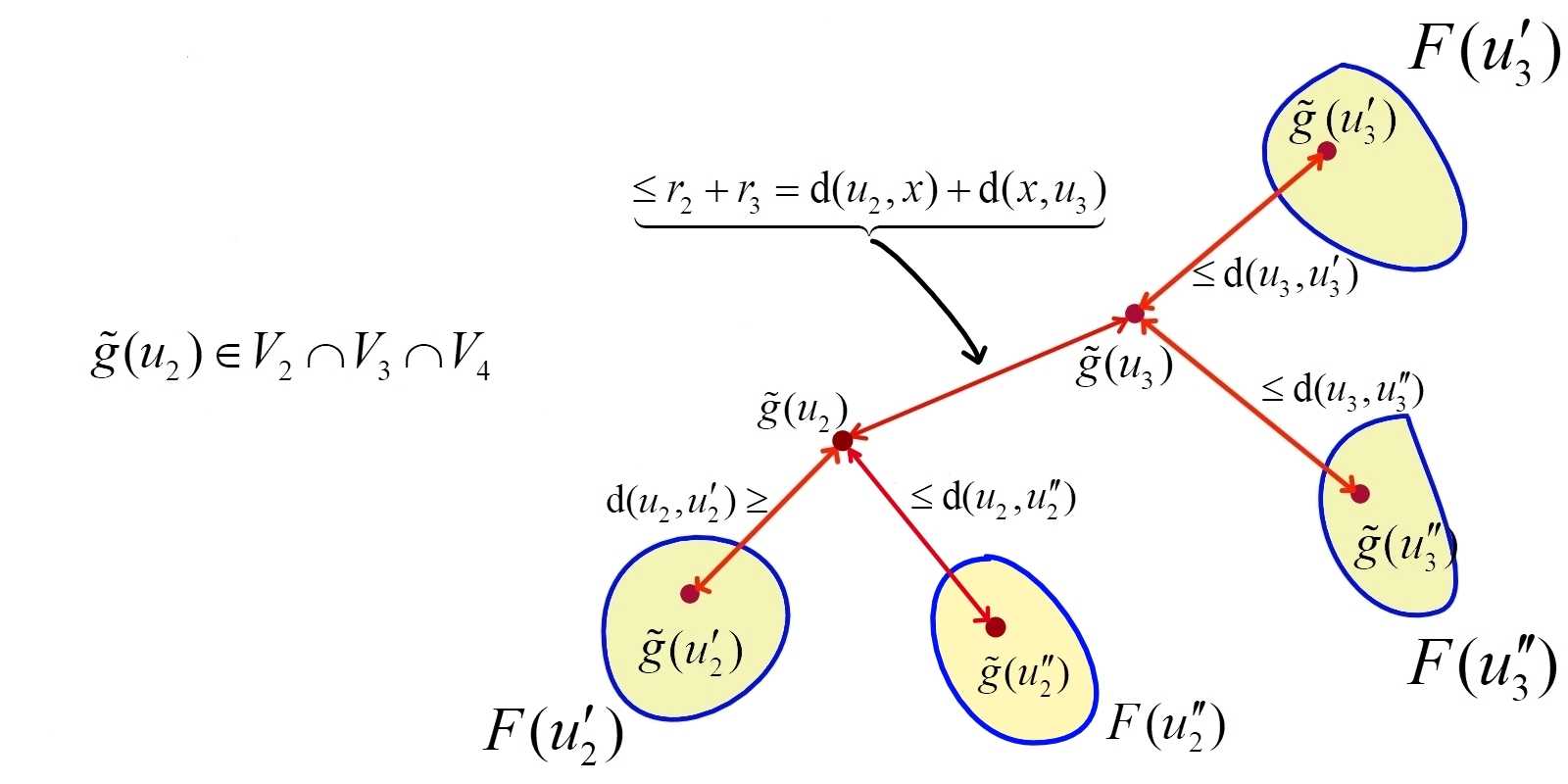}}
\caption{The existence of $\tg:\tS\to X$ which agrees with this diagram implies \rf{V-234}.}
\end{figure}
\par Claim \reff{CL-D} tells us that there exists a $\ds$-Lipschitz mapping $\tf_S:\tS\to\X$ with Lipschitz seminorm $\|\tf_S\|_{\Lip((\tS,\,\ds),\X)}\le 1$ such that $\tf_S(u'_i)\in F(u'_i)$ and  $\tf_S(u''_i)\in F(u''_i)$, $i=2,3$.
\par Then, following the scheme of the proof of \rf{F1-AV}, we show that $\tf_S(u_2)$ belongs to $V_2\cbg V_3\cbg V_4$ proving the required property \rf{V-234}.
\msk
\par Finally, following the same approach, we prove that
\begin{align}\lbl{V-F}
V_1\cbg V_2\cbg V_3=[G(u'_1)+(r_1+r_2)\BX]\cbg G(u'_2)\cbg G(u''_2)\ne \emp.
\end{align}
\smsk
\par More specifically, we set $S=\{u'_1,u_2,u'_2,u''_2\}$. Then, thanks to the hypothesis of Theorem \reff{MAIN-RT}, there exists a $\ds$-Lipschitz mapping $\tf_S:S\to\X$ with $\|\tf_S\|_{\Lip((S,\,\ds),\X)}\le 1$ such that
$$
\tf_S(u'_1)\in F(u'_1),~~\tf_S(u'_2)\in F(u'_2)~~~
\text{and}~~\tf_S(u''_2)\in F(u''_2).
$$
Following the proof of \rf{F1-AV}, we show that $\tf_S(u_2)\in V_1\cbg V_2\cbg V_3$ completing the proof of \rf{V-F}.
\par Thus, \rf{VI-INT} is proven, so that $W_1\cbg W_3\cbg W_4\ne\emp$.
\msk

\par {\it STEP 3}. First, using a similar approach, we show that $W_2\cbg W_3\cbg W_4\ne\emp$.
\par Next, we prove that
\begin{align}\lbl{W-124}
W_1\cbg W_2\cbg W_4=
G(u'_1)\cbg G(u''_1)\cbg
\{[G(u'_3)\cbg G(u''_3)]+(r_1+2r_2+r_3)\BX\}
\ne\emp.
\end{align}
\par To  see this, we set $S=\{u'_1,u''_1,u'_3,u''_3\}$ and $\tS=\{u_1,u'_1,u''_1,u_3,u'_3,u''_3\}$. Note that \rf{W-124} holds whenever the distance between the set $G(u_1')\cap G(u_1'')$ and the set $G(u_3')\cap G(u_3'')$ is bounded by $r_1+r_2+r_3$. See Fig. 20.

\hspace*{-4mm}
\begin{figure}[h!]
\centering{\includegraphics[scale=0.36]{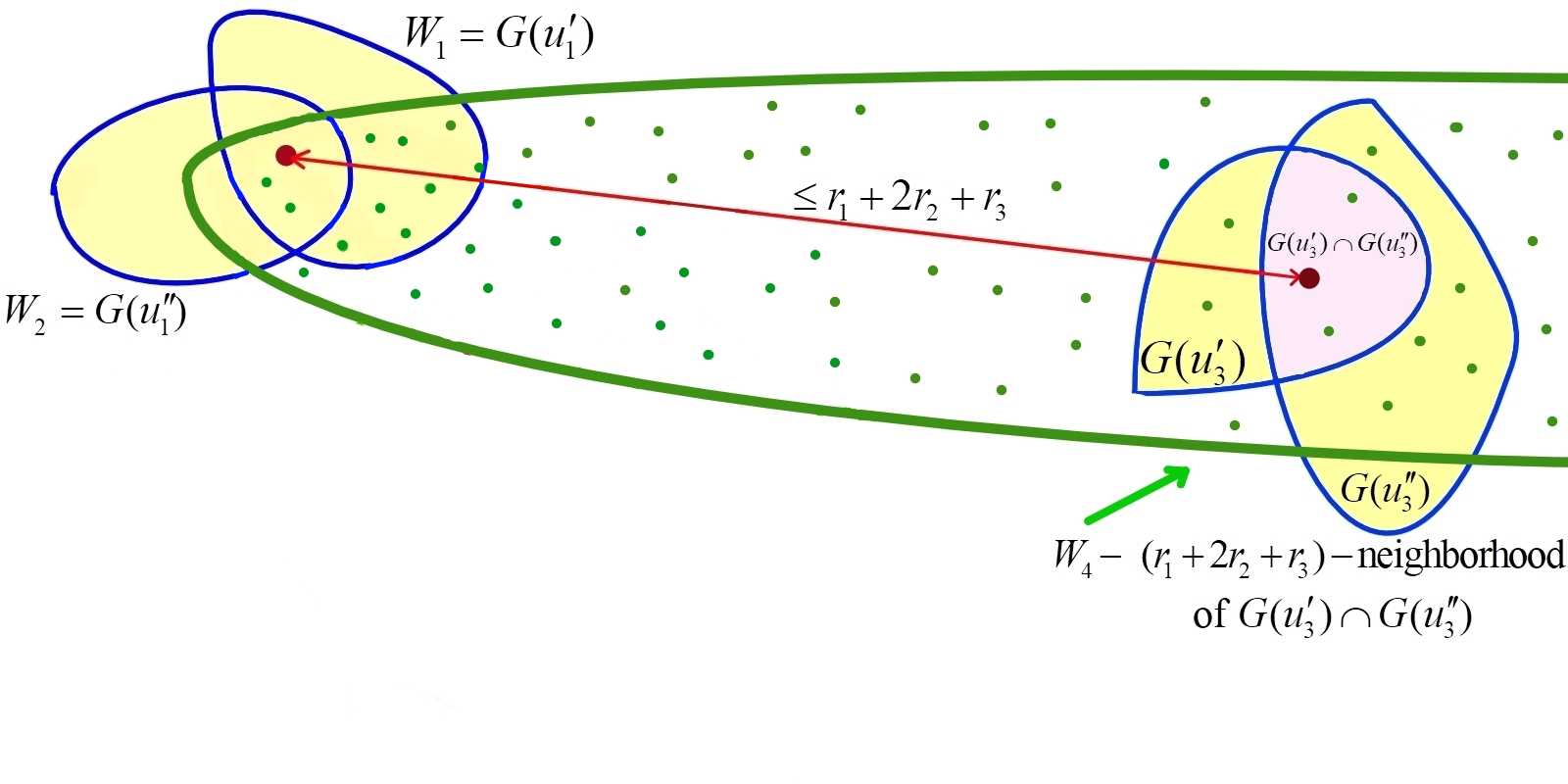}}
\vspace*{-15mm}
\caption{The sets $W_1$, $W_2$ and $W_4$.}
\end{figure}
\par We note that Claim \reff{CL-D} tells us that there exists a $\ds$-Lipschitz mapping $\tf_S:\tS\to\X$ with $\|\tf_S\|_{\Lip((\tS,\,\ds),\X)}\le 1$ such that $\tf_S(u'_i)\in F(u'_i)$ and  $\tf_S(u''_i)\in F(u''_i)$, $i=1,3$. The reader can easily see that the point  $\tf_S(u_1)$ belongs to $W_1\cbg W_2\cbg W_4$ so that \rf{W-124} holds.
\msk
\par In a similar way we prove that
\begin{align}\lbl{W-123}
W_1\cbg W_2\cbg W_3=
G(u'_1)\cbg G(u''_1)\cbg
\{[G(u'_2)\cbg G(u''_2)]+(r_1+r_2)\BX\}
\ne\emp.
\end{align}
\par More specifically, we set $S=\{u'_1,u''_1,u'_2,u''_2\}$, $\tS=\{u_1,u'_1,u''_1,u_2,u'_2,u''_2\}$, and apply Claim \reff{CL-D} to $S$ and $\tS$. Thanks to this claim,
there exists a $\ds$-Lipschitz mapping $\tf_S$ on $\tS$ with $\|\tf_S\|_{\Lip((\tS,\,\ds),\X)}\le 1$ such that $\tf_S(u'_i)\in F(u'_i)$ and  $\tf_S(u''_i)\in F(u''_i)$, $i=1,2$. One can readily see that the point $\tf_S(u_1)\in W_1\cbg W_2\cbg W_3$ proving the required property \rf{W-123}.
\smsk
\par The proof of the proposition is complete.\bx
\msk

\par We turn to the proof of inequality \rf{HD-RT}.

\par Note that, thanks to formula \rf{G-XP}, for every $x,y\in\Mc$ we have
\begin{align}\lbl{G-Y}
F^{[2]}(x)=\bigcap_{u,u',u''\in\Mc} T_x(u,u',u'')~~~\text{and}~~~
F^{[2]}(y)=\bigcap_{u,u',u''\in\Mc} T_y(u,u',u'').
\end{align}

\begin{lemma}\label{L-AW} For every $\tau>0$ and every $x\in\Mc$ the following representation
\begin{align}\lbl{GX-1}
F^{[2]}(x)+\tau\,\BX=
\bigcap\,\,
\left\{\,\left[T_x(u,u',u'')\cbg T_x(v,v',v'')\right]
+\tau\,\BX\,\right\}
\end{align}
holds. Here the first intersection in the right hand side of \rf{GX-1} is taken over all elements
$u,u',u'',v,v',v''\in\Mc$.
\end{lemma}
\par {\it Proof.} The lemma is immediate from representation \rf{G-Y}, Lemma \reff{H-IN} and Proposition \reff{N-EM}.\bx

\begin{proposition}\label{HD-G1} For every $x,y\in\Mc$ the following inequality
\begin{align}\lbl{HD-KJ}
\dhf(F^{[2]}(x),F^{[2]}(y))\le \gamma_0(L)\ds(x,y)
\end{align}
holds. Here $\gamma_0(L)=L\,\theta(L)^2$ where $\theta(L)$ is the constant from Theorem \reff{N-S}.
\end{proposition}
\par {\it Proof.} Let $x,y\in\Mc$ and let  $\tau=\gamma_0(L)\ds(x,y)$. Prove that
\begin{align}\lbl{FIL}
F^{[2]}(x)+\gamma_0(L)\ds(x,y)\BX=F^{[2]}(x)+\tau\BX
\supset F^{[2]}(y)\,.
\end{align}
\par Lemma \reff{L-AW} tells us that this inclusion holds provided
\begin{align}\lbl{A-1}
A=[T_x(u,u',u'')\cbg T_x(v,v',v'')]
+\tau\,\BX\,\supset F^{[2]}(y)
\end{align}
for arbitrary $u,u',u'',v,v',v''\in\Mc$. We note that the set $A$ is the orbit of $y$ with respect to the diagram shown in Fig. 21. (I.e., the set $A$ consists of all points $a=g(y)$ where $g$ runs over all mappings which agree with the diagram in Fig. 21.)
\begin{figure}[h!]
\centering{\includegraphics[scale=0.6]{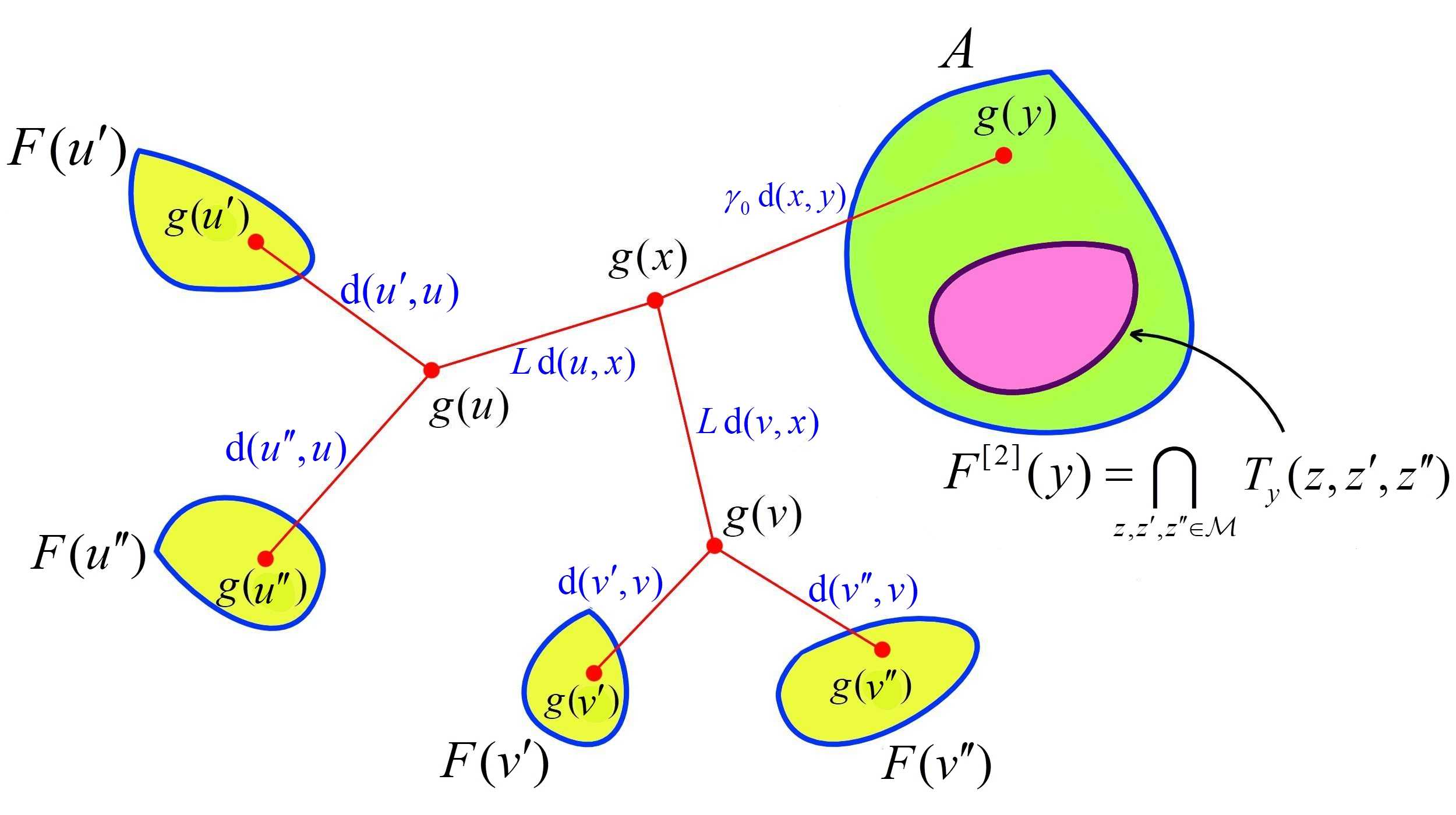}}
\caption{The set $A$ is the orbit of $y$ with respect to this diagram.}
\end{figure}
\par In fact, we prove a stronger inclusion than \rf{A-1}, namely that
\be
A&\supset& T_y(u,u',u'')\bigcap T_y(v,v',v'')\bigcap T_y(x,u',v')\bigcap\nn\\
\nn\\
&& T_y(x,u',v'')\bigcap T_y(x,u'',v')
\bigcap T_y(x,u'',v'').\nn
\ee
\par To prove the above inclusion, we introduce the following sets:
\begin{align}\lbl{CI-S}
C_1=F(u')+\ds(u',u)\BX,~~~C_2=F(u'')+\ds(u'',u)\BX,~~~
C=T_x(v,v',v'')\,.
\end{align}
Let
\begin{align}\lbl{DL}
\ve=L\,\theta(L)\ds(x,y)~~~~~\text{and}~~~~~r=\ds(x,u)\,.
\end{align}
\par Then $\tau=\gamma_0(L)\ds(x,y)=\theta(L)\,\ve$,
and
$$
A=[T_x(u,u',u'')\cbg T_x(v,v',v'')]+\tau\,\BX
=[(C_1\cbg C_2)+Lr\BX]\cbg C+\theta(L)\,\ve\,\BX.
$$
%
See Fig. 22.

\hspace*{20mm}
\begin{figure}[h!]
\centering{\includegraphics[scale=0.6]{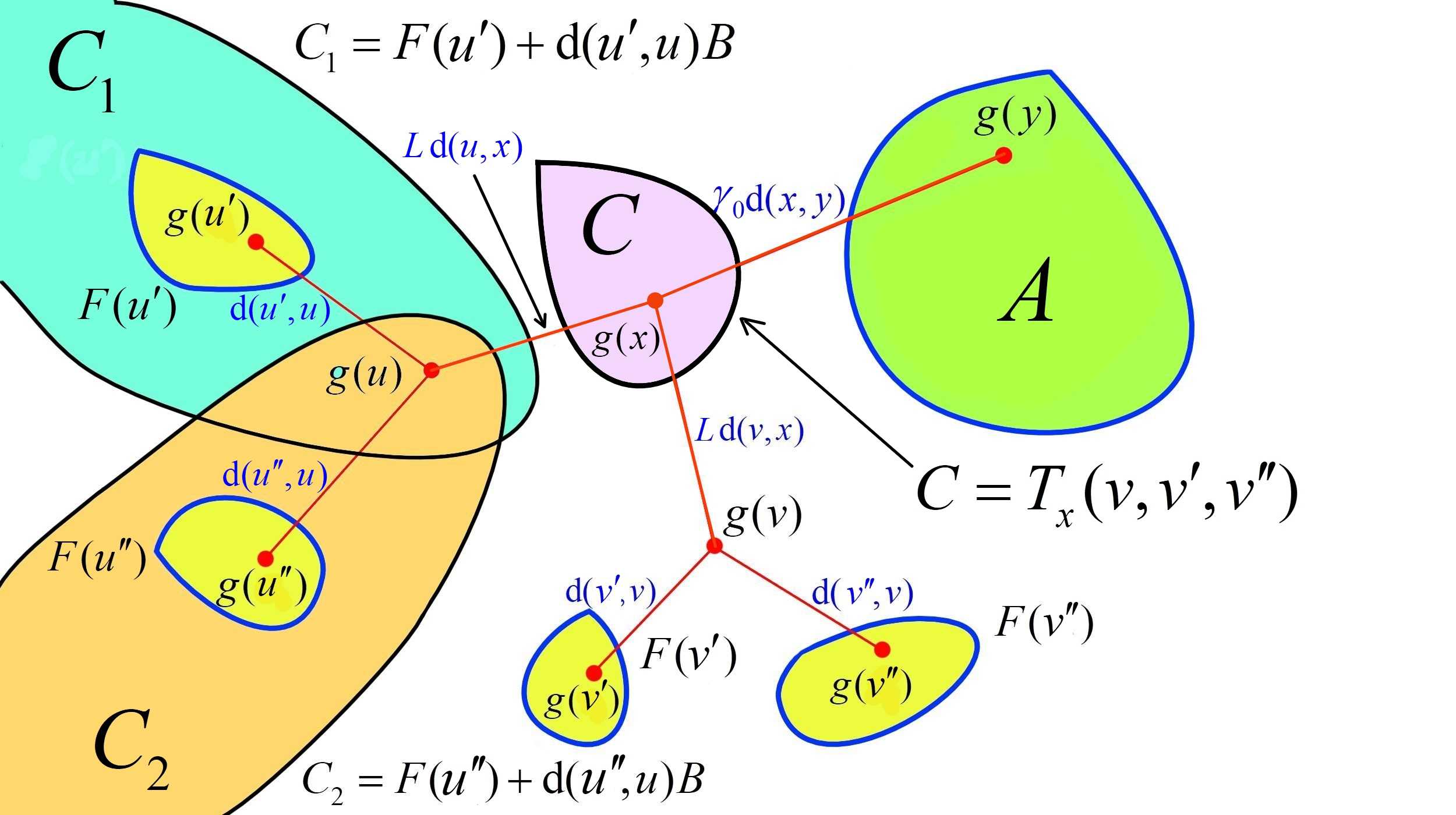}}
\caption{The sets $C$, $C_1$ and $C_2$.}
\end{figure}

\par We want to apply Proposition \reff{P-F3} to the set $A$. To do this, we have to verify condition \rf{A-PT} of this proposition, i.e., to show that
\begin{align}\lbl{C-A}
C_1\cbg C_2\cbg(C+r\BX)\ne\emp.
\end{align}
\par Let $S=\{u',u'',v',v''\}$ and let $\tS=\{x,u,u',u'',v,v',v''\}$. Claim \reff{CL-D} tells us that there exists a $\ds$-Lipschitz mapping $\tf_{S}:\tS\to\X$ with $\ds$-Lipschitz seminorm $\|\tf_{S}\|_{\Lip((\tS,\,\ds),X)}\le 1$
such that
$$
\tf_S(u')\in F(u'),~~\tf_S(u'')\in F(u''),~~
\tf_S(v')\in F(v')~~~\text{and}~~~\tf_S(v'')\in F(v'').
$$
\par Prove that $\tf_S(u)$ belongs to the left hand side of  \rf{C-A}. Indeed, thanks to the inequality $\|\tf_{S}\|_{\Lip((\tS,\,\ds),X)}\le 1$, we have
$$
\|\tf_S(u')-\tf_S(u)\|\le \ds(u',u),~~~
\|\tf_S(u'')-\tf_S(u)\|\le \ds(u'',u)~~~
\|\tf_S(x)-\tf_S(u)\|\le \ds(x,u)=r.
$$
and $\|\tf_S(x)-\tf_S(u)\|\le \ds(x,u)=r$. Thanks to these properties and \rf{CI-S}, $\tf_S(u)\in C_1\cbg C_2$.
\smsk
\par In a similar way we show that $\tf_S(x)\in T_x(v,v',v'')=C$, see \rf{CI-S} and \rf{H-D}. Hence, we have $\tf_S(u)\in C+r\BX$. Thus,
$$
C_1\cbg C_2\cbg(C+r\BX) \ni \tf_S(u)
$$
proving \rf{C-A}.
\smsk
\par We see that property \rf{A-PT} of Proposition \reff{P-F3} holds. This proposition tells us that
\begin{align}
A&=[(C_1\cbg C_2)+Lr\BX]\cbg C+\theta(L)\,\ve\,\BX
\nn\\
&\supset [(C_1\cbg C_2)+(Lr+\ve)\BX]
\cbg \{[(C_1+r\BX)\cbg C] +\ve \BX\}
\cbg \{[(C_2+r\BX)\cbg C] +\ve \BX\}
\nn\\
&=A_1\cbg A_2\cbg A_3.
\nn
\end{align}
\par Prove that
\begin{align}\lbl{S-AI}
A_i\supset F^{[2]}(y)~~~\text{for every}~~~i=1,2,3.
\end{align}
\par We begin with the set
$$
A_1=[C_1\cbg C_2]+(Lr+\ve)\BX.
$$
Thus,
$$
A_1=
[(F(u')+\ds(u',u)\BX)\cbg (F(u'')+\ds(u'',u)\BX)]
+(L\ds(u,x)+L\,\theta(L)\ds(x,y))\BX\,.
$$
See \rf{CI-S}. By the triangle inequality,
$$
\ds(u,x)+\theta(L) \ds(x,y)\ge \ds(u,x)+\ds(x,y)\ge \ds(u,y).
$$
Hence,
$$
A_1\supset
[(F(u')+\ds(u',u)\BX)\cbg (F(u'')+\ds(u'',u)\BX)]
+L\ds(u,y)\BX=T_y(u,u',u'')\,.
$$
But $T_y(u,u',u'')\supset F^{[2]}(y)$, see \rf{G-Y}, so that $A_1\supset F^{[2]}(y)$.
\msk

\par We turn to the proof of the inclusion
$A_2\supset F^{[2]}(y)$. Note that $A_2$ is defined by
\begin{align}\lbl{DF-S2}
A_2=[(C_1+r\BX)\cbg C] +\ve \BX.
\end{align}
\par By the triangle inequality,
\begin{align}\lbl{N-10}
C_1+r\BX=F(u')+\ds(u',u)\BX+\ds(u,x)\BX\supset
F(u')+\ds(u',x)\BX\,.
\end{align}

\par Let
\begin{align}\lbl{CI-DF}
\tC=F(u')+\ds(u',x)\BX,~~~\tC_1=F(v')+\ds(v',v)\BX,~~~
\tC_2=F(v'')+\ds(v'',v)\BX,
\end{align}
and let
\begin{align}\lbl{TR}
\tr=\ds(v,x).
\end{align}
\par In these settings
$$
C=T_x(v,v',v'')=[\tC_1\cbg\tC_2]+L\tr \BX.
$$
\par Let
\begin{align}\lbl{TC-1}
\tA=\{[(\tC_1\cbg\tC_2)+L\tr \BX]\cbg \tC\}+\ve\BX.
\end{align}
Then, thanks to \rf{DF-S2} and \rf{N-10},
$$
A_2\supset
[(F(u')+\ds(u',x)\BX)\cbg C]+\ve \BX=
\{[(\tC_1\cbg\tC_2)+L\tr \BX]\cbg \tC\}+\ve \BX=\tA\,.
$$
See Fig. 23.

\begin{figure}[h!]
\centering{\includegraphics[scale=0.6]{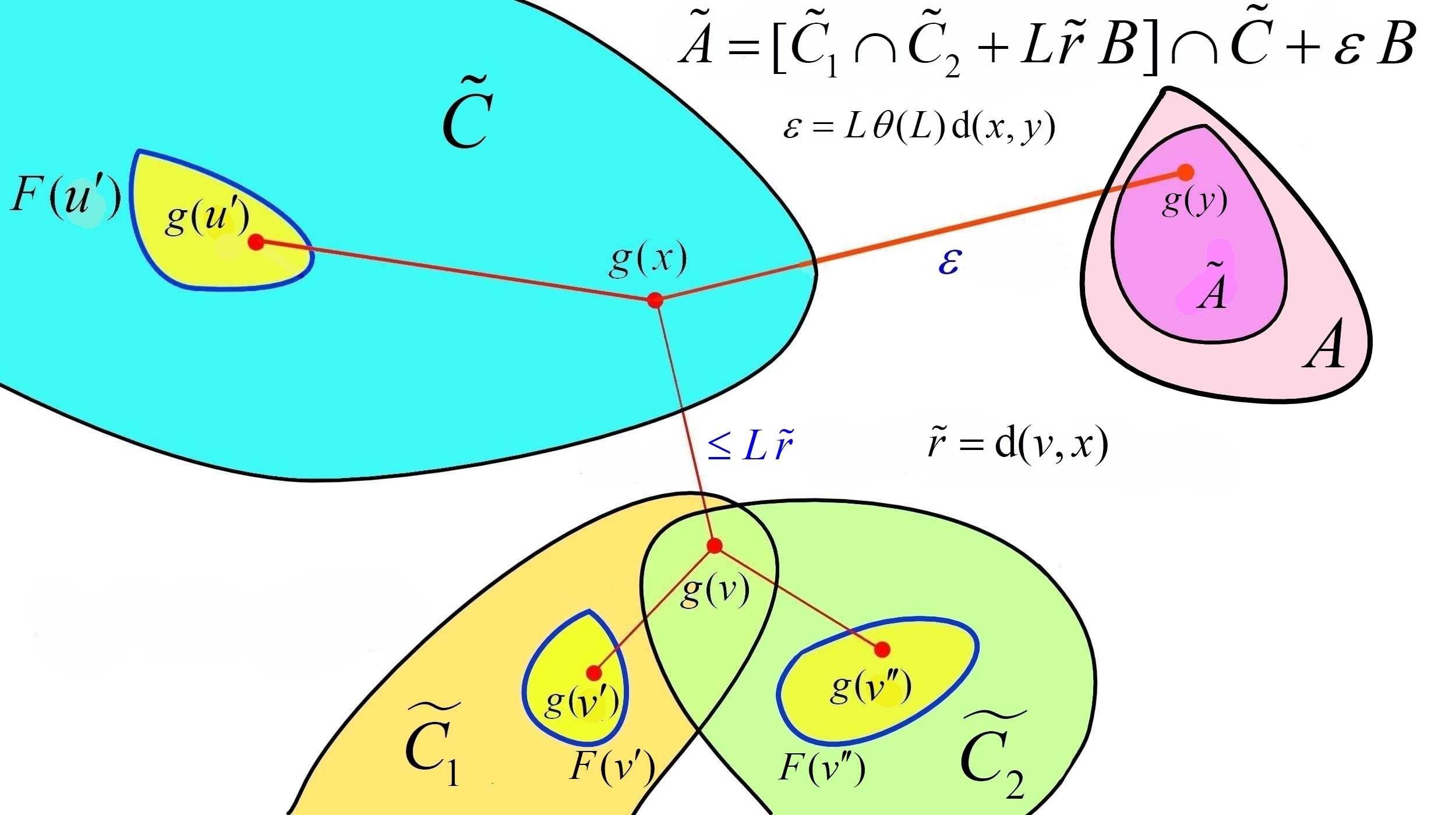}}
\caption{The sets $\tC$, $\tC_1$ and $\tC_2$.}
\end{figure}

\par Prove that $\tA\supset F^{[2]}(y)$. We will do this by applying Proposition \reff{P-F3} to the set $\tA$. But first we must check the hypothesis of this proposition, i.e., to show that
\begin{align}\lbl{PR-1E}
\tC_1\cbg\tC_2\cbg(\tC+\tr\BX)\ne\emp\,.
\end{align}

\par To establish this property, we set $S=\{u',v',v''\}$
and $\tS=\{x,u',v,v',v''\}$.
Claim \reff{CL-D} tells us that there exists a $\ds$-Lipschitz mapping $\tf_{S}:\tS\to\X$ with $\ds$-Lipschitz seminorm $\|\tf_{S}\|_{\Lip((\tS,\,\ds),X)}\le 1$ such that
$$
\tf_{S}(u')\in F(u'),~~\tf_{S}(v')\in F(v')~~~\text{and}~~~ \tf_{S}(v'')\in F(v'').
$$
\par Combining these properties of $\tf_{S}$ with  definitions \rf{CI-DF} and \rf{TR}, we conclude that
$\tC_1\cbg\tC_2\cbg(\tC+\tr\BX)\ni \tf_{S}(v)$
proving \rf{PR-1E}.
\smsk
\par We recall that $\ve=L\,\theta(L)\ds(x,y)$, see \rf{DL}, so that
$$
\tA=
[\{(\tC_1\cbg \tC_2)+L\tr \BX\}\cbg \tC]
+L\,\theta(L)\ds(x,y)\BX, ~~~~\text{(see \rf{TC-1}).}
$$
We apply Proposition \reff{P-F3} to $\tA$ and obtain the following:
\begin{align}
\tA&\supset \{(\tC_1\cbg \tC_2)+(L\tr+L\ds(x,y))\BX\}\nn\\
&\cbg \{[(\tC_1+\tr \BX)\cbg \tC] +L\ds(x,y)\BX\}
\cbg \{[(\tC_2+\tr \BX)\cbg \tC] +L\ds(x,y)\BX\}\nn\\
&=\tA_1\cbg\tA_2\cbg\tA_3.
\nn
\end{align}
\par Prove that $\tA_i\supset F^{[2]}(y)$ for every $i=1,2,3$. First, let us show that
\begin{align}\lbl{S1-ING}
\tA_1=(\tC_1\cbg \tC_2)+(L\tr+L\,\ds(x,y))\BX\supset F^{[2]}(y).
\end{align}
By \rf{TR} and the triangle inequality,
$\tr+\ds(x,y)=\ds(v,x)+\ds(x,y)\ge \ds(v,y)$ so that
\begin{align}
\tA_1&\supset
(\tC_1\cbg \tC_2)+L\ds(v,y)\BX\nn\\
&=
[(F(v')+\ds(v',v)\BX)\cbg
(F(v'')+\ds(v'',v)\BX)]+L\ds(v,y)\BX\nn\\
&=T_y(v,v',v'').\nn
\end{align}
See \rf{CI-DF} and \rf{H-D}. This inclusion and \rf{G-Y} imply \rf{S1-ING}.
\smsk
\par Next, let us show that
\begin{align}\lbl{S2-IN}
\tA_2=[(\tC_1+\tr \BX)\cbg \tC] +L\ds(x,y)\BX\supset F^{[2]}(y)).
\end{align}
\par Thanks to \rf{CI-DF}, \rf{TR} and the triangle inequality,
$$
\tC_1+\tr \BX=F(v')+\ds(v',v)\BX+\ds(v,x)\BX\supset
F(v')+\ds(v',x)\BX
$$
so that
$$
\tA_2\supset
[(F(v')+\ds(v',x)\BX)\cbg (F(u')+\ds(u',x)\BX)]
+L\,\ds(x,y)\BX=T_y(x,u',v').
$$
See \rf{H-D}. From this and \rf{G-Y}, we have $\tA_2\supset T_y(x,u',v')\supset F^{[2]}(y)$, proving   \rf{S2-IN}.
\msk
\par In the same way we show that
$$
\tA_3=[(\tC_2+\tr \BX)\cbg \tC] +L\ds(x,y) \BX\supset
T_y(x,u',v'')\supset F^{[2]}(y)\,.
$$
\par Combining this with \rf{S1-ING} and \rf{S2-IN}, we obtain the required inclusion $\tA_i\supset F^{[2]}(y)$ for every $i=1,2,3$. In turn, this proves that
$$
\tA\supset F^{[2]}(y)~~~~~\text{because}~~~~
\tA\supset \tA_1\cbg\tA_2\cbg\tA_3\supset F^{[2]}(y).
$$
\par We know that $A_2\supset\tA$, so that $A_2\supset F^{[2]}(y)$. In the same fashion we show that
$$
A_3=[(C_2+r\BX)\cbg C] +L\,\ve\,\BX\supset F^{[2]}(y)
$$
proving \rf{S-AI}. Hence,
$A\supset A_1\cbg A_2\cbg A_3\supset F^{[2]}(y)$ so that
\rf{A-1} holds.
\smsk
\par Thus, \rf{FIL} is proved. By interchanging the roles of elements $x$ and $y$ in this inclusion we obtain the inclusion
$F^{[2]}(y)+\gamma_0(L)\ds(x,y)\BX\supset F^{[2]}(x)$.
These two inclusions imply inequality \rf{HD-KJ} proving the proposition.\bx
\msk
\smsk

\par We are in a position to complete the proof of Theorem \reff{MAIN-RT}.
\par Recall that $\lambda_1,\lambda_2$ are parameters satisfying \rf{GM-FN}, and $L=\lambda_2/\lambda_1$. Thus, $L$ and $\lambda_1$ satisfy \rf{L-GE3}. We also recall that $\ds=\lambda_1\,\rho$, see \rf{D-AL}. Let $\gamma$ be a parameter satisfying \rf{GM-FN}.
\par In these settings, the mappings $F^{[1]}$ and $F^{[2]}$ defined by formulae \rf{F-1G} and \rf{F-2G} are the first and the second order $(\{\lambda_1,\lambda_2\},\rho)$-balanced refinements of $F$ respectively. See Definition \reff{F-IT}.
\par Proposition \reff{N-EM} tells us that, under these conditions, $F^{[2]}(x)\ne\emp$ on $\Mc$. In turn, Proposition \reff{HD-G1} states that
$$
\dhf(F^{[2]}(x),F^{[2]}(y))\le \gamma_0(L)\ds(x,y)~~~~~\text{for all}~~~~~x,y\in\Mc.
$$
\par Recall that $\gamma_0(L)=L\cdot\theta(L)^2$
where $\theta=\theta(L)=(3L+1)/(L-1)$, see \rf{TH-L}. Hence,
$\theta(L)=(3\lambda_2+\lambda_1)/(\lambda_2-\lambda_1)$ so that
$$
\dhf(F^{[2]}(x),F^{[2]}(y))\le\gamma_0(L)\ds(x,y)\le L\cdot\theta(L)^2\,\ds(x,y))
=\lambda_2\,\frac{(3\lambda_2+\lambda_1)^2}
{(\lambda_2-\lambda_1)^2}
\,\rho(x,y).
$$
\par Combining this inequality with the third inequality in \rf{GM-FN}, we obtain \rf{HD-RT} proving Theorem \reff{MAIN-RT} for the parameters $\lambda_1,\lambda_2$ and $\gamma$ satisfying \rf{GM-FN}.
\par In particular, one can set $\lambda_1=4/3$, $\lambda_2=4$, and $\gamma=100$. Indeed, in this case $e(\mfM,\X)\le \lambda_1=4/3$, see Remark \reff{EXT-CN}, so that $\lambda_1,\lambda_2$ and $\gamma$ satisfy \rf{GM-FN}.
\smsk
\par Next, let $X$ be a two dimensional Euclidean space, and let $\lambda_1,\lambda_2$ and $\gamma$ satisfy \rf{X-HSP}. In this case, we prove \rf{F2-NEM} and \rf{HD-RT} by replacing in the proof of Theorem \reff{MAIN-RT} the function $\theta=\theta(L)$ defined by \rf{TH-L} with the function $\theta(L)=1+2\,L/\sqrt{L^2-1}$. We leave the details to the interested reader.
\par In particular, we can set $\lambda_1=4/\pi$, $\lambda_2=12/\pi$, $\gamma=38$. Indeed, in this case,  $e(\mfM,X)\le 4/\pi$, see Remark \reff{EXT-CN}, which implies \rf{X-HSP} for these values of parameters $\lambda_1$, $\lambda_2$ and $\gamma$.
\par Finally, suppose that $\X$ is a Euclidean space, $\Mc$ is a subset of a Euclidean space $E$, and $\rho$ is the Euclidean metric in $E$. In this case  $e(\mfM,\X)=1$, see Remark \reff{EXT-CN}, so that one can set $\lambda_1=1$, $\lambda_2=3$ and $\gamma=25$. Clearly, in this case inequalities \rf{X-HSP} hold.
\smsk
\par The proof of Theorem \reff{MAIN-RT} is complete.\bx

\SECT{4. Balanced refinements of line segments in a Banach space.}{4}

\indent
\par In this section we prove Theorem \reff{X-LSGM}.
\par Let $(\Mc,\rho)$ be a pseudometric space, and let  $(\X,\|\cdot\|)$ be a Banach space. We assume that $\dim X>1$. Let us recall that $\Kc_1(\X)$ is the family of all non-empty convex compacts in $\X$ of dimension at most $1$
(i.e., the family of all points and all bounded closed line segments in $\X$).
\par We need the following version of one dimensional Helly's Theorem.
\begin{theorem}\label{HT-IX} Let $\Kc$ be a collection of closed convex subsets of $\X$ containing a set $K_0\in\Kc_1(\X)$. If $K_0$ has a common point with any two members of $\Kc$, then there exists a point common to all of the collection $\Kc$.
\end{theorem}
\par {\it Proof.} We introduce a family $\tKc=\{K\cbg K_0:K\in\Kc\}$, and apply to $\tKc$ one dimensional Helly's Theorem. (See part {\it (a)} of Lemma \reff{H-R}.) \bx
\smsk
\par We also need the following variant of Proposition \reff{P-F3} for the family $\Kc_1(\X)$.
\begin{proposition}\label{C123} Let $\X$ be a Banach space, and let $r\ge 0$. Let $C,C_1,C_2\subset \X$ be convex closed subsets, and let $C_1\in\Kc_1(\X)$. Suppose that
\begin{align}\lbl{A-PT-1}
C_1\cbg C_2\cbg(C+r\BX)\ne\emp.
\end{align}
\par Then for every $L>1$ and every $\ve>0$ the following inclusion
$$
[\{(C_1\cbg C_2)+Lr\BX\}\cbg C]+\theta(L)\,\ve\BX
\supset[(C_1\cbg C_2)+(Lr+\ve)\BX]
\cbg[\{(C_1+r\BX)\cbg C\}+\ve\BX]
$$
holds. Here $\theta(L)$ is the same as in Theorem \reff{N-S}, i.e., $\theta(L)=(3L+1)/(L-1)$ for an arbitrary $X$, and $\theta(L)=1+2L/\sqrt{L^2-1}$ whenever $\X$ is a  Euclidean space.
\end{proposition}

\par {\it Proof.} Let
\bel{A-IN1}
a\in
[(C_1\cap C_2)+(Lr+\ve)\BX]
\cap
\{[(C_1+r\BX)\cap C]+\ve\BX\}.
\ee
\par Prove that
\bel{A-7}
a\in \{[(C_1\cap C_2)+Lr\BX]\cap C\}+\theta(L)\,\ve\BX\,.
\ee
\par First, let us show that
\bel{N-EM1-1}
C_1\cap C_2\cap(C+r\BX)\cap \BX(a,Lr+\ve)\ne\emp.
\ee
\par Recall that $C_1\in\Kc_1(\X)$. Helly's Theorem \reff{HT-IX} tells us that it is suffices to show that any two sets in the left hand size of \rf{N-EM1-1} have a common point with $C_1$.
\par First we note that $C_1$, $C_2$ and $C+r\BX$ have a common point. See \rf{A-PT-1}. We also know that
$$
a\in (C_1\cap C_2)+(Lr+\ve)\BX,
$$
see \rf{A-IN1}, so that $C_1\cap C_2\cap \BX(a,Lr+\ve)\ne\emp$.
\par Let us prove that
\bel{N-51}
C_1\cap (C+r\BX)\cap \BX(a,Lr+\ve)\ne\emp\,.
\ee
Property \rf{A-IN1} tells us that
$a\in [(C_1+r\BX)\cap C] +\ve\BX$. Let $b\in X$ be a point nearest to $a$ on $(C_1+r\BX)\cap C$, and let $b_1\in X$ be a point nearest to $b$ on $C_1$. Then
$\|b_1-b\|\le r$ and $\|a-b\|\le \ve$.
Moreover, $b_1\in C_1\cap(C+r\BX)$ and
$$
\|a-b_1\|\le \|a-b\|+\|b-b_1\|\le \ve +r\le \ve +Lr,
$$
so that $b_1\in \BX(a,Lr+\ve)$. Hence,
$$
b_1\in C_1\cap(C+r\BX)\cap \BX(a,Lr+\ve)
$$
proving \rf{N-51}.
\smsk
\par Thus, \rf{N-EM1-1} holds proving the existence of a point $x\in \X$ such that
\bel{N-61}
x\in C_1\cap C_2\cap (C+r\BX)\cap \BX(a,Lr+\ve)\,.
\ee
\par In particular, $x\in C+r\BX$ so that
$\BX(x,r)\cap C\ne\emp$ proving that condition \rf{C-PR} of Proposition \reff{N-S} is satisfied. We apply this proposition to $x$, $r$ and the set $C$ and get:
$$
[C\cap \BX(x,Lr)]+\theta(L)\,\ve\BX\supset
(C+\ve\BX)\cap
(\BX(x,Lr)+\ve \BX)=
(C+\ve\BX)\cap \BX(x,Lr+\ve)\,.\nn
$$
\par From \rf{N-61} we learn that $a\in \BX(x,Lr+\ve)$. In turn, \rf{A-IN1} tells us that
$$
a\in [(C_1+r\BX)\cap C]+\ve\BX\subset
C+\ve\BX.
$$
Hence, $(C+\ve\BX)\cap \BX(x,Lr+\ve)\ni a$ proving that
$[C\cap \BX(x,Lr)]+\theta(L)\,\ve\BX\ni a$.
\smsk
\par Finally, \rf{N-61} tells us that $x\in C_1\cap C_2$ proving the required inclusion \rf{A-7}.\bx

\msk
\par We recall that $\dim\X>1$ so that the finiteness number $N(1,\X)=\min\{2^{2},2^{\dim \X}\}=4$. Let $F:\Mc\to \Kc_1(\X)$ be a set-valued mapping. We suppose that $F$ satisfies the hypothesis of Theorem \reff{X-LSGM}, i.e., the following statement is true.
\smsk
\begin{claim}\label{A-XK1} For every subset $\Mc'\subset\Mc$ with $\#\Mc'\le 4$ the restriction $F|_{\Mc'}$ of $F$ to $\Mc'$ has a Lipschitz selection $f_{\Mc'}:\Mc'\to \X$ with $\|f\|_{\Lip(\Mc',\X)}\le 1$.
\end{claim}
\par Let $\vl=\{\lambda_1,\lambda_2\}$, and let $F^{[1]}$ and $F^{[2]}$ be the first and the second order $(\vl,\rho)$-balanced refinements of $F$.
See Definition \reff{F-IT}. Our aim is to show that if
\begin{align}\lbl{GM-FN-P}
\lambda_1\ge 1, ~~~~~~\lambda_2\ge 3\lambda_1,
~~~~~~
\gamma\ge \lambda_2\,(3\lambda_2+\lambda_1)/
(\lambda_2-\lambda_1),
\end{align}
then the set-valued mapping $F^{[2]}$ is a $\gamma$-core of $F$ (with respect to $\rho$), i.e.,
$$
F^{[2]}(x)\ne\emp~~\text{for every}~~x\in\Mc,~~~\text{and}~~~
\dhf(F^{[2]}(x),F^{[2]}(y))\le \gamma\rho(x,y)~~\text{for all}~~x,y\in\Mc.
$$
\par To prove this, we set $L=\lambda_2/\lambda_1$ and  introduce a new pseudometric on $\Mc$ defined by
$$
\ds(x,y)=\lambda_1\rho(x,y),~~~~~x,y\in \Mc.
$$
Note that, thanks to \rf{GM-FN-P},
\begin{align}\lbl{L-DS2}
L\ge 3~~~~~~\text{and}~~~~~~ \rho\le \ds ~~~~\text{on}~~~~\Mc.
\end{align}
\par We also note that in these settings, for every $x\in\Mc$,
\begin{align}\lbl{F12-X1}
F^{[1]}(x)=
\bigcap_{z\in\Mc}\,
\left[F(z)+\ds(x,z)\BX\right]~~\text{and}~~
F^{[2]}(x)=\bigcap_{z\in\Mc}\,
\left[F^{[1]}(z)+L\ds(x,z)\BX\right].
\end{align}
\par Next, we need the following analog of Lemma \reff{H-IN}.
\begin{lemma}\label{H-1} Let $\Kc$ be a collection of convex closed subsets of $\X$ containing a set $K_0\in\Kc_1(\X)$. Suppose that $\cbg\{K:K\in\,\Kc\}\ne\emp$. Then for every $r\ge 0$, we have
$$
\left(\,\bigcap_{K\in\,\Kc} K\right) +r\BX
=\bigcap_{K\in\,\Kc}\,
\left\{\,\left[\,K \,\cbig\, K_0\right]+r\BX\,\right\}.
$$
\end{lemma}
\par {\it Proof.} Let $\tKc=\{K\cbg K_0:K\in\Kc\}$. Clearly, $\tKc\subset \Kc_1(\X)$. It is also clear that the statement of the lemma is equivalent to the equality
$$
\left(\,\bigcap_{\tK\in\,\tKc} \tK\right) +r\BX
=\bigcap_{\tK\in\,\tKc}\,
\left\{\,\tK +r\BX\,\right\}
$$
provided $\cbg\{\tK:\tK\in\Kc\}\ne\emp$. We prove this equality by a slight modification of the proof of Lemma \reff{H-IN}. More specifically, we obtain the result by using in that proof Helly's Theorem \reff{HT-IX} instead of Theorem \reff{H-TH}. We leave the details to the interested reader.\bx
\smsk
\begin{lemma}\label{FZ-E} For every $x\in\Mc$ the set $F^{[1]}(x)$ belongs to the family $\Kc_1(X)$. Moreover, for every $x,z\in\Mc$, we have
\begin{align}\lbl{FD-2}
F^{[1]}(z)+L\ds(x,z)\BX=\bigcap_{v\in\Mc}
\left\{\left[(F(v)+\ds(z,v)\BX)\cbg F(z)\right] +L\ds(x,z)\BX\right\}.
\end{align}
\end{lemma}
\msk
\par {\it Proof.} Let $\Kc=\{F(z)+\ds(z,x)\BX:z\in\Mc\}.$
Clearly, $\Kc$ is a family of all bounded closed convex subsets of $\X$ containing the set $F(x)\in\Kc_1(\X)$.
\par Theorem \reff{HT-IX} tells us that the set
$
F^{[1]}(x)=\cbg\{K:K\in\Kc\}
$
is non-empty whenever for every $z',z''\in\Mc$ the set
\begin{align}\lbl{W-1}
E(x,z',z'')=F(x)\cbg [F(z')+\ds(z',x)\BX]\cbg [F(z'')+\ds(z'',x)\BX]
\ne\emp.
\end{align}
\par Fix $z',z''\in\Mc$ and set $\Mc'=\{x,z',z''\}$. Thanks to Claim \reff{A-XK1}, there exists a function $f_{\Mc'}:\Mc'\to\X$ satisfying the following conditions:
$f_{\Mc'}(x)\in F(x)$, $f_{\Mc'}(z')\in F(z')$, $f_{\Mc'}(z'')\in F(z'')$
and
$$
\|f_{\Mc'}(z')-f_{\Mc'}(x)\|\le \rho(z',x)\le \ds(z',x),
~~~\|f_{\Mc'}(z'')-f_{\Mc'}(x)\|\le \rho(z'',x)\le \ds(z'',x).
$$
See \rf{L-DS2}. Then $f_{\Mc'}(x)\in E(x,z',z'')$ so that \rf{W-1} holds. Hence, $F^{[1]}(x)\ne\emp$ proving that $F^{[1]}(x)\in\Kc_1(X)$.
\smsk
\par It remains to note that equality \rf{FD-2} is immediate from \rf{F12-X1} and Lemma \reff{H-1}.
\smsk
\par The proof of the lemma is complete.\bx
\begin{lemma}\label{GW-2} For every $x\in\Mc$, the following equality
$$
F^{[2]}(x)=\bigcap_{u,u'\in\Mc}
\left\{\left[(F(u')+\ds(u',u)\BX)\cbg F(u)\right] +L\ds(u,x)\BX\right\}
$$
holds.
\end{lemma}
\par {\it Proof.} The lemma is immediate from representation \rf{F12-X1} and Lemma \reff{FZ-E}.\bx
\par Given $x,u,u'\in\Mc$ we put
\begin{align}\lbl{TH-D}
\tT_x(u,u')=
[(F(u')+\ds(u',u)\BX)\cbg F(u)]+L\ds(u,x)\,\BX.
\end{align}
\par Note that, in our terminology, the set $\tT_x(u,u')$ is the orbit of $x$ with respect to the diagram shown in Fig. 24. Thus, $\tT_x(u,u')$ is the family of all points $a=g(x)$ where $g$ runs over all mappings which agree with the diagram in Fig. 24 below.
\newpage
\begin{figure}[h]
\centering{\includegraphics[scale=0.57]{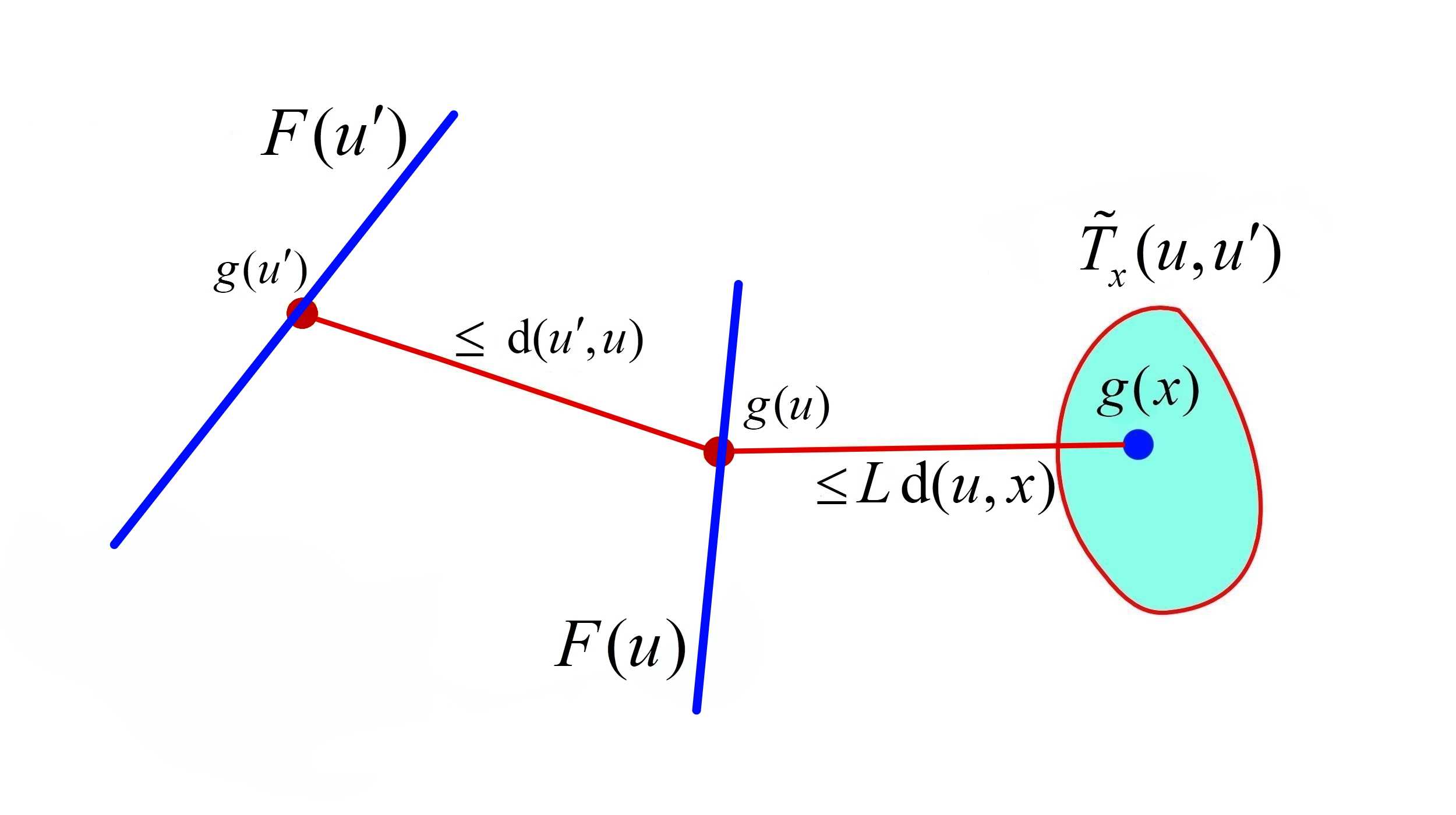}}
\caption{$\tT_x(u,u')$ is the orbit of $x$ with respect to this diagram.}
\end{figure}
\msk
\par Now, Lemma \reff{GW-2} reformulates as follows:
\begin{align}\lbl{G-XP1}
F^{[2]}(x)=\bigcap_{u,u'\in\Mc} \tT_x(u,u').
\end{align}
\begin{proposition}\label{TN-EM} For every $x\in\Mc$ the set $F^{[2]}(x)$ is non-empty.
\end{proposition}
\par {\it Proof.} Recall that $F(x)\in\Kc_1(\X)$. Furthermore, by \rf{TH-D}, $F(x)=\tT_x(x,x)$. Therefore, by \rf{G-XP1} and Helly's Theorem \reff{HT-IX}, the set $F^{[2]}(x)\ne\emp$ whenever for every $u_i,u'_i\in\Mc$, $i=1,2$, we have
\begin{align}\lbl{TH-INT}
F(x)\cbg\tT_x(u_1,u'_1)\cbg \tT_x(u_2,u'_2)\ne\emp.
\end{align}
\par We note that \rf{TH-INT} holds provided there exists a mapping $g$ which agrees with the diagram shown in Fig. 25.
\begin{figure}[h!]
\centering{\includegraphics[scale=0.3]{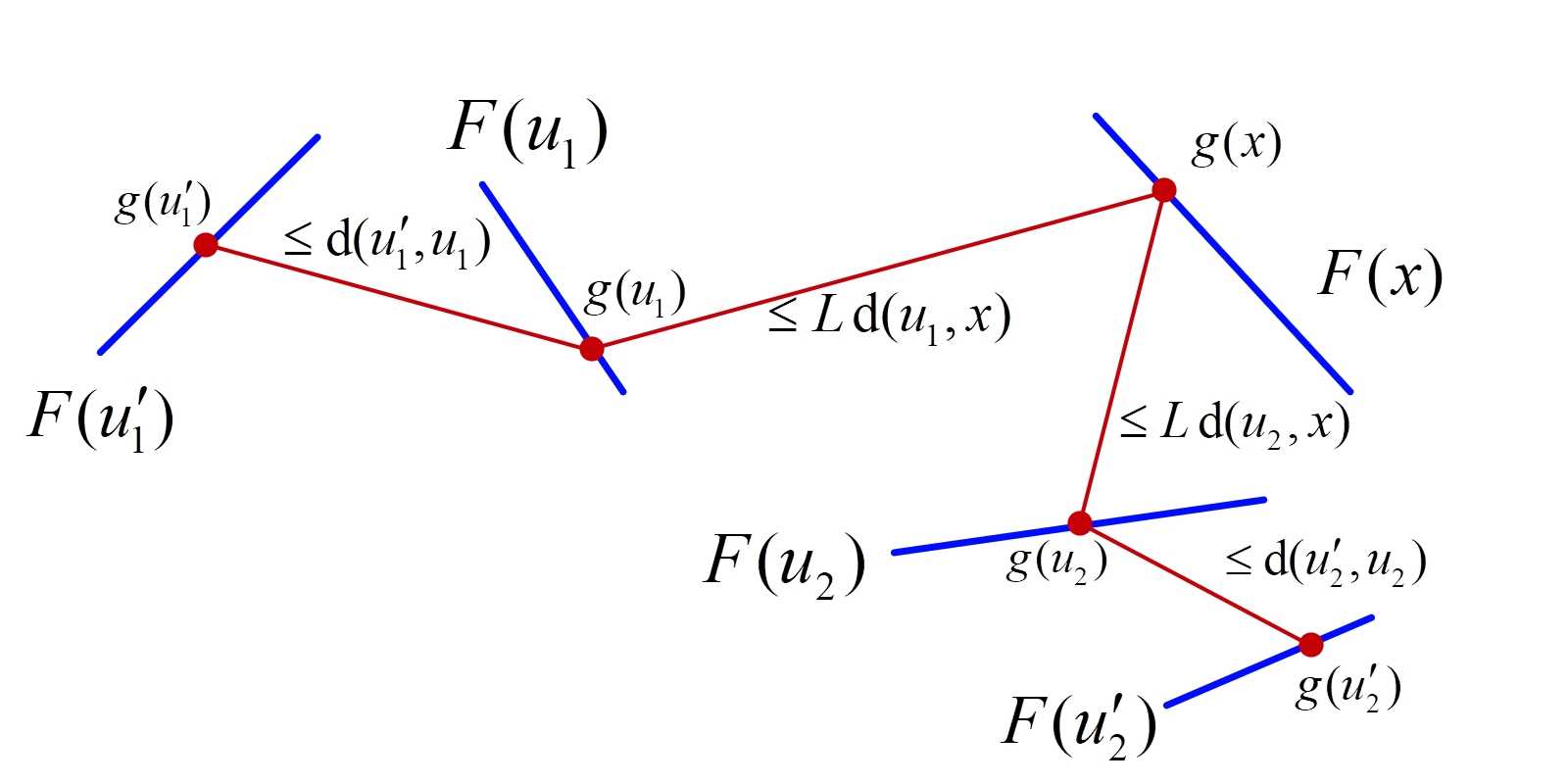}}
\caption{The mapping $g$ agrees with this diagram.}
\end{figure}

\par Thanks to \rf{TH-D},
\begin{align}\lbl{H-DLP}
\tT_x(u_i,u'_i)=[(F(u_i')+\ds(u_i',u_i)\BX)\cbg F(u_i)] +L\ds(u_i,x)\,\BX,~~~i=1,2.
\end{align}
\par Fix elements $u_1,u'_1,u_2,u'_2\in\Mc$ and prove that property \rf{TH-INT} holds.
\par First we note that, without loss of generality, one may assume that $\rho(u_1,x)\ge \rho(u_2,x)$.
Next, we introduce the following sets:
\begin{align}\lbl{G-S123}
G_1=F(u_2),~~~~
G_2=F(u_2')+\rho(u_2,u_2')\BX,~~~~
G_3=F(x)+\rho(u_2,x)\BX,
\end{align}
and
\begin{align}\lbl{G-S4}
G_4=[(F(u_1')+\rho(u_1',u_1)\BX)\cbg F(u_1)]+\rho(u_1,u_2,)\BX\,.
\end{align}
\par Prove that if
\begin{align}\lbl{G-ALL}
G_1\cbg G_2\cbg G_3\cbg G_4\ne\emp
\end{align}
then \rf{TH-INT} holds.
\par (Note that \rf{G-ALL} holds whenever there exists a mapping $g$ which agrees with the diagram in Fig. 26.)

\begin{figure}[h!]
\centering{\includegraphics[scale=0.3]{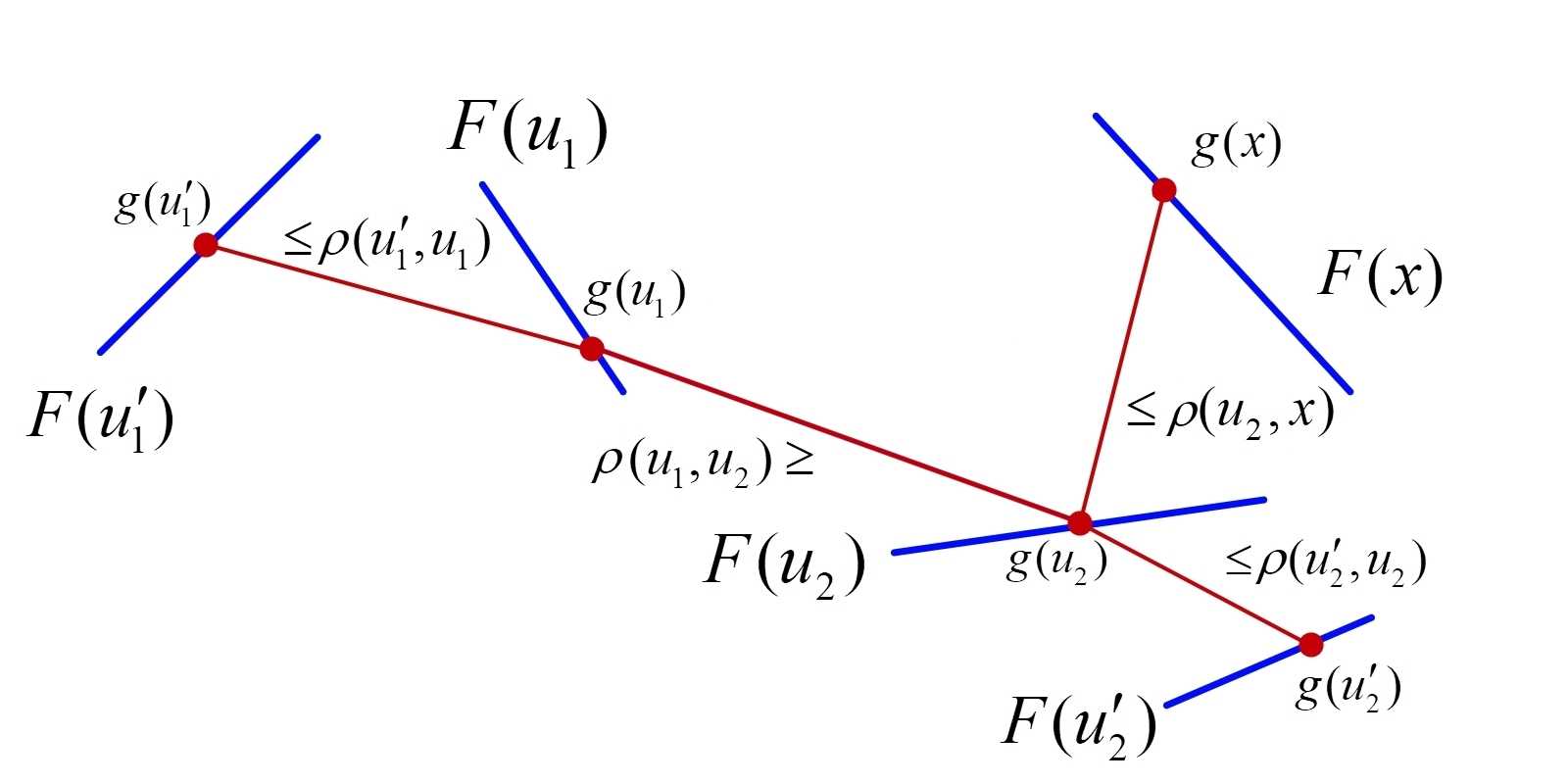}}
\caption{The existence of a mapping $g$ which agrees with this diagram implies \rf{G-ALL}.}
\end{figure}
\par Indeed, let $\tMc=\{u_1',u_1,x,u_2,u_2'\}$. This property and definitions \rf{G-S123}, \rf{G-S4} imply the existence of a mapping $g:\tMc\to X$ with the following properties: $g(v)\in F(v)$ on $\tMc$,
\begin{align}\lbl{G-LP}
\|g(u_1)-g(u_1')\|\le \rho(u_1,u_1'),~~~~
\|g(u_1)-g(u_2)\|\le \rho(u_1,u_2),
\end{align}
and
\begin{align}\lbl{G-U2X}
\|g(u_2)-g(u_2')\|\le \rho(u_2,u_2'),~~~~\|g(u_2)-g(x)\|\le \rho(u_2,x).
\end{align}
\par We establish \rf{TH-INT} by showing that
\begin{align}\lbl{GX-TH}
g(x)\in F(x)\cbg\tT_x(u_1,u'_1)\cbg \tT_x(u_2,u'_2).
\end{align}
\par In fact, from the above properties of $g$ it follows that $g(x)\in F(x)$. We also know that $g(u_2)\in F(u_2)$, $g(u_2')\in F(u_2')$. Thanks to \rf{G-LP}, \rf{G-U2X} and \rf{L-DS2},
$$
\|g(u_2)-g(u_2')\|\le \rho(u_2,u_2')\le \ds(u_2,u_2')
~~~~\text{and}~~~~
\|g(u_2)-g(x)\|\le \rho(u_2,x)\le L\ds(u_2,x).
$$
\par From these properties of $g$ and definition \rf{H-DLP}, we have
$$
g(x)\in
[(F(u_2')+\ds(u_2',u_2)\BX)\cbg F(u_2)]+L\ds(u_2,x)\,\BX=\tT_x(u_2,u'_2).
$$
\par It remains to show that $g(x)\in \tT_x(u_1,u'_1)$. As we know,
\begin{align}\lbl{GU-IN}
g(x)\in F(x),~~~~~g(u_1)\in F(u_1),~~~~~\text{and}~~~~~ g(u_1')\in F(u_1').
\end{align}
Furthermore, thanks to \rf{G-LP} and \rf{L-DS2},
\begin{align}\lbl{GU1}
\|g(u_1)-g(u_1')\|\le \rho(u_1,u_1')\le \ds(u_1,u_1').
\end{align}
\par Let us estimate $\|g(u_1)-g(x)\|$. From \rf{G-LP}, \rf{G-U2X} and the triangle inequality, we have
\begin{align}
\|g(u_1)-g(x)\|&\le \|g(u_1)-g(u_2)\|+\|g(u_2)-g(x)\|\le \rho(u_1,u_2)+\rho (u_2,x)\nn\\
&\le \rho(u_1,x)+\rho(x,u_2)+\rho (u_2,x)=
\rho(u_1,x)+2\rho(x,u_2).\nn
\end{align}
Recall that $\rho(u_1,x)\ge \rho(u_2,x)$. This and \rf{L-DS2} yield
$$\|g(u_1)-g(x)\|\le 3\rho(u_1,x)\le L\ds(u_1,x).$$
\par This inequality, \rf{GU-IN}, \rf{GU1} and \rf{H-DLP} imply the required property $g(x)\in \tT_x(u_1,u'_1)$ proving \rf{GX-TH}.
\smsk

\par Thus, to complete the proof of the proposition, we have to prove that the sets $G_i$, $i=1,...,4$, have a common point.

\par See Fig. 27.

\begin{figure}[h!]
\centering{\includegraphics[scale=0.3]{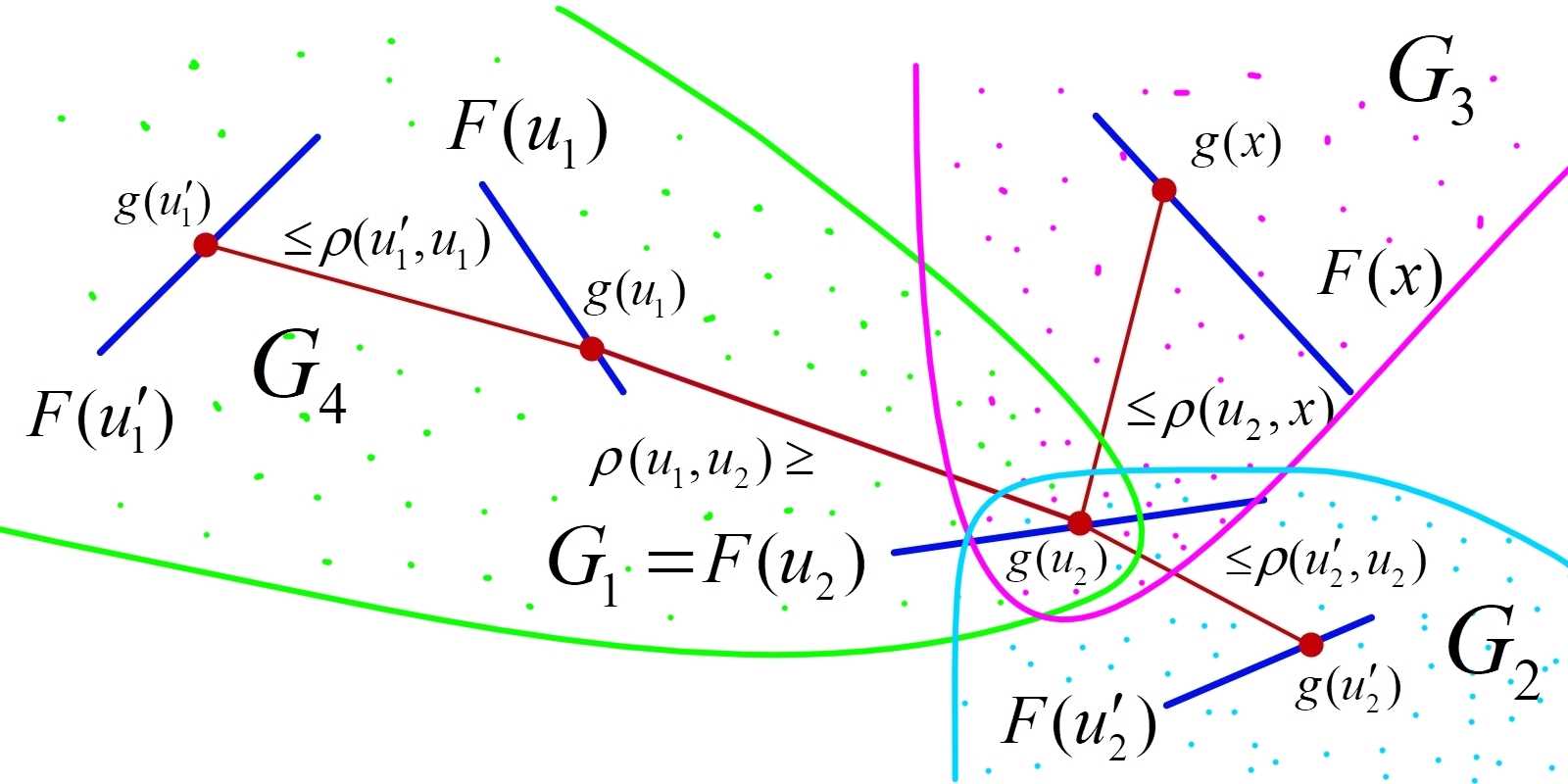}}
\caption{The sets $G_1$, $G_2$, $G_3$ and $G_4$.}
\end{figure}

\par Note that $G_1=F(u_2)\in\Kc_1(X)$ so that, by  Theorem \reff{HT-IX}, this property holds provided
\begin{align}\lbl{G-S-HT}
G_1\cbg G_i\cbg G_j\ne\emp~~~~\text{for every}~~~~2\le i,j\le 4,~i\ne j.
\end{align}
\par First prove that $G_1\cbg G_2\cbg G_3\ne\emp.$
\par Let $\Mc_1=\{u_2',u_2,x\}$. Thanks to Claim \reff{A-XK1}, there exists a mapping $f_1:\Mc_1\to\X$ with the following properties:  $f_1(x)\in F(x)$, $f_1(u_2)\in F(u_2)$, $f_1(u_2')\in F(u_2')$,
$$
\|f_1(u_2)-f_1(x)\|\le \rho(u_2,x)~~~~\text{and}~~~~
\|f_1(u_2)-f_1(u_2')\|\le \rho(u_2,u_2').
$$
These properties of $f_1$ and definition \rf{G-S123} tell us that $f_1(u_2)\in G_1\cbg G_2\cbg G_3$ proving that the sets $G_1$, $G_2$ and $G_3$ have a common point.
\smsk
\par Prove that $G_1\cbg G_2\cbg G_4\ne\emp.$ 
\par Let  $\Mc_2=\{u_1',u_1,u_2',u_2\}$. Using Claim \reff{A-XK1}, we produce a mapping $f_2:\Mc_2\to\X$ such that $f_2(u_i)\in F(u_i)$, $f_2(u_i')\in F(u_i')$ for every $i=1,2$,~ $\|f_2(u_1)-f_2(u_1')\|\le \rho(u_1,u_1')$,
$$
\|f_2(u_1)-f_2(u_2)\|\le \rho(u_1,u_2)~~~~~~
\text{and}~~~~~~
\|f_2(u_2)-f_2(u_2')\|\le \rho(u_2,u_2').
$$
These properties of $f_2$ and \rf{G-S123}, \rf{G-S4} yield $f_2(u_2)\in G_1\cbg G_2\cbg G_4$ proving that the sets $G_1$, $G_2$ and $G_4$ have a common point.
\smsk
\par In the same way we show that $G_1\cbg G_3\cbg G_4\ne\emp$. (We set $\Mc_3=\{u_1',u_1,x,u_2\}$, produce a corresponding function $f_3:\Mc_3\to\X$ and show that $f_3(u_2)\in G_1\cbg G_3\cbg G_4$.)
\par Thus, \rf{G-S-HT} holds, proving that the sets $G_i$ have a common point.
\par The proof of the proposition is complete.\bx
\msk
\par In this section we set
$\gamma_0=\gamma_0(L)=L\,\theta(L)$.
\begin{proposition}\label{HD-1D} For every $x,y\in\Mc$ the following inequality
$$\dhf(F^{[2]}(x),F^{[2]}(y))\le \gamma_0(L)\ds(x,y)$$
holds.
\end{proposition}
\par {\it Proof.} Let $x,y\in\Mc$. Thanks to \rf{G-XP1},
\begin{align}\lbl{GXY-2}
F^{[2]}(x)=\bigcap_{u,u'\in\Mc} \tT_x(u,u')~~~\text{and}~~~
F^{[2]}(y)=\bigcap_{u,u'\in\Mc} \tT_y(u,u').
\end{align}
Recall that
\begin{align}\lbl{H-RM}
\tT_x(u,u')=
[(F(u')+\ds(u',u)\BX)\cbg F(u)]+L\ds(u,x)\,\BX.
\end{align}
\par Recall that the set $\tT_x(u,u')$ is the orbit of $x$ with respect to the diagram of Fig. 24.

\par We know that $F^{[2]}(x)\ne\emp$, see Proposition \reff{TN-EM}, and $\tT_x(x,x)=F(x)\in\Kc_1(\X)$. These properties, \rf{GXY-2} and Lemma \reff{H-1} tell us that
\begin{align}\lbl{THXY}
F^{[2]}(x)+\gamma_0(L)\ds(x,y)\,\BX=
\bigcap_{u,u'\in\Mc}\,\,
\left\{\,\left[\tT_x(u,u')\cbg F(x)\right] +\gamma_0(L)\ds(x,y)\,\BX\,
\right\}.
\end{align}
\par We fix $u,u'\in\Mc$ and introduce a set
$$
\tA=[\tT_x(u,u')\cbg F(x)]+\gamma_0(L)\ds(x,y)\,\BX.
$$
\par Clearly, this set is the orbit of $y$ with respect to the diagram in Fig. 28.
\begin{figure}[h!]
\centering{\includegraphics[scale=0.30]{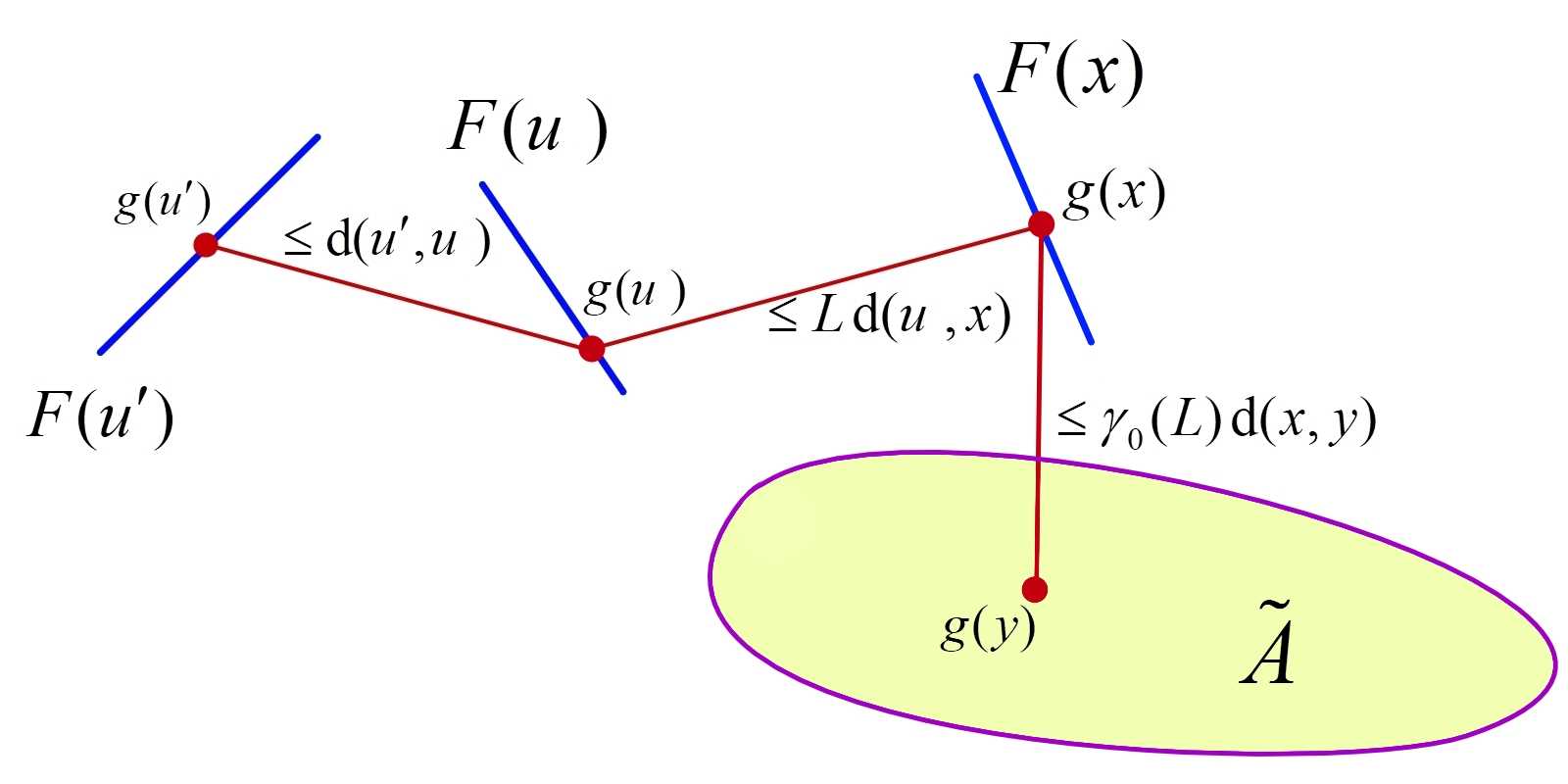}}
\vspace*{-6mm}
\caption{$\tA=\{g(y)\}$ where $g$ runs over all mappings  which agree with this diagram.}
\end{figure}

We also introduce sets
\begin{align}\lbl{C-K1}
C_1=F(u),~~~C_2=F(u')+\ds(u',u)\BX,~~~~\text{and}~~~~C=F(x).
\end{align}
Let
\begin{align}\lbl{DL1}
\ve=L\ds(x,y)~~~\text{and}~~~r=\ds(x,u)\,.
\end{align}
\par In these settings, $\gamma_0(L)\ds(x,y)=\theta(L)\,\ve$
and
$$
\tA=[\tT_x(u,u')\cbg F(x)]+\gamma_0(L)\ds(x,y)\,\BX
=[\{(C_1\cbg C_2)+Lr\BX\}\cbg C]+\theta(L)\,\ve\,\BX.
$$

\par Our aim is to prove that
$\tA\supset F^{[2]}(y)=\cbg\{\tT_y(z,z'):z,z'\in\Mc\}$.
See Fig. 29.
\msk

\begin{figure}[h!]
\centering{\includegraphics[scale=0.30]{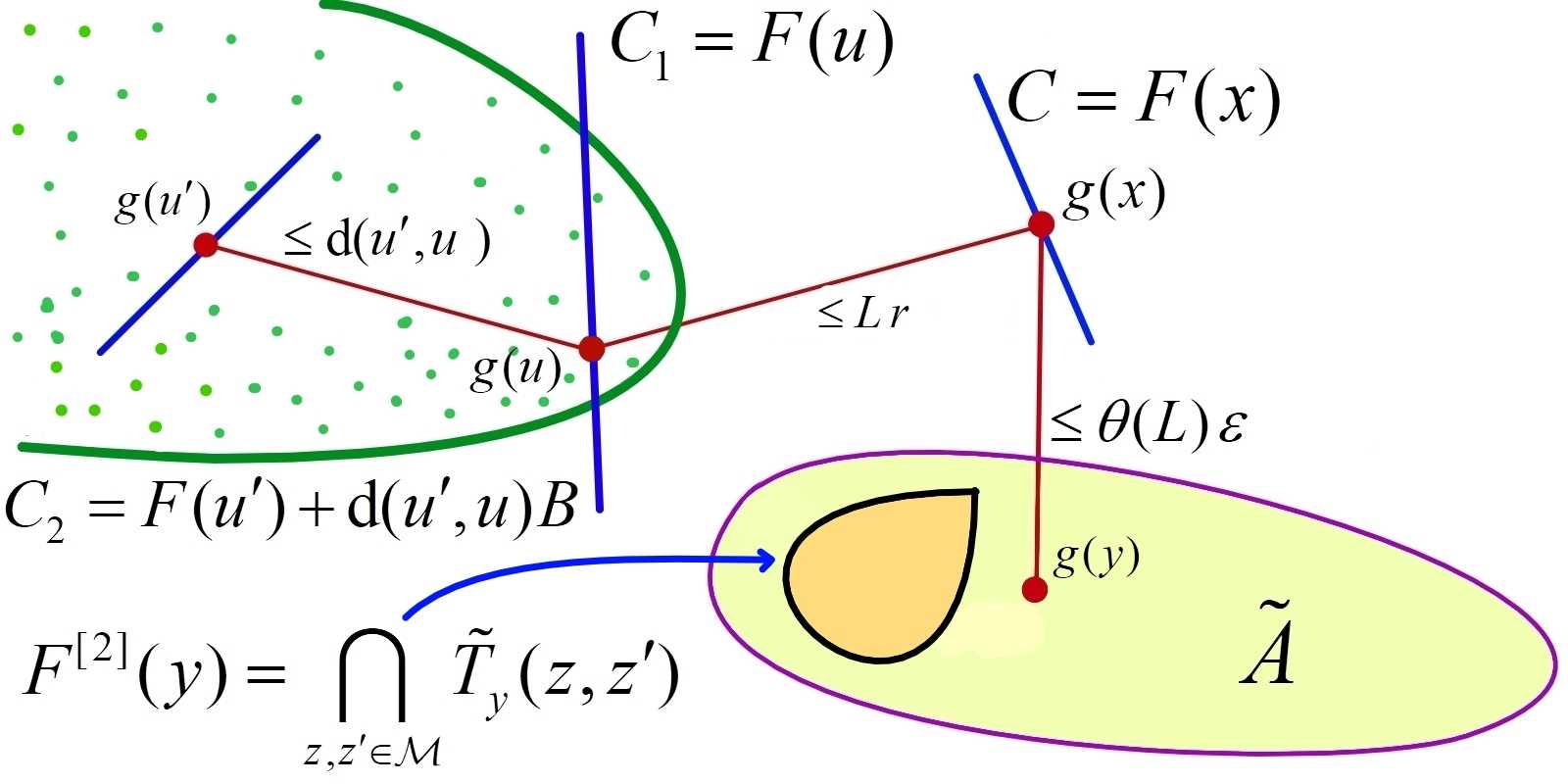}}
\caption{The sets $C_1$, $C_2$ and $C$.}
\end{figure}

\par In fact, we prove that
$$
\tA\supset \tT_y(u,u')\cap \tT_y(x,u').
$$
\par Let us apply Proposition \reff{C123} to the set $\tA$. To do this, we have to verify condition \rf{A-PT-1}, i.e., to show that
\begin{align}\lbl{C-INT}
C_1\cbg C_2\cbg(C+r\BX)\ne\emp.
\end{align}
\par Let $\tMc=\{x,u,u'\}$. Thanks to Claim \reff{A-XK1}, there exists a $\rho$-Lipschitz selection $f_{\tMc}$ of the restriction $F|_{\tMc}$ with $\|f_{\tMc}\|_{\Lip((\tMc;\rho),\X)}\le 1$. Thus, $f_{\tMc}(u')\in F(u')$, $f_{\tMc}(u)\in F(u)$, $f_{\tMc}(x)\in F(x)$,
$$
\|f_{\tMc}(u')-f_{\tMc}(u)\|\le \rho(u',u)~~~
\text{and}~~~
\|f_{\tMc}(x)-f_{\tMc}(u)\|\le \rho(x,u).
$$
\par Let us prove that
\begin{align}\lbl{FM-U}
f_{\tMc}(u)\in C_1\cbg C_2\cbg(C+r\BX).
\end{align}
\par Indeed, $f_{\tMc}(u)\in F(u)=C_1$, see \rf{C-K1}. Furthermore, $f_{\tMc}(u')\in F(u')$ and, thanks to \rf{L-DS2}, $\rho\le \ds$ on $\Mc$. Hence,
$$
\|f_{\tMc}(u')-f_{\tMc}(u)\|\le \rho(u',u)\le \ds(u',u)
$$
proving that $f_{\tMc}(u)\in C_2$, see \rf{C-K1}. Finally, by \rf{C-K1} and \rf{DL1}, $f_{\tMc}(x)\in F(x)=C$ and
$$
\|f_{\tMc}(x)-f_{\tMc}(u)\|\le \rho(x,u)\le \ds(x,u)=r,
~~~~~\text{so that}~~~~f_{\tMc}(u)\in C+r\BX.
$$
\par Thus, \rf{FM-U} is true so that property \rf{C-INT} holds. Furthermore, $C_1=F(u)\in\Kc_1(\X)$, so that all conditions of the hypothesis of Proposition \reff{C123} are satisfied. By this proposition,
\begin{align}
\tA&=[\{(C_1\cbg C_2)+Lr\BX\}\cbg C]+\theta(L)\,\ve\,\BX
\nn\\
&\supset [(C_1\cbg C_2)+(Lr+\ve)\BX]
\cbg [\{(C_1+r\BX)\cbg C\} +\ve\BX]
\nn\\
&=\tA_1\cbg \tA_2.
\nn
\end{align}
\par Prove that $\tA_i\supset F^{[2]}(y)$ for every $i=1,2$.
\smsk
\par We begin with the set 
$$
\tA_1=C_1\cbg C_2+(Lr+\ve)\BX.
$$
\par Thanks to \rf{C-K1} and \rf{DL1},
$$
\tA_1=
[\{F(u')+\ds(u',u)\BX\}\cbg F(u)]
+(L\ds(u,x)+L\ds(x,y))\BX\,.
$$
By the triangle inequality,
$$\ds(u,x)+\ds(x,y)\ge \ds(u,y)$$
so that
$$
\tA_1\supset
[\{F(u')+\ds(u',u)\BX\}\cbg F(u)]
+L\ds(u,y)\BX=\tT_y(u,u'),~~~~~\text{see \rf{H-RM}}.
$$
But, by \rf{GXY-2},  $\tT_y(u,u')\supset F^{[2]}(y)$ which implies the required inclusion $\tA_1\supset F^{[2]}(y)$.

\par See Fig. 30.

\begin{figure}[h!]
\centering{\includegraphics[scale=0.30]{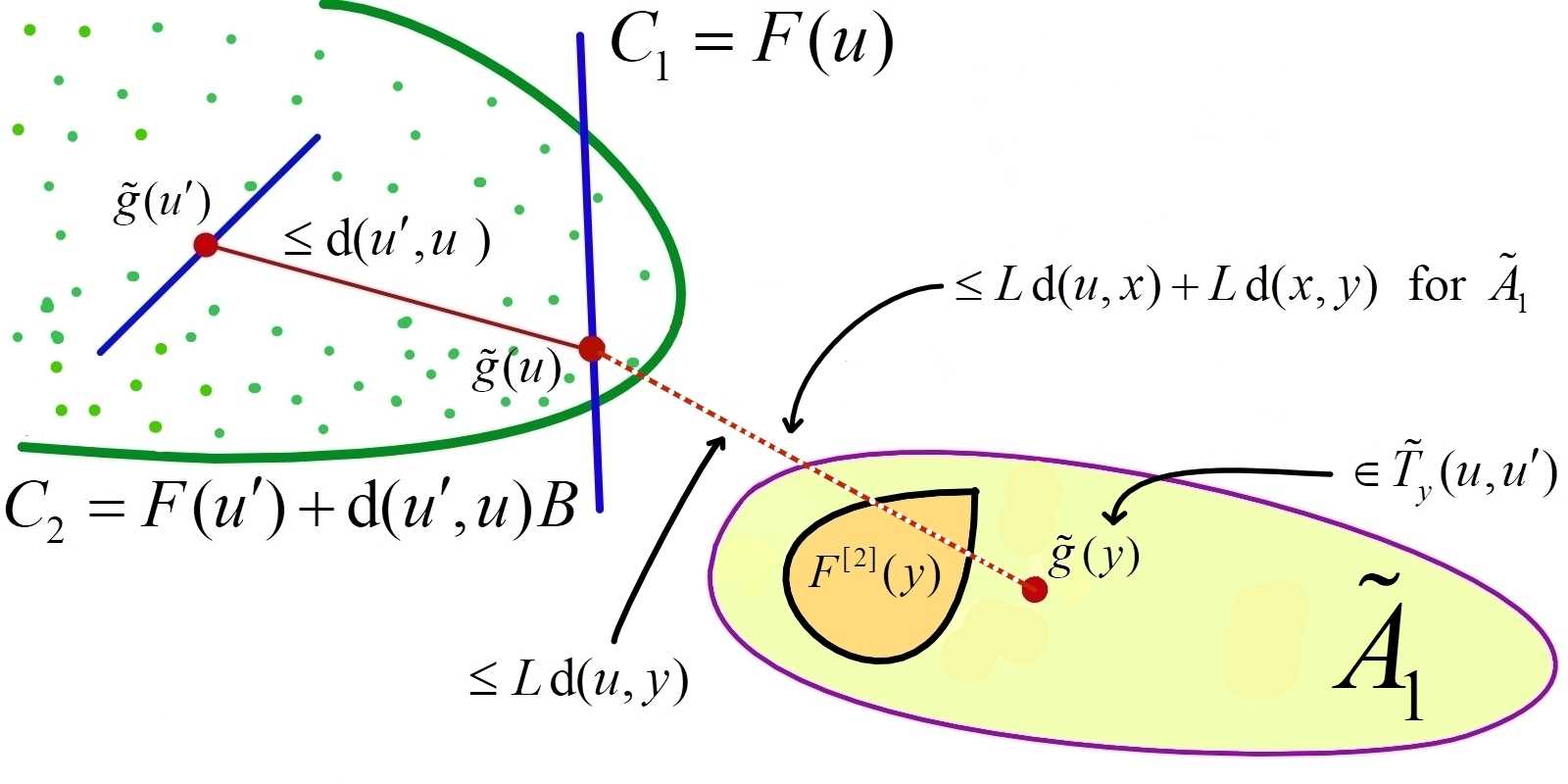}}
\caption{The set $\tA_1=C_1\cbg C_2+(Lr+\ve)\BX$ containing the set $\tT_y(u,u')$.}
\end{figure}

\par We turn to the set $\tA_2=
[(C_1+r\BX)\cbg C] +\ve\BX$. By \rf{TH-D}, \rf{C-K1} and  \rf{DL1},
$$
\tA_2=[(F(u)+\ds(u,x)\BX)\cbg F(x)]+L\ds(x,y)\BX=T_y(x,u).
$$
Thanks to \rf{GXY-2}, $\tT_y(u,x)\supset F^{[2]}(y)$ proving that $\tA_2\supset F^{[2]}(y)$.
\smsk
\par Thus,
$$
\tA=[\tT_x(u,u')\cbg F(x)]+\gamma_0(L)\ds(x,y)\,\BX
\supset \tA_1\cbg \tA_2\supset F^{[2]}(y)~~~~\text{for every}~~~~u,u'\in\Mc.
$$
From this and representation \rf{THXY}, we have
$
F^{[2]}(x)+\gamma_0(L)\ds(x,y)\,\BX\supset F^{[2]}(y)\,.
$
\par By interchanging the roles of $x$ and $y$ we obtain also
$$
F^{[2]}(y)+\gamma_0(L)\ds(x,y)\,\BX\supset F^{[2]}(x)\,.
$$
\par These two inclusions and \rf{HFD} imply the statement of the proposition.\bx
\msk
\par We complete the proof of Theorem \reff{X-LSGM} as follows. We fix $\lambda_1$, $\lambda_2$ and $\gamma$ satisfying inequalities \rf{GM-FN-1}. Then, by Proposition \reff{TN-EM}, the set $F^{[2]}(x)\ne\emp$ for every $x\in\Mc$.
\par In turn, Proposition \reff{HD-1D} tells us that   $\dhf(F^{[2]}(x),F^{[2]}(y))\le \gamma_0(L)\ds(x,y)$ on $\Mc$. We recall that $L=\lambda_2/\lambda_1$, $\ds=\lambda_1\rho$, $\gamma_0(L)=L\,\theta(L)$ and
$\theta(L)=(3L+1)/(L-1)$. Hence,
$$
\dhf(F^{[2]}(x),F^{[2]}(y))\le
\gamma_0(L)\ds(x,y)=L\left(\frac{3L+1}{L-1}\right)\ds(x,y)
=\frac{\lambda_2\,(3\lambda_2+\lambda_1)}
{(\lambda_2-\lambda_1)}\,\rho(x,y).
$$
This inequality and \rf{GM-FN-1} imply \rf{HD-RT}.
\par Thus, \rf{F2-NEM} and \rf{HD-RT} hold provided $\lambda_1$, $\lambda_2$ and $\gamma$ satisfy inequalities \rf{GM-FN-1}. In par\-ti\-cular, we can set
$\lambda_1=1$, $\lambda_2=3$ and
$\gamma=\lambda_2\,(3\lambda_2+\lambda_1)/
(\lambda_2-\lambda_1)=15$.
\par Let now $\X$ be a Euclidean space, and let $\lambda_1$, $\lambda_2$ and $\gamma$ be the parameters satisfying inequalities \rf{X1-H}.

\par In this case, we replace in the above calculations the constant $\theta(L)=(3L+1)/(L-1)$ with
$\theta(L)=1+2L/\sqrt{L^2-1}$. This leads us to the required estimate
$$
\dhf(F^{[2]}(x),F^{[2]}(y))\le
\left\{\lambda_2+2\lambda_2^2\,\slbig
\left(\lambda_2^2-\lambda_1^2\right)^{\frac12}\right\}
\,\rho(x,y)\le \gamma\,\rho(x,y)
$$
proving that \rf{F2-NEM} and \rf{HD-RT} hold for $\lambda_1$, $\lambda_2$ and $\gamma$ satisfying \rf{X1-H}.
\par The proof of Theorem \reff{X-LSGM} is complete.\bx

\SECT{5. The main theorem in $\LTI$.}{5}

\noindent {\bf 5.1 The case $X=\R$.}

\begin{proposition}\label{X-1DIM} Let $\MR$ be a pseudometric space. Let $m=1$ and let $X=\R$; thus, $\ell=\ell(m,X)=1$, see \rf{NMY-1}. In this case Conjecture \reff{BR-IT} holds for every $\lambda_1\ge 1$ and $\gamma\ge 1$.
\par Thus, the following statement is true: Let $F$ be a set-valued mapping from $\Mc$ into the family $\Kc(\R)$ of all closed bounded intervals in $\R$. Suppose that for every $x,y\in\Mc$ there exist points $g(x)\in F(x)$ and $g(y)\in F(y)$ such that $|g(x)-g(y)|\le \rho(x,y)$.
\smsk
\par Let $F^{[1]}(x)$, $x\in\Mc$, be the $\lambda_1$-balanced refinement of the mapping $F$, i.e., the set
\begin{align}\lbl{F1-PR1}
F^{[1]}(x)=
\bigcap_{z\in\Mc}\,\left[F(z)+
\lambda_1\,\rho(x,z)\,\BXR\right]~~~~
\text{where}~~~~I_0=[-1,1].
\end{align}
\par Then $F^{[1]}(x)\ne\emp$ for every $x\in\Mc$, and
$$
\dhf(F^{[1]}(x),F^{[1]}(y))\le \gamma\,\rho(x,y)~~~~\text{for all}~~~~x,y\in\Mc.
$$
\end{proposition}
\smsk
\par Note that Proposition \reff{X-1DIM} easily follows from the one dimensional Helly Theorem and a formula for a neighborhood of the intersection of intervals in $\R$. We formulate these statements in the following
\begin{lemma}\label{H-R} Let $\Kc\subset\Kc(\R)$ be a collection of closed bounded intervals in $\R$.
\smsk
\par (a) (Helly's Theorem in $\R$.) If the intersection of every two intervals from $\Kc$ is non-empty, then there exists a point in $\R$ common to all of the family $\Kc$.
\par (b) Suppose that $\cbig\{K:K\in\Kc\}\ne\emp$. Then for every $r\ge 0$ we have
$$
\left(\,\bigcap_{K\in\,\Kc} K\right) +r\BXR
=\bigcap_{K\in\,\Kc}\,
\left\{\,K+r\BXR\,\right\}.
$$
\end{lemma}
\par {\it Proof of part (b).} In Lemma \reff{H-IN} we proved an analog of property {\it (b)} for $\RT$. We prove part {\it (b)} by replacing in that proof the Helly Theorem \reff{H-TH} with the one dimensional Helly Theorem formulated in part {\it (a)} of the present lemma. We leave the details to the interested reader.\bx
\bsk

\par {\it Proof of Proposition \reff{X-1DIM}.} We have to prove that the set $F^{[1]}(x)$ is non-empty for each $x\in\Mc$, and for every $x,y\in\Mc$
\bel{TM-1}
\dhf(F^{[1]}(x),F^{[1]}(y))\le \gamma\,\rho(x,y).
\ee
We know that the restriction $F|_{\Mc'}$ of $F$ to every two point subset $\Mc'\subset\Mc$ has a Lipschitz selection $f_{\Mc'}:\Mc'\to\R$ with $\|f_{\Mc'}\|_{\Lip(\Mc',\R)}\le \gamma$. Thus, for every $z,z'\in\Mc$ there exist points
\bel{GZ-ZP}
g(z)\in F(z),~g(z')\in F(z')~~~\text{such that}~~~ |g(z)-g(z')|\le \gamma\,\rho(z,z').
\ee
\par We recall that the set-valued mapping $F^{[1]}$ is defined by formula \rf{F1-PR1}.
\par Prove that $F^{[1]}(x)\ne\emp$\, for every $x\in\Mc$.
Indeed, thanks to \rf{F1-PR1} and Helly's Theorem for intervals (part (i) of Lemma \reff{H-R}), $F^{[1]}(x)\ne\emp$ provided
\bel{FFP}
(F(z)+\gamma\,\rho(x,z)\,\BXR)\cap(F(z')+\gamma\,\rho(x,z')
\,\BXR)\ne\emp
\ee
for every $z,z'\in\Mc$.
\par We know that there exist points $g(z)$ and $g(z')$ satisfying \rf{GZ-ZP}. Let
$$
a=\min\{g(z)+\gamma\,\rho(z,x),g(z')+\gamma\,\rho(z',x)\}.
$$
Thanks to the inequality $|g(z)-g(z')|\le \gamma\,\rho(z,z')$, we have
$$
g(z)=\min\{g(z),g(z')+\gamma\,\rho(z',z)\}
$$
so that, by the triangle inequality,
$$
|a-g(z)|\le \max\{\gamma\,\rho(z,x),
|\gamma\,\rho(z',x)-\gamma\,\rho(z',z)|\}
=\gamma\,\rho(z,x).
$$
We also know that $g(z)\in F(z)$, see \rf{GZ-ZP}, so that $a\in F(z)+\gamma\,\rho(x,z)\,\BXR$.
\par In the same way we show  that
$a\in F(z')+\gamma\,\rho(x,z')\,\BXR$ proving property
\rf{FFP}.
\smsk
\par Prove that
\bel{FXY}
F^{[1]}(x)+\gamma\,\rho(x,y)\,\BXR\supset F^{[1]}(y)~~~~
\text{for every}~~~~x,y\in\Mc.
\ee
\par We know that $F^{[1]}(x)\ne\emp$ which enables us to apply part (b) of Lemma \reff{H-R} to the left hand side of
\rf{FXY}. This lemma and definition \rf{F1-PR1} tell us that
\begin{align}
F^{[1]}(x)+\gamma\,\rho(x,y)\,\BXR&=
\left\{\bigcap_{z\in \Mc}\,
\left[F(z)+\gamma\,\rho(x,z)\,\BXR\right]\right\}
+\gamma\,\rho(x,y)\,\BXR\nn\\
&=\bigcap_{z\in \Mc}\,
\left[F(z)+(\gamma\,\rho(x,z)+\gamma\,\rho(x,y))\,\BXR\right]
\nn
\end{align}
so that, thanks to the triangle inequality,
$$
F^{[1]}(x)+\gamma\,\rho(x,y)\,\BXR\supset
\bigcap_{z\in \Mc}\,
\left[F(z)+\gamma\,\rho(y,z)\,\BXR\right]=F^{[1]}(y)
$$
proving \rf{FXY}. By interchanging the roles of $x$ and $y$ we obtain also
$$
F^{[1]}(y)+\gamma\,\rho(x,y)\,\BXR\supset F^{[1]}(x).
$$
These two inclusions prove the required inequality \rf{TM-1}.
\smsk
\par The proof of Proposition \reff{X-1DIM} is complete.\bx
\bsk\bsk

\noindent {\bf 5.2 Rectangular hulls of plane convex sets.}

\indent
\par Let us fix some additional notation.  We let
$\RCT$ denote the family of all bounded closed rectangles in $\RT$ {\it with sides parallel to the coordinate axes} $Ox_1$ and $Ox_2$.
\par Let $Q_0=\BX$ be the unit ball of the Banach space $X=\LTI$, i.e., the square $Q_0=[-1,1]^2$. Given $a\in\RT$ and $r\ge 0$, we set $rQ_0=[-r,r]^2$ and $Q(a,r)=rQ_0+a$.

\begin{definition} \lbl{DF-RH}{\em Let $S$ be a non-empty bounded convex closed subset in $\RT$. We set
$$
\HR[S]=\cbg\{\Pi: \Pi\in\RCT, \Pi\supset S\},
$$
and refer to $\HR[S]$ as the {\it ``rectangular hull``} of the set $S$.}
\end{definition}
\msk
\par Note the following useful representation of the rectangular hull which easily follows from Definition \reff{DF-RH}:
\begin{align}\lbl{H-DP}
\HR[S]=(S+Ox_1)\cbg (S+Ox_2).
\end{align}
\par In the next section we need the following auxiliary
\begin{lemma}\label{2K-H} Let $K_1,K_2\in\Kc(\RT)$ be two convex compacts in $\RT$ with non-empty intersection. Let $\tau\ge 0$ and let $Q=[-\tau,\tau]^2$. Then
\begin{align}\lbl{K-12}
(K_1\cbg K_2)+ Q=(K_1+ Q)\cbg (K_2+ Q)\cbg \HR[(K_1\cbg K_2)+Q].
\end{align}
\end{lemma}
\par {\it Proof.} Obviously, the right hand side of \rf{K-12} contains its left hand side.
\par Let us prove the converse statement. Fix a point
\begin{align}\lbl{X-H}
a\in
(K_1+ Q)\cbg (K_2+ Q)\cbg \HR[K_1\cbg K_2+Q].
\end{align}
\par Our aim is to prove that $a\in (K_1\cbg K_2)+ Q$.
Clearly, this property holds if and only if
$Q(a,\tau)\cbg K_1\cbg K_2\ne\emp$. It is also clear that
$Q(a,\tau)=\Pi_1(a)\cbg \Pi_2(a)$ where
$$
\Pi_i(a)=Q(a,\tau)+Ox_i,~~~~~~~i=1,2.
$$
\par Thus, $a\in (K_1\cbg K_2)+ Q$ provided
$$
K_1\cbg K_2\cbg \Pi_1(a)\cbg \Pi_2(a)\ne\emp.
$$
See Fig. 31.
\begin{figure}[h!]
\centering{\includegraphics[scale=0.16]{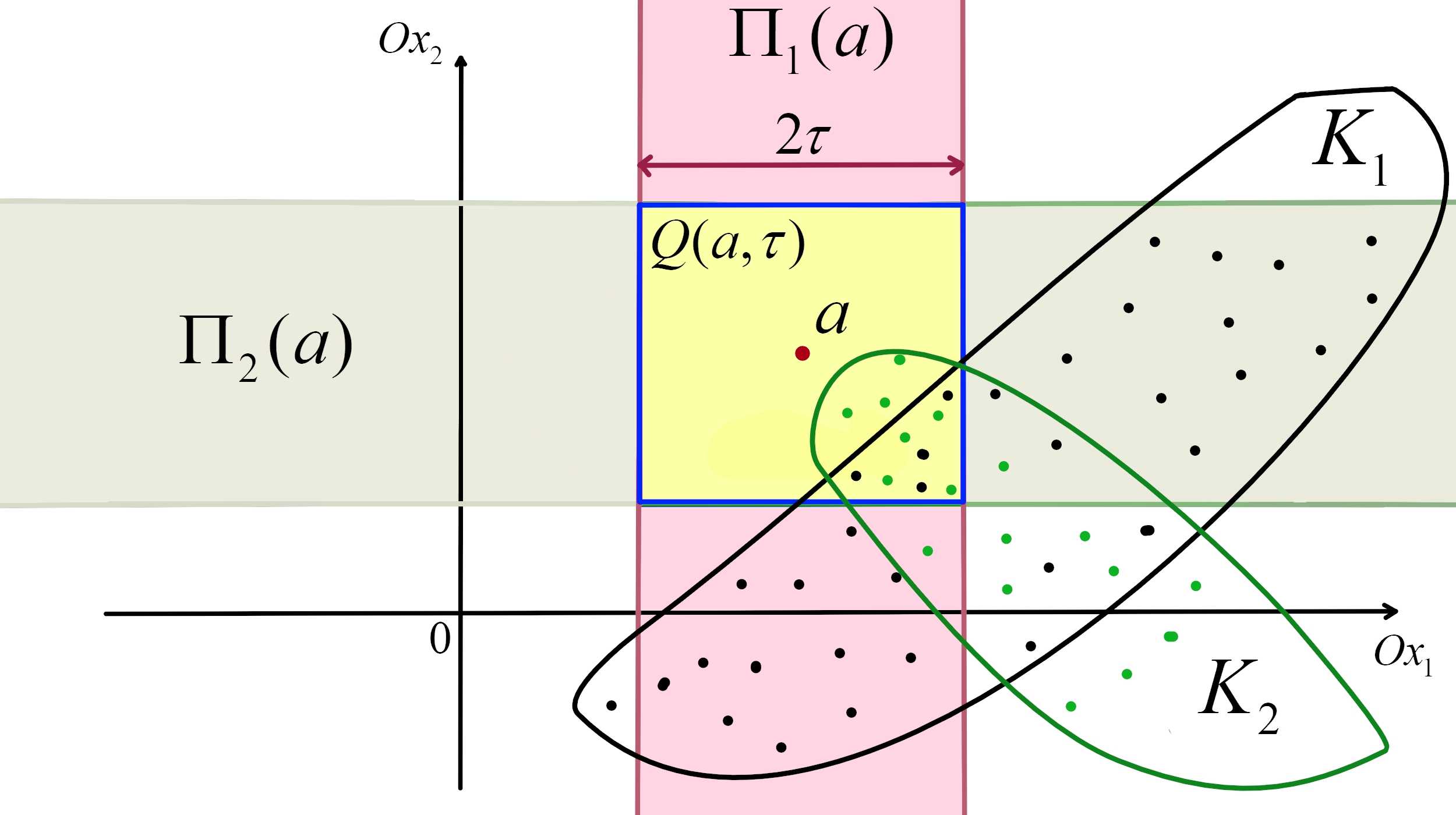}}
\caption{TEXT}
\end{figure}
\bsk

\par Thanks to Theorem \reff{H-TH}, the family of sets $\{K_1,K_2,\Pi_1(a),\Pi_2(a)\}$ has a common point provided  any three members of this family have a non-empty intersection.
\par Prove that it is true for $a$ satisfying \rf{X-H}. Indeed, $a\in K_i+Q$, so that $K_i\cbg Q(a,\tau)\ne\emp$, $i=1,2$. Hence,
$$
K_i\cbg \Pi_1(a)\cbg \Pi_2(a)=K_i\cbg Q(a,\tau)\ne\emp, ~~~~i=1,2.
$$
\par Next, thanks to \rf{H-DP} and \rf{X-H}, for every $i=1,2$,
$$
a\in\HR[(K_1\cbg K_2)+Q]\subset (K_1\cbg K_2)+Q+Ox_i.
$$
Hence, $K_1\cbg K_2\cbg \Pi_i(x)\ne\emp$, $i=1,2$, and the proof of the lemma is complete.\bx
\msk

\bsk\bsk

\noindent {\bf 5.3 Balanced refinements of set-valued mappings in $\LTI$.}

\indent
\par In this section we refine the result of Theorem \reff{MAIN-RT} for the space  $\X=\LTI$.
\begin{theorem}\label{LTI-M} In the settings of Theorem \reff{MAIN-RT} properties \rf{F2-NEM} and \rf{HD-RT} hold
provided $X=\LTI$\,,
\begin{align}\lbl{PR-LTI}
\lambda_1\ge 1,~~~~~\lambda_2\ge 3\lambda_1,~~~~\text{and}~~~~~\gamma\ge \lambda_2\,(3\lambda_2+\lambda_1)/
(\lambda_2-\lambda_1).
\end{align}
\par In particular, \rf{F2-NEM} and \rf{HD-RT} hold whenever $\lambda_1=1$, $\lambda_2=3$ and $\gamma=15$.
\end{theorem}
\par {\it Proof.} We mainly follow the scheme of the proof
of Theorem \reff{MAIN-RT} given in Section 3.
\par Let $F:\Mc\to\Kc(\RT)$ be a set-valued mapping satisfying the hypothesis of Theorem \reff{LTI-M}. Thus, the next statement is true.
\begin{claim}\label{LTI-A} For every $\Mc'\subset\Mc$, $\#\Mc'\le 4$, the restriction $F|_{\Mc'}$ of $F$ to $\Mc'$ has a $\rho$-Lipschitz selection $f_{\Mc'}:\Mc'\to \LTI$ with $\rho$-Lipschitz seminorm $\|f_{\Mc'}\|_{\Lip((\Mc',\,\rho),\LTI)}\le 1$.
\end{claim}

\par Let $\lambda_1$ and $\lambda_2$ be positive constants satisfying inequalities \rf{PR-LTI}. We set $L=\lambda_2/\lambda_1$; thus, $L\ge 3$. Then we introduce a pseudometric on $\Mc$ defined by
$\ds(x,y)=\lambda_1\,\rho(x,y)$, $x,y\in\Mc$.
\par We let $F^{[1]}$ and  $F^{[2]}$  denote the first and the second order $(\{1,L\},\ds)$-balanced refinements of $F$ respectively. See Definition \reff{F-IT}. Thus,
$$
F^{[1]}(x)=
\bigcap_{z\in\Mc}\,
\left[F(z)+\ds(x,z)\,Q_0\right]~~~\text{and}~~~
F^{[2]}(x)=\bigcap_{z\in\Mc}\,
\left[F^{[1]}(z)+L\ds(x,z)\,Q_0\right],~~x\in\Mc.
$$
\par We also recall that $e(\mfM,\LTI)=1$. In this case, Lemma \reff{G-NE1} and Proposition \reff{N-EM} tell us that $F^{[1]}(x)\ne\emp$ and $F^{[2]}(x)\ne\emp$ for every $x\in\Mc$.
\par Let
$$
\tgm(L)=L\,\theta(L)~~~\text{where}~~~ \theta(L)=(3L+1)/(L-1).
$$
\par Prove that
\begin{align}\lbl{H-LTI}
\dhf(F^{[2]}(x),F^{[2]}(y))\le \tgm(L)\,\ds(x,y)~~~
\text{for every}~~~x,y\in\Mc.
\end{align}
\par We recall that, thanks to formula \rf{G-XP},
$F^{[2]}(x)=\cbg\{T_x(u,u',u''):u,u',u''\in\Mc\}$
where, given $u,u',u''\in\Mc$,
\begin{align}\lbl{H-D1}
T_x(u,u',u'')=
\{[F(u')+\ds(u',u)Q_0]\cbg [F(u'')+\ds(u'',u)Q_0]\}+L\ds(u,x)\,Q_0.
\end{align}
\par In particular, $T_x(u,u',u'')\ne \emp$ for all $u,u',u''\in\Mc$ (because $F^{[2]}(x)\ne\emp$).
\par The next lemma is a refinement of the formula \rf{GX-1} for the special case of $X=\LTI$.
\begin{lemma}\label{F2-NB} Let $\tau>0$ and let $Q=\tau Q_0=[-\tau,\tau]^2$. Then, for every $x\in\Mc$, we have
\begin{align}\lbl{F-KNT}
F^{[2]}(x)+Q=
\bigcap_{v,u,u',u''\in\Mc}\,\,
\left\{\,
\left[T_x(u,u',u'')\cbg\left(F(v)+\ds(x,v)Q_0\right)\right]+Q
\right\}.
\end{align}
\end{lemma}
\par {\it Proof.} Lemma \reff{L-AW} tells us that
\begin{align}\lbl{F2-P}
F^{[2]}(x)+Q=
\bigcap\,\,
\left\{\,\left[T_x(u_1,u'_1,u''_1)\cbg T_x(u_2,u'_2,u''_2)\right]+Q\,\right\}
\end{align}
where the intersection is taken over all $u_i,u'_i,u''_i\in\Mc$, $i=1,2$. Note also that, by \rf{H-D1},
$$
F(v)+\ds(x,v)Q_0=T_x(x,v,v).
$$
From this and \rf{F2-P} it follows that the right hand side of \rf{F-KNT} contains its left hand side.
\par Prove the converse statement. Fix a point
\begin{align}\lbl{Y-DQ}
a\in
\bigcap_{v,u,u',u''\in\Mc}\,\,
\left\{\,
\left[T_x(u,u',u'')\cbg\left(F(v)+
\ds(x,v)Q_0\right)\right]+Q
\right\}
\end{align}
and show that $a\in F^{[2]}(x)+Q$. In view of formula \rf{F2-P}, it suffice to prove that for every $u_1,u'_1,u''_1,u_2,u'_2,u''_2\in\Mc$ {\it the point $a$ belongs to the set $A$ defined by}
\begin{align}\lbl{A-TT}
A=\left[T_x(u_1,u'_1,u''_1)\cbg T_x(u_2,u'_2,u''_2)\right]+Q.
\end{align}
%
\par To see this, given $i\in\{1,2\}$ we introduce the following sets: $Q_i=L\ds(u_i,x)\,Q_0$,
\begin{align}\lbl{K-FQ}
K'_{i}=F(u'_i)+\ds(u_i,u'_i)Q_0~~~~
\text{and}~~~~K''_{i}=F(u''_i)+\ds(u_i,u''_i)Q_0.
\end{align}
\par In these settings,
$T_x(u_i,u'_i,u''_i)=K'_{i}\cbg K''_{i}+Q_i$, $i=1,2$. See \rf{H-D1}.
\par Note that $K'_{i}\cbg K''_{i}\ne\emp$ because $T_x(u_i,u'_i,u''_i)\ne\emp$. Therefore, thanks to Lemma \reff{2K-H},
\begin{align}\lbl{TT-HP}
T_x(u_i,u'_i,u''_i)=(K'_{i}+Q_i)\cbg (K''_{i}+Q_i)\cbg
\Hc[T_x(u_i,u'_i,u''_i)],~~~~~i=1,2.
\end{align}
\par Now, let us introduce the following families of sets:
$$
\Kc^+=\{K'_{i}+Q_i,K''_{i}+Q_i:i=1,2\},
~~\Kc^{++}=\{\Hc[T_x(u_i,u'_i,u''_i)]:i=1,2\},
~~\Kc=\Kc^+\cupbig \Kc^{++}.
$$
Then, thanks to \rf{A-TT} and \rf{TT-HP}, $A=[\cbg\{K:K\in\Kc\}]+Q$.
\par We recall that, thanks to Proposition \reff{N-EM}, the set $F^{[2]}(x)\ne\emp$ so that the left hand side of \rf{F2-P} is non-empty as well. From this and \rf{A-TT} it follows that $A\ne\emp$ proving that $\cbg\{K:K\in\Kc\}\ne\emp$. Therefore, thanks to Lemma \reff{H-IN},
$$
A=\cbg\{[K\cbg K']+Q:K,K'\in\Kc\}.
$$
\par Thus, to prove that $a\in A$ it suffice to show that $a\in K\cbg K'+Q$ for every $K,K'\in\Kc$.
\par To do this, first let us note that, thanks to \rf{K-FQ},
$$
K'_{i}+Q_i=F(u'_i)+\ds(u_i,u'_i)Q_0+L\ds(u_i,x)\,Q_0
\supset F(u'_i)+(\ds(u_i,u'_i)+\ds(u_i,x))\,Q_0
$$
for every $i=1,2$. Therefore, thanks to the triangle inequality,
\begin{align}\lbl{M-1}
K'_{i}+Q_i\supset F(u'_i)+\ds(u'_i,x)\,Q_0.
\end{align}
In the same way we prove that
\begin{align}\lbl{M-2}
K''_{i}+Q_i\supset F(u''_i)+\ds(u''_i,x)\,Q_0,~~~~i=1,2.
\end{align}
\par Furthermore, we know that
\begin{align}\lbl{M-3}
\Hc[T_x(u_i,u'_i,u''_i)]\supset T_x(u_i,u'_i,u''_i), ~~~~i=1,2.
\end{align}
\par On the other hand, property \rf{Y-DQ} tells us that
\begin{align}\lbl{A-FP}
a\in
T_x(u,u',u'')\cbg[F(v)+\ds(x,v)Q_0)]+Q~~~~
\text{for every}~~~~u,u',u'',v\in\Mc.
\end{align}
\par Combining this property with \rf{M-1}, \rf{M-2} and \rf{M-3}, we conclude that
$$
a\in K\cbg K'+Q~~~\text{whenever either}~~~
K\in\Kc^+,K'\in\Kc^{++}~~~\text{or}~~~K,K'\in\Kc^+.
$$
\par It remains to prove that
\begin{align}\lbl{ABC}
a\in H_1\cbg H_2+Q~~~\text{where}
~~~H_i=\Hc[T_x(u_i,u'_i,u''_i)], ~~~~~i=1,2.
\end{align}
\par It is immediate from Lemma \reff{H-R}, part (b),
that
$$
H_1\cbg H_2+Q=(H_1+Q)\cbg (H_2+Q)
$$
so that
$$
H_1\cbg H_2+Q=\left\{\Hc[T_x(u_1,u'_1,u''_1)]+Q\right\}\cbg \left\{\Hc[T_x(u_2,u'_2,u''_2)]+Q\right\}.
$$
From this and \rf{M-3}, we have
\begin{align}\lbl{FN-Y}
H_1\cbg H_2+Q\supset \left\{T_x(u_1,u'_1,u''_1)+Q\right\}\cbg \left\{T_x(u_2,u'_2,u''_2)+Q\right\}.
\end{align}
\par But, thanks to \rf{A-FP}, $a\in T_x(u_i,u'_i,u''_i)+Q$, $i=1,2$. Combining this property
with \rf{FN-Y}, we obtain the required property \rf{ABC}
completing the proof of the lemma.\bx
\bsk
\par We are in a position to prove inequality \rf{H-LTI}.
Our proof will follow the scheme of the proof of Proposition \reff{HD-G1}.
\smsk
\par Let $x,y\in\Mc$, and let $\tau=\tgm(L)\ds(x,y)$. (Recall that $\tgm(L)=L\,\theta(L)$ and $\ds=\lambda_1\,\rho$.)
\par Lemma \reff{F2-NB} tells us that
\begin{align}\lbl{F-LRP}
F^{[2]}(x)+\tau Q_0=
\bigcap_{v,u,u',u''\in\Mc}\,\,
\left\{\,
\left[T_x(u,u',u'')\cbg\left(F(v)+\ds(x,v)Q_0\right)\right]
+\tau Q_0\right\}.
\end{align}
\par Let us fix elements $u,u',u'',v\in\Mc$ and a set
\begin{align}\lbl{A-LRP1}
\tA=\left[T_x(u,u',u'')\cbg(F(v)+\ds(x,v)Q_0)\right]
+\tau Q_0.
\end{align}
\par Our goal is to show that
$$
\tA\supset F^{[2]}(y)= \bigcap_{z,z',z''\in\Mc}\,T_y(z,z',z'').~~~~~~~~\text{See \rf{G-XP}.}
$$
\par In fact, we will prove a stronger imbedding:
$$
\tA\,\supset\,\, T_y(u,u',u'')\,\bigcap\, T_y(x,u',v)\,\bigcap\, T_y(x,u'',v).
$$
\par We also note that the set $\tA$ is the orbit of the element $y$ with respect to the diagram shown in Fig. 32.

\begin{figure}[h!]
\centering{\includegraphics[scale=0.67]{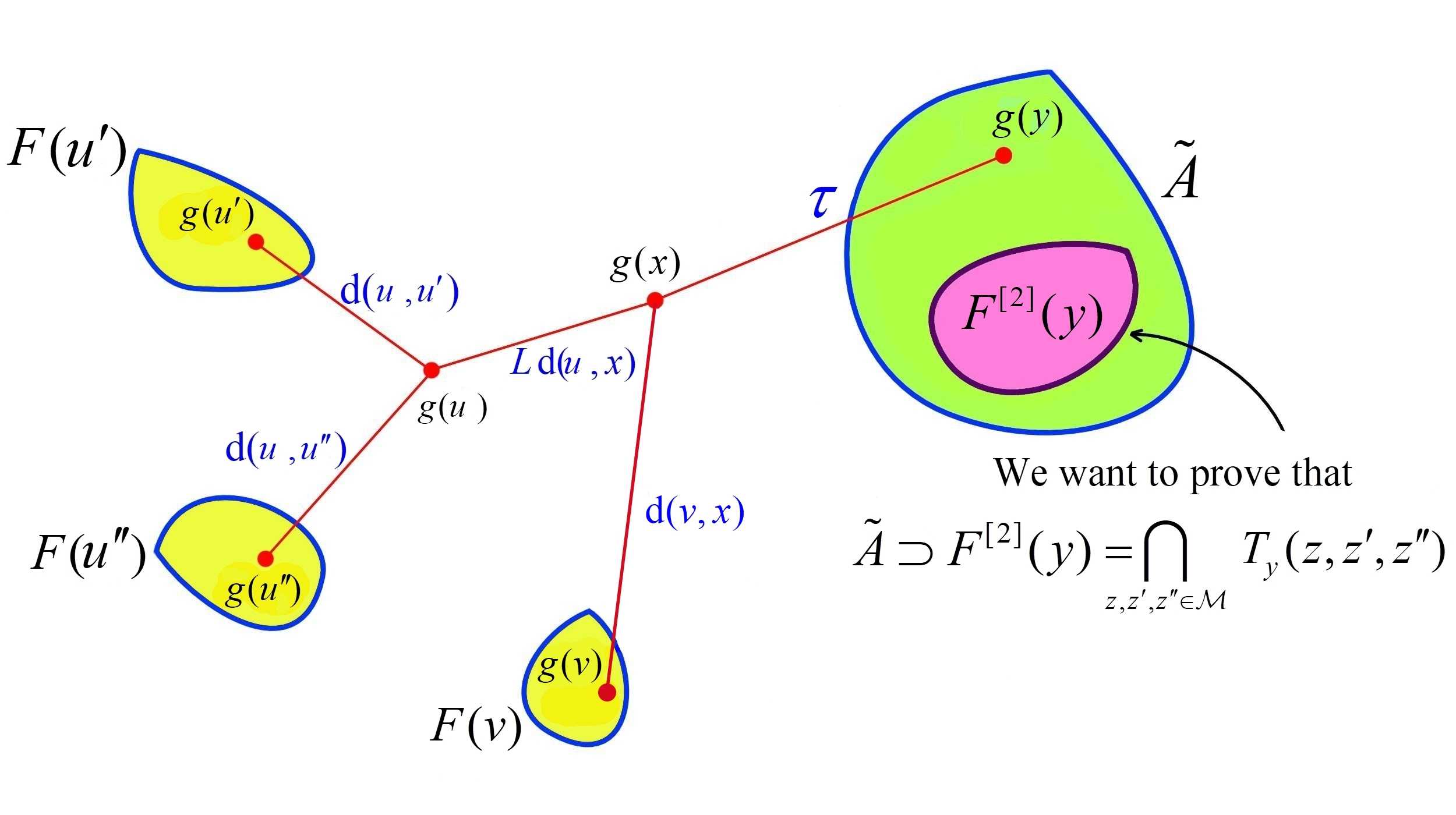}}
\caption{$\tA=\{g(y)\}$ where $g$ runs over all mappings  which agree with this diagram.}
\end{figure}
\par Let
\begin{align}\lbl{C123-2}
C_1=F(u')+\ds(u',u)Q_0,~~~C_2=F(u'')+\ds(u'',u)Q_0,~~~
C=F(v)+\ds(x,v)Q_0,
\end{align}
and let $\ve=L\,\ds(x,y)$ and $r=\ds(u,x)$.
Then
$$
\tau=\tgm(L)\ds(x,y)=L\,\theta(L)\ds(x,y)
=\theta(L)\,\ve.
$$
In these settings, $T_x(u,u',u'')=(C_1\cbg C_2)+L\,r\, Q_0$, see \rf{H-D1}, and
$$
\tA=[\{(C_1\cbg C_2)+L\,r\, Q_0\}\cbg C] +\theta(L)\,\ve\,Q_0.~~~~~~~~~
\text{See \rf{A-LRP1}.}
$$

%
%
\par Let us apply Proposition \reff{P-F3} to the sets $C_1,C_2$ and $C$ defined by \rf{C123-2}. To do this, we have to verify condition \rf{A-PT} of that proposition, i.e., to prove that
\begin{align}\lbl{PT-2}
C_1\cbg C_2\cbg(C+rQ_0)\ne\emp.
\end{align}
\par Let $\Mc'=\{u,u',v\}$. Then, thanks to Claim \reff{LTI-A}, there exists a $\rho$-Lipschitz selection $f_{\Mc'}:\Mc'\to\LTI$ of the restriction $F|_{\Mc'}$ with $\|f_{\Mc'}\|_{\Lip((\Mc',\,\rho),\LTI)}\le 1$.
\par Because $e(\mfM,\LTI)=1$ and $\ds=\lambda_1\rho\ge\rho$, the mapping $f_{\Mc'}:\Mc'\to\LTI$
can be extended to a $\ds$-Lipschitz mapping $\tf:\Mc\to\LTI$ defined on all of $\Mc$ with $\ds$-Lipschitz seminorm
$$
\|\tf\|_{\Lip((\Mc,\ds),\LTI)}\le
\|f_{\Mc'}\|_{\Lip((\Mc',\,\rho),\LTI)}\le 1.
$$
\par Thus, $\tf(u')=f_{\Mc'}(u')\in F(u')$, $\tf(u'')=f_{\Mc'}(u'')\in F(u'')$, $\tf(v)=f_{\Mc'}(v)\in F(v)$,
$$
\|\tf(u')-\tf(u)\|\le \ds(u',u),~~~
\|\tf(u'')-\tf(u)\|\le \ds(u'',u)
$$
and
$$
\|\tf(x)-\tf(u)\|\le \ds(u,x)=r,~~~~~
\|\tf(x)-\tf(v)\|\le \ds(v,x).
$$
Hence, $\tf(u)\in C_1\cbg C_2$ and $\tf(x)\in C$, so that
$C_1\cbg C_2\cbg(C+rQ_0)\ni\tf(u)$ proving \rf{PT-2}.
\smsk
\par This enables us to apply Proposition \reff{P-F3} to the sets $C_1$, $C_2$ and $C$. By this proposition,
\begin{align}
\tA&=[\{(C_1\cbg C_2)+LrQ_0\}\cbg C]+\theta(L)\,\ve\,Q_0
\nn\\
&\supset [(C_1\cbg C_2)+(Lr+\ve)Q_0]
\cbg [\{(C_1+rQ_0)\cbg C\} +\ve Q_0]
\cbg [\{(C_2+rQ_0)\cbg C\} +\ve Q_0]
\nn\\
&=S_1\cbg S_2\cbg S_3.
\nn
\end{align}
\par Prove that $S_i\supset F^{[2]}(y)$ for every $i=1,2,3$. We begin with the set
\begin{align}
S_1&=(C_1\cbg C_2)+(Lr+\ve)Q_0\nn\\
&=
[\{F(u')+\ds(u',u)Q_0\}\cbg \{F(u'')+\ds(u'',u)Q_0\}]
+(L\ds(u,x)+L\ds(x,y))Q_0\,.
\nn
\end{align}
%
See \rf{C123-2}. By the triangle inequality, $\ds(u,x)+\ds(x,y)\ge \ds(u,y)$ so that
$$
S_1\supset
[\{F(u')+\ds(u',u)Q_0\}\cbg \{F(u'')+\ds(u'',u)Q_0\}]
+L\ds(u,y)Q_0=T_y(u,u',u'')\,.
$$
See \rf{H-D1}. But, thanks to \rf{G-XP}, $T_y(u,u',u'')\supset F^{[2]}(y)$ proving the required inclusion $S_1\supset F^{[2]}(y)$.
\smsk
\par Prove that $S_2\supset F^{[2]}(y)$. We have
\begin{align}
S_2&=
[(C_1+rQ_0)\cbg C] +\ve Q_0\nn\\
&=
[\{(F(u')+\ds(u',u)Q_0)+\ds(x,u)Q_0\}\cbg \{F(v)+\ds(x,v)Q_0\}]+L\ds(x,y)Q_0\,.
\nn
\end{align}
%
Therefore, thanks to the triangle inequality, \rf{H-D1} and \rf{G-XP}
$$
S_2\supset
[(F(u')+\ds(u',x)Q_0)\cbg (F(v)+\ds(x,v)Q_0)]
+L\ds(x,y)Q_0=T_y(x,u',v)\supset F^{[2]}(y).
$$
\par In the same way we show that $S_3\supset F^{[2]}(y)$.
Hence, $\tA\supset S_1\cbg S_2\cbg S_3\supset F^{[2]}(y)$.
\smsk
\par Combining this inclusion with definition \rf{A-LRP1} and representation \rf{F-LRP}, we conclude that
$F^{[2]}(x)+\tau Q_0\supset F^{[2]}(y)$. By interchanging the roles of $x$ and $y$ we obtain also the inclusion $F^{[2]}(y)+\tau Q_0\supset F^{[2]}(x)$. These two inclusions imply inequality
$$
\dhf(F^{[2]}(x),F^{[2]}(y))\le \tau=\tgm(L)\,\ds(x,y)=\lambda_1\,\tgm(L)\,\rho(x,y)
$$
proving \rf{HD-RT} with $\gamma=\lambda_1\,L(3L+1)/(L-1)$.
We recall that $L=\lambda_2/\lambda_1$, so that inequality \rf{HD-RT} holds with any $\gamma\ge \lambda_2\,(3\lambda_2+\lambda_1)/
(\lambda_2-\lambda_1)$.
\par The proof of Theorem \reff{LTI-M} is complete.\bx
\msk


\end{document}